\documentclass[a4paper,12pt]{amsart} 
\title[Chern--Weil and Hilbert--Samuel formulae]{Chern--Weil and Hilbert--Samuel formulae for singular hermitian line bundles}

\usepackage{amsmath, amssymb, mathrsfs, amsthm, geometry}
\usepackage{mathtools} 
\usepackage{hyperref}
\usepackage{float}
\usepackage[all]{xy}

\usepackage{enumitem}

\usepackage{cleveref}

\crefformat{equation}{(#2#1#3)} 

\AtBeginDocument{\renewcommand{\ref}[1]{\Cref{#1}}}

\usepackage{tikz-cd} 
\usetikzlibrary{arrows,matrix,positioning}
\tikzstyle{dmatrix}=[matrix of math nodes,row sep=2.5em, column sep=2.5em,
text height=1.5ex, text depth=0.25ex] 
\usepackage{color}

\newcommand{\on}[1]{\operatorname{#1}}
\newcommand{\bb}[1]{{\mathbb{#1}}} 

\newcommand{\ca}[1]{{\mathcal{#1}}}

\newcommand{\abs}[1]{\lvert#1\rvert}

\newcommand{\hra}{\hookrightarrow}
\newcommand{\sub}{\subseteq}

\theoremstyle{definition}
\newtheorem{definition}{Definition}[section]

\theoremstyle{plain}
\newtheorem{proposition}[definition]{Proposition}
\newtheorem{lemma}[definition]{Lemma}
\newtheorem{theorem}[definition]{Theorem}
\newtheorem{corollary}[definition]{Corollary}

\newtheorem{introtheorem}{Theorem}

\theoremstyle{remark}

\newtheorem{remark}[definition]{Remark}

\newtheorem{example}[definition]{Example}

\newcounter{temp}

\date{\today}

\newcounter{nootje}
\newcommand{\beq}{\begin{equation}}
\newcommand{\eeq}{\end{equation}}
\newcommand{\beqs}{\begin{equation*}}
\newcommand{\eeqs}{\end{equation*}}

\newcommand{\SL}{\operatorname{SL}}
\newcommand{\divisor}{\operatorname{div}}
\newcommand{\vol}{\operatorname{vol}}
\newcommand{\id}{\operatorname{Id}}
\renewcommand{\Im}{\operatorname{Im}}

\def\rr{\mathbb{R}}

\def\zz{\mathbb{Z}}
\def\cc{\mathbb{C}}
\def\qq{\mathbb{Q}}

\def\L{\mathcal{L}}
\def\M{\mathcal{M}}
\def\Div{\operatorname{Div}_{\mathbb{R}}}
\def\WbDiv{\operatorname{W-b-Div}_{\rr}}
\def\CbDiv{\operatorname{C-b-Div}_{\rr}}
\newcommand*\isom{\xrightarrow{\sim}}

\newcommand{\tmatrix}{A}
\newcommand{\tvector}{c}

\newcommand{\class}[1]{\on{cl}(#1)}

\numberwithin{equation}{section}
\subjclass{14C17, 14E99, 32U05, 32U25, 52A39}
 \thanks{A. Botero was supported by the collaborative research 
center SFB 1085 \emph{Higher Invariants - Interactions between
  Arithmetic Geometry and Global Analysis} funded by the Deutsche
Forschungsgemeinschaft. J.~I.~Burgos was partially supported by MINISTERIO
DE CIENCIA E INNOVACION research projects PID2019-108936GB-C21 and 
ICMAT Severo Ochoa project CEX2019-000904-S.
D. Holmes was supported by NWO grant VI.Vidi.193.006 and NWO grant
613.009.103.
}

\begin{document}
\author[A.~M.~Botero]{Ana Mar\'ia Botero}
\author[J.~I.~Burgos Gil]{Jos\'e Ignacio Burgos Gil}
\author[D.~Holmes]{David Holmes}
\author[R.~de~Jong]{Robin de Jong}

\begin{abstract}
We show a Chern--Weil type statement and a Hilbert--Samuel formula for
a large class of singular plurisubharmonic metrics on a line bundle
over a smooth projective complex variety. For this we use the theory
of b-divisors and the so-called multiplier ideal volume function. We
apply our results to the line bundle of Siegel--Jacobi forms over the
universal abelian variety endowed with its canonical invariant
metric. This generalizes the results of \cite{bkk} to higher degrees.
\end{abstract}

\maketitle 

\tableofcontents

\section{Introduction}

\subsection{Motivation and background}

Let $X$ be a differentiable manifold and let $\ca E$ be
a vector bundle on $X$. The \emph{Chern classes} of $\ca E$ are
topological invariants $c_i(\ca E) \in H^{2i}(X, \zz)$ that measure,
in some sense, how far $\ca E$ is from being trivial.
\emph{Chern-Weil theory} tells us that these Chern classes can
be represented in the de Rham cohomology of $X$ by means of smooth differential forms
obtained by applying an invariant polynomial to the curvature matrix
associated to a smooth connection on $\ca E$. 

When $X$ is a complex manifold and $\ca E$ is a holomorphic vector bundle on $X$
equipped with a smooth hermitian metric, there is a unique smooth connection on $\ca E$ compatible at the same time
with the hermitian metric and with the holomorphic
structure. We see that for a holomorphic vector bundle $\ca E$ on a complex
manifold, a smooth hermitian metric $h$ uniquely determines \emph{Chern
  forms} $c_{1}(\ca E,h),\dots,c_{i}(\ca E,h),\dots$ representing the Chern
classes $c_{1}(\ca E),\dots,c_{i}(\ca E),\dots$ in de Rham cohomology. 

Assume that $X$ is compact, that $\dim X = n$, and let $x_1,\ldots,x_n$ be variables where $x_i$ has weight $i$.
Let $P \in \bb{Z}[\ldots,x_i,\ldots]$ be a homogeneous polynomial of weight $n$. As a particular instance of Chern--Weil theory we find the equality
\begin{equation}
 P(\ldots, c_i(\ca E),\ldots) = \int_X P(\ldots, c_i(\ca E,h), \ldots )
\end{equation}
in $\zz$, where the expression $P(\ldots, c_i(\ca E),\ldots)$ on the left hand side is interpreted as an integer
upon taking the degree, and where $P(\ldots, c_i(\ca E,h), \ldots )$ is a closed $(n,n)$-form on $X$. As a special case, we note that if  $\ca L$ is a holomorphic \emph{line bundle} on $X$ equipped with
a smooth hermitian metric $h$ then we have the equality
 \begin{equation} \label{eq:line_bundle_smooth}
c_1(\ca L)^n = \int_X c_1(\ca L,h)^n 
\end{equation}
in $\zz$, with $c_1(\ca L)^n \in \zz$ the \emph{degree} of the line bundle $\ca L$.

Now although smooth hermitian metrics always exist, they can be difficult to write
down explicitly. In fact in many situations the \emph{natural} hermitian
metric $h$ associated to the problem at  hand is only smooth on a dense open subset but singular along say a normal crossings divisor $D$ on $X$. In this context Mumford \cite{hi} has introduced the notion of \emph{good} metrics. The condition of being ``good'' is phrased in terms of conditions on the asymptotic behavior of $h$ and its first and second derivatives near $D$; we refer to
\ref{exm:6} for the precise definition. 

Good metrics are a class of singular metrics that for the purpose of
Chern--Weil theory are as good as smooth metrics. More precisely,  let $X$ be a
compact complex manifold, $\ca E$ a holomorphic vector bundle on $X$ and $U\subset X$ a dense
open subset whose complement is a normal crossings divisor. Let $h$ be a smooth metric on $\ca E|_{U}$ such that $h$ is
good in the sense of Mumford as a singular metric on $\ca E$. Then for every polynomial $P\in
\bb{Z}[\dots,x_{i},\dots]$ the smooth closed differential form
\begin{displaymath}
  P(\ca E,h)=P(c_{1}(\ca E,h),\dots,c_{i}(\ca E,h),\dots)
\end{displaymath}
on $U$ is locally integrable as a differential form on $X$, and thus 
determines a closed current $[P(\ca E,h)]$ on $X$. In fact, the closed current $[P(\ca E,h)]$
represents the class
$P(\ca E)=P(\dots,c_{i}(\ca E),\dots)$ in de Rham cohomology of $X$. In particular, if $x_i$ has weight $i$ and 
$P \in \bb{Z}[\ldots,x_i,\ldots]$ is a homogeneous polynomial of weight $n$, we have the equality
\begin{equation} \label{eq:24}
 P(\ldots, c_i(\ca E),\ldots) = \int_U P(\ldots, c_i(\ca E,h), \ldots )
\end{equation}
in $\zz$, where again the expression $P(\ldots, c_i(\ca E),\ldots)$ on the left hand side is interpreted as an integer
upon taking the degree. If  $\ca L$ is a holomorphic \emph{line bundle} on $X$ equipped with
a smooth metric $h$ on $\ca L|_{U}$ such that $h$ is
good as a singular metric on $\ca L$, then we have the equality
 \begin{equation} \label{eq:line_bundle}
c_1(\ca L)^n = \int_U c_1(\ca L,h)^n 
\end{equation}
in $\zz$. 
 
As equations \ref{eq:24} and \ref{eq:line_bundle} suggest, good metrics are better
understood in the context of \emph{extending} a vector bundle from a
non-compact manifold to a compactification. In fact, the notion of good metric appears naturally in the study of
compactifications of \emph{locally symmetric domains}. 

 Let $D$ be a
bounded symmetric domain and $\Gamma $ a neat arithmetic group acting
on $D$ by isometries. Then the quotient $U=D/\Gamma $ is a locally symmetric space
and has a natural structure of a quasi-projective variety.  In \cite{AVRT:compact} the authors introduce
a family of toroidal compactifications of $U$. We have that $D$ is a quotient $D=K\backslash G$, where $G$ is a
semisimple adjoint group and $K$ is a maximal compact
subgroup. Any unitary representation of $K$ defines a (fully
decomposable) holomorphic vector bundle $\ca E^o$ on $U$ together with a smooth
hermitian metric $h$.

As Mumford proves in \cite{hi}, when $X$ is a smooth toroidal compactification of $U=D/\Gamma$ in the sense of \cite{AVRT:compact}, 
there exists a unique extension of $\ca E^o$ to a vector bundle $\ca E$ on $X$ such that the metric $h$ extends to a singular
  good hermitian metric. In particular, when $n=\dim X$ and the variables $x_i$ have weight $i$ then 
  for every homogeneous polynomial $P \in \bb{Z}[\ldots,x_i,\ldots]$ of weight $n$ we have the
  equality \ref{eq:24} in $\zz$. As a consequence, we see that the element
  $P(\ldots, c_i(\ca E),\ldots) \in \zz$ does not depend on the choice of toroidal compactification $X$.

Mumford's result applies notably to the case of \emph{pure Shimura
varieties}, as the connected components of their associated
analytic spaces are locally symmetric spaces. 

It makes sense to ask whether an extension of the Chern--Weil result \ref{eq:24} holds in the context
of \emph{mixed Shimura varieties}, for example moduli
spaces of abelian varieties with marked points. Such spaces again come equipped with natural 
automorphic line bundles, with natural smooth hermitian metrics on them. Moreover, such spaces have natural
smooth toroidal compactifications.

However, as was shown in \cite{bkk} in an example dealing with certain modular elliptic surfaces, in the 
case of mixed Shimura varieties new types
of metric singularities appear, that are \emph{not} good in the sense of Mumford. In fact,
a naive extension of the Chern--Weil result \ref{eq:24} turns out to fail already in this simple situation.

Let us give some details.
Let $U$ denote the moduli space $E(N)$ of elliptic curves with two marked points and with level~$N$. Let $\ca L$
be the line bundle of Jacobi forms on $U$ of weight four and index four, and let $h$ denote the natural invariant hermitian metric on $\ca L$.
As is proved in \cite{bkk}, for any smooth toroidal compactification $X$ of $U$, there exists a
subset $S \subset X$ of codimension two as well as an integer $r\ge 1$ such that
the smooth hermitian line bundle $\ca L^{\otimes r} $ over $U$ can be
extended to a line bundle $\ca L_{X\setminus S,r}$ over 
$  X\setminus S$ in such a way that $h^{\otimes r}$ extends to a  good singular
metric in the sense of Mumford. Since $S$ has codimension two, the
line bundle $\ca L_{X\setminus S,r}$ can be
extended uniquely to a line bundle $\ca L_{ X,r}$ over
$ X$.  However, the metric $h$ is no longer 
good at some of the points of $S$ and in fact, the intersection number
\begin{displaymath}
  \frac{1}{r^{2}} c_{1}(\ca L_{ X,r})^{2}
\end{displaymath}
turns out to depend on the choice of toroidal compactification $
X$. It follows that \ref{eq:24} can not be directly
extended to the case of mixed Shimura varieties. 

As is shown in \cite{bkk}, to recover a Chern--Weil type result in
the setting of the mixed Shimura variety $U=E(N)$  we should not
extend the line bundle $\ca L$ of Jacobi forms as a line bundle, but as a different object. The result
is best stated in the language of \emph{b-divisors}. Loosely speaking, a b-divisor is a collection of divisors
on a tower of modifications of a given compactification $X$
of $U$, compatible with push-forward. The precise definition is given
in \ref{sec:b-div}.

Let $s$ be a non-zero rational
section of the holomorphic line bundle $\ca L$, and let  $D$ be the
divisor of $s$, so that  $\ca L\simeq \ca{O}(D)$. In \cite{bkk},  a
natural
b-divisor $  D(\ca L,h,s)$ on $X$ with rational coefficients extending the divisor $D$ is constructed, which  is then shown
to have a natural self-intersection number $ D(\ca L,h,s)^2$ in $\bb R$. Here, the self-intersection number $ D(\ca L,h,s)^2$ is computed as a limit of
self-intersections of the ``incarnations'' of $ D(\ca L,h,s)$ on all models in the tower.
Finally, the equality 
\begin{equation}\label{eq:BKK}
   D(\ca L,h,s)^{2}=\int _{U} c_{1}(\ca L,h)^{\wedge 2}
\end{equation}
is shown to hold in $\bb R$. Note that \ref{eq:BKK} can be seen as a
Chern--Weil type result, where the intersection number on the left is
taken in the sense of b-divisors on $X$, and the integral on the right
is taken over the open dense subset $U \subset X$.

In \cite{bkk} the above Chern--Weil result is complemented with a
result of Hilbert--Samuel type that we also recall. Let $\bb H$ be the
upper half plane. Jacobi forms of weight $k$ and index
$m$ are holomorphic
functions on $\bb H\times \bb C$ satisfying certain transformation
properties with respect to a subgroup $\Gamma \subset \SL(2,\bb Z)$
and the abelian group $\bb Z^{2}$ as well as a growth condition when
approaching the boundary of $\bb H$. The basic reference for the theory of Jacobi
forms is the book \cite{EZ}. Let $J_{4r,4r}(\Gamma (N))$ be the space
of Jacobi forms of weight $4r$ and index $4r$ with respect to the
principal congruence subgroup $\Gamma (N)$. The elements of
$J_{4r,4r}(\Gamma (N))$ can be seen as sections  of the line bundle $\ca
L^{\otimes r}$ over $U$ satisfying some growth conditions along $X\setminus
U$. The second main result of \cite{bkk} is the formula
\begin{equation}\label{eq:26}
  \lim_{r\to \infty} \frac{\dim J_{4r,4r}(\Gamma (N))}{r^{2}/2}= D(\ca L,h,s)^{2}.
\end{equation}
That is, the asymptotic growth of the dimensions of the spaces of Jacobi
forms is not given by the self-intersection of a divisor but of a
b-divisor. Therefore to recover a Hilbert--Samuel type of result for
the space of Jacobi forms  we need to extend $\ca L$ as a b-divisor
and not as an ordinary line bundle.

Note that by combining equations \ref{eq:BKK} and \ref{eq:26}  we
find that the asymptotic growth of the space of Jacobi forms is
governed by an integral of a smooth differential form over~$U$.

The proof of \ref{eq:BKK}  and \ref{eq:26} in \cite{bkk} consists of
an explicit computation of both sides of each equation. The  present
paper finds its origin in the following two questions: whether
there could be a more intrinsic approach to  proving the equalities
\ref{eq:BKK} and \ref{eq:26}, and whether these equalities  might
admit a generalization to the setting of automorphic line bundles on
mixed Shimura varieties in \emph{higher dimensions}. 
We will answer both questions in the affirmative. 
As we will see in a forthcoming paper,  in particular this will give us a formula for the asymptotic growth of the dimensions of
spaces of Siegel--Jacobi forms.

\subsection{Statement of the main results}

Let $X$ be a
projective complex manifold and let $\ca L$ be a holomorphic line
bundle on $X$. We start out by defining a suitable class of singular
hermitian metrics on $\ca L$, namely those with \emph{almost
  asymptotically algebraic singularities}. 

The precise definition of almost asymptotically algebraic
singularities is given in \ref{def:6}, but in essence it asserts that
our metrics be \emph{plurisubharmonic}, with local potentials that are
bounded away from a normal crossings divisor $D$ on $X$ and
that can be approximated well by logarithms of sums of square norms of
holomorphic functions near~$D$ (i.e., by metrics with \emph{algebraic
  singularities}).

Many natural singular hermitian metrics turn out to be almost
asymptotically algebraic. For example, we prove, using Demailly's
celebrated regularization theorem for plurisubharmonic metrics (see
\cite[Theorem 4.2]{dk}, \cite[Theorem 3.2]{DP}), that a \emph{toroidal
  plurisubharmonic metric} is almost asymptotically algebraic; see
\ref{prop:toroidal_implies_aaas}. Also, as we will see in
\ref{ex:good_implies_aaas}, a plurisubharmonic metric which is good in
the sense of Mumford has almost asymptotically algebraic
singularities.

Our main results are generalizations of the Chern--Weil formula
\ref{eq:BKK} and of the Hilbert--Samuel formula \ref{eq:26}
to the setting of plurisubharmonic metrics with almost asymptotically
algebraic singularities. 
In order to formulate our Chern--Weil result, let $h$ be a
plurisubharmonic metric on $\ca L$. 

We can make sense of $c_1(\ca L, h)$ as a closed, positive
$(1,1)$-current on $X$ simply by taking $-dd^c\log h(s)$ for a local
generating section $s$ of $\ca L$. Let $n= \dim X$. A generalisation
of Bedford-Taylor calculus by Boucksom, Eyssidieux, Guedj and Zeriahi
\cite{begz} allows us to define the non-pluripolar Monge--Amp\`ere
measure $\left\langle c_1(\ca L, h)^n \right\rangle$, a closed,
positive $(n, n)$-current on $X$.

The right hand side of our Chern--Weil formula will be the \emph{non-pluripolar volume}
\begin{equation} \label{RHS}
\int_X \left\langle c_1(\ca L, h)^n \right\rangle 
\end{equation}
of the plurisubharmonic line bundle $(\ca L,h)$. We note that, if $h$ has bounded local potentials on a Zariski  dense open
set of $U$ of $X$, then the equality
\[ \int_X \left\langle c_1(\ca L, h)^n \right\rangle = \int_U c_1(\ca L, h)^n  
\]
holds in $\bb R_{\ge 0}$. 
In this case we say that $h$ has \emph{Zariski unbounded
  locus}. In particular such a metric has ``small unbounded locus''
in the sense of \cite{begz}. Most examples of singular metrics arising
from algebraic geometry have Zariski unbounded locus; for example,
plurisubharmonic metrics with almost asymptotically algebraic
singularities have Zariski unbounded locus. 

Fix a smooth reference metric $h_0$ on $\ca L$.
Let $\theta =c_{1}(\mathcal{L},h_{0})$ be the first Chern form of $h_0$ and set
  $\varphi=-\log(h(s)/h_{0}(s))$, where $s$ is any non-zero rational section of $\ca L$. Let $ \mathcal{J}(\varphi)$ denote the multiplier ideal associated to $\varphi$ (see \ref{def:multiplier_sheaf}).
  The non-pluripolar volume \ref{RHS} turns out to be intimately related to the \emph{multiplier ideal volume} of $(\ca L,h)$, given by the limit
 \begin{displaymath}
   \vol_{\ca J}(\ca L,h) :=  \lim _{k\to \infty}\frac{h^{0}(X,\mathcal{L}^{\otimes k}\otimes
      \mathcal{J}(k\varphi))}{k^{n}/n!}  \, .
  \end{displaymath}
Note that the right hand side is indeed independent of the choice of $h_0$. The quantity $ \vol_{\ca J}(\ca L,h)$ is called \emph{arithmetic volume} in \cite{DX}.

More precisely, Darvas and Xia show in \cite{DX}  that if $\L$ is ample we have the lower bound
\begin{equation} \label{ineq:arithm_NPP}
 \vol_{\ca J}(\ca L,h)  \ge \int_X \left\langle c_1(\ca L, h)^n \right\rangle \, . 
\end{equation}
Moreover, they show that if one assumes that the non-pluripolar volume
is strictly positive, then equality holds in \ref{ineq:arithm_NPP} if
and only if $h$ is in some precise sense very well approximable by
algebraic singularities, see \cite[Theorem 5.5]{DX} for the precise
conditions. This can be seen as an analytic Hilbert--Samuel type
statement.  We note that these results were very recently extended and
generalized to the case where $\ca L$ is only assumed to be
pseudoeffective, see \cite[Theorem 1.1]{DX21}.

Our notion of almost asymptotically algebraic singularities is
probably stronger than the list of equivalent conditions on
singularity type that is mentioned in \cite[Theorem~5.5]{DX} (see
\ref{cor:env} and the discussion afterwards), however our notion of
almost asymptotically algebraic singularities has the advantage of
being easy to verify in concrete cases. 

As a first illustration of our viewpoint, we give a direct proof of equality in \ref{ineq:arithm_NPP} in the case of almost asymptotically algebraic singularities, without the assumption that $\ca L$ is ample.
\begin{introtheorem} \label{thm:asym:intro} (\ref{th:asym}) 
 Assume that the plurisubharmonic metric $h$ has almost asymptotically
 algebraic singularities. Then the analytic Hilbert--Samuel type
 formula
\[
\vol_{\ca J}(\ca L,h) =  \int_X \langle c_1(\ca L,h)^n \rangle 
\]
holds in $\rr_{\ge 0}$.  
\end{introtheorem}

In the spirit of \cite{bkk}, the left hand side of our Chern--Weil
formula will be an intersection product of b-divisors. Assume that the plurisubharmonic metric $h$ has Zariski unbounded locus. Given a
non-zero rational section $s$ of $\ca L$, we construct an associated
$\rr$-b-divisor $D(\ca L,h,s)$ on $X$ using the \emph{Lelong numbers}
of the current $c_1(\L,h)$ at all prime divisors on all modifications
of $X$ (see \ref{b-div-lelong}). This extends a construction due to
Boucksom, Favre and Jonsson in the local case in \cite{bfj-val}.

Unfortunately, the set of all $\rr$-b-divisors does not admit a natural intersection product.  Dang and Favre \cite{Da-Fa} have shown though that  the set of so-called \emph{approximable nef} b-divisors (see \ref{nef-app}) does admit a natural intersection product with values in $\bb R$.

\begin{introtheorem} \label{thm:approx_nef_intro} (\ref{thm:psh_implies_approximable})
Assume that the plurisubharmonic metric $h$ has Zariski unbounded locus. Then the associated $\rr$-b-divisor $D(\ca L,h,s)$ on $X$ is approximable nef. 
\end{introtheorem}

In particular, the degree $D(\ca L,h,s)^n \in \bb R$ is well-defined. 
This allows us to state the following Chern--Weil type result.
\begin{introtheorem} \label{thm:b-div-chern-weil:intro} (\ref{thm:b-div-chern-weil})
 Assume that the plurisubharmonic metric $h$ has almost asymptotically algebraic singularities. Then the equality
\[
D(\ca L,h,s)^n =   \int_X \langle c_1(\ca L,h)^n \rangle 
\]
holds in $\rr_{\ge 0}$.  
\end{introtheorem}

Combining \ref{thm:asym:intro} and \ref{thm:b-div-chern-weil:intro} we obtain the equality 
\begin{equation} \label{eq:HS_intro}
 \vol_{\ca J}(\ca L,h) = D(\ca L,h,s)^n  
\end{equation}
for plurisubharmonic metrics with almost asymptotically algebraic
singularities, which can be seen as a b-divisorial version of the classical Hilbert--Samuel formula
for nef line bundles. 

It is shown in \cite{bo} and \cite{BoteroBurgos} that in the toric and
toroidal settings the b-divisorial degree of a nef b-divisor can be
computed by combinatorial means in terms of a Monge--Amp\`ere measure
associated to a weakly concave function on a suitable polyhedral
space. Also we have that the \emph{volume} (see \ref{eq:vol-b} for a precise definition) of a toric or toroidal b-divisor
agrees with its degree (see \cite[Theorem 5.11]{bo} and \cite[Theorem
5.13]{BoteroBurgos}). Such results seem to be special to the toric and toroidal cases: 
in the Appendix we give an example to show that in general the volume function is \emph{not} 
continuous on the space of big and approximable nef b-divisors (whereas the degree function is continuous).

We see that in the case of an almost asymptotically algebraic psh metric whose associated b-divisor is toroidal, the multiplier ideal volume agrees with the volume of the associated b-divisor, as both agree with the degree of the b-divisor. It would be interesting to know whether the equality between the
  multiplier ideal volume of a psh metric and the volume of the
  associated b-divisor continues to be true in the general case of almost asymptotically
  algebraic psh metrics (see \ref{rem:app}.3). \\

\noindent To illustrate our general results we shall consider the following example dealing with universal abelian varieties. This example generalizes the set-up of \cite{bkk}.

Let $g \in \zz_{\ge 1}$ and $N \in \zz_{\ge 3}$ and let $A_{g,N}$ denote the fine moduli space of principally polarized complex abelian varieties of dimension~$g$ and level~$N$. 
 Let $\pi\colon U_{g,N} \to A_{g,N}$ be the universal abelian variety,
 and $\overline{U}_{g,N}$ and $\overline{A}_{g,N}$ be any projective smooth
 toroidal compactifications of $\overline{U}_{g,N}$ and
 $\overline{A}_{g,N}$, respectively, together with a map
 $\overline{\pi}\colon \overline{U}_{g,N} \to\overline{A}_{g,N}$
 extending $\pi$, as discussed for example in
 \cite[Chapter~XI]{fc}. Let $k, m \in \zz_{\ge 0}$ and let $\L_{k,m}$
 denote the \emph{line bundle of Siegel-Jacobi forms} on $U_{g,N}$ of
 \emph{weight}~$k$ and \emph{index}~$m $. It is endowed with a
 canonical smooth invariant hermitian metric $h_{k,m}$. The first
 Chern form $c_1(\L_{k,m}, h_{k,m} )$ is a semipositive $(1,1)$-form
 on $U_{g,N}$.
 
 The next result is a special case of \ref{thm:biext_toroidal}.
\begin{introtheorem} \label{thm:SJ_is_aaas_intro} The smooth hermitian line bundle $(\L_{k,m},h_{k,m})$ has a canonical extension $( \overline{\L}_{k,m}, \overline{h}_{k,m} )$ as a $\qq$-line bundle with a plurisubharmonic metric with toroidal, and hence almost asymptotically algebraic, singularities over $\overline{U}_{g,N}$. 
\end{introtheorem}

Let $n = \dim U_{g,N} = g+ g(g+1)/2$. Let $s$ be any non-zero rational section of $\L_{k,m}$, and let $D(\overline{\L}_{k,m},\overline{h}_{k,m},s) $ be the b-divisor associated to $s$ and the metric $\overline{h}_{k,m} $. As is verified by~\ref{prop:3}, the b-divisor $D(\overline{\L}_{k,m},\overline{h}_{k,m},s) $ is independent of the choice of the compactification $\overline{U}_{g,N}$.

 \ref{thm:SJ_is_aaas_intro} together with \ref{thm:asym:intro} and \ref{thm:b-div-chern-weil:intro} imply the following. 
\begin{introtheorem}(\ref{thm:siegel-jacobi-CW})
Let notations be as above. Then the equalities
\[
\int_{U_{g,N}} c_1(\L_{k,m}, h_{k,m} )^n = D(\overline{\L}_{k,m},\overline{h}_{k,m},s)^n = \vol_{\ca J}(\overline{\L}_{k,m},\overline{h}_{k,m})
\]
hold in $\rr_{\ge 0}$.
\end{introtheorem}
Note that this generalizes \ref{eq:BKK} and \ref{eq:26} to higher degrees. In our follow-up paper \cite{future} we shall prove that the
b-divisors associated to $( \overline{\L}_{k,m}, \overline{h}_{k,m} )
$ are toroidal b-divisors on the smooth toroidal variety
$\overline{U}_{g,N}$ in the sense of \cite{BoteroBurgos}. Along the
way we show that these b-divisors are \emph{not} given by a divisor on
any single model of $\overline{U}_{g,N}$. This shows that it is
really necessary to consider these limits of divisors on all
models. Moreover, we will use the techniques developed in the present paper
to compute the asymptotic growth of the dimensions of spaces of Siegel--Jacobi
forms.

\subsection{Structure of the paper}

The purpose of \ref{anal} is to review several analytic notions such as singular hermitian metrics, (quasi-)plurisubharmonic functions and metrics, Lelong numbers, multiplier ideals and the multiplier ideal volume. We discuss various notions of singularity type for quasi-plurisubharmonic functions and review Demailly's regularization theorem as well as the notion of non-pluripolar products. We state an important monotonicity result for non-pluripolar products due to Boucksom, Eyssidieux, Guedj and Zeriahi.

In \ref{sec:psh-functions-that} we introduce the notion of almost asymptotically algebraic singularities and discuss various examples. We will see that good plurisubharmonic and toroidal plurisubharmonic metrics have almost asymptotically algebraic singularities.  We prove that in the case of almost asymptotically algebraic singularities, the multiplier ideal volume and the non-pluripolar volume coincide.

In \ref{sec:b-div} we recall the notion of b-divisors and review part of the work of Dang and Favre on approximable nef b-divisors.
In \ref{sec:b-div-psh} we show how to associate a b-divisor to a line bundle with a plurisubharmonic metric with Zariski unbounded locus. We show that such b-divisors are approximable nef, and prove our Chern--Weil type result for almost asymptotically algebraic singularities. 

Finally in \ref{sec:Siegel-Jacobi} we discuss the biextension line
bundle, and more generally the line bundle of Siegel--Jacobi forms, on
the universal abelian variety as an example of a line bundle with a
plurisubharmonic metric with almost asymptotically algebraic
singularities.

\subsection*{Acknowledgements}
We are grateful to ICMAT and to the Department of Mathematics of
Regensburg University for their hospitality during the preparation
of this paper.

\section{Analytic preliminaries}\label{anal}
In this section $X$ denotes a complex manifold of
pure dimension $n$.   

\subsection{Singular hermitian metrics}
 
We refer to \cite{dem-reg}, \cite{dem}  and \cite{ochia} (see also
\cite{bouck}, \cite{bfj-sing} and \cite{kim}) for definitions and
proofs of the analytic properties given in this section.
We take our notion of plurisubharmonic (psh) functions from
\cite[Chapter 3]{ochia}.
\begin{definition}\label{def:5}
  Let $U$ be an open coordinate subset of $X$ that we identify with an
  open subset of $\bb C^{n}$. 
  A function $\varphi\colon U\to
  \bb R\cup \{-\infty\}$ is called \emph{plurisubharmonic (psh)} if it
  satisfies the following two conditions:
  \begin{enumerate}
  \item $\varphi$ is upper semicontinuous and is not identically
    $-\infty$ on  any connected component of $U$;
  \item for every $z\in U$ and $a\in \bb C^{n}$ the function in one
    complex variable
    \begin{displaymath}
      \zeta \longmapsto \varphi(z+a\zeta)\in \bb R\cup \{-\infty\} 
    \end{displaymath}
    is either identically $-  \infty$ or subharmonic in each connected
    component of the open set $\{\zeta \in \bb C\mid z+a\zeta\in
    U\}$. 
  \end{enumerate}
  A function $\varphi \colon U\to
  \bb R\cup \{-\infty\}$ on an  arbitrary open subset $U$ of $X$ is
  called $psh$ if $U$ can be covered by open coordinate subsets
  $U_{i}$ and each $\varphi|_{U_{i}}$ is psh.
\end{definition}
We write $dd^{c}$ for the operator $\frac{i}{\pi
}\partial\bar \partial$. The following
characterization of psh functions, that does not refer to coordinate charts, is used frequently.
\begin{proposition}
  Let $U\subset X$ be an open set and let $\varphi\colon  U \to \bb R\cup \{-\infty\}$ be a measurable function. Then $\varphi$ is psh if and only
if the following 
two conditions are satisfied:
\begin{enumerate}
\item The function $\varphi$ is \emph{strongly upper semicontinuous}. That
  is, for all $V\subset U$ of total Lebesgue measure, and all $x \in U$, the condition
  \begin{displaymath}
    \varphi(x)=\limsup_{\substack{y\to x\\y\in V}}\varphi(y)
  \end{displaymath}
  holds.
\item The function $\varphi$ is locally integrable and the $(1,1)$-current
  $dd^{c}\varphi$ 
  is positive. 
\end{enumerate}
\end{proposition}
Examples of psh functions are given by the functions $\frac{1}{2}\log( |f_1|^2+ \cdots + |f_N|^2)$ where $f_1, \ldots, f_N \in \ca O_X(U)$ are non-zero holomorphic functions.
\begin{remark} \label{rem:local_potentials}
If $T$ is a closed positive $(1,1)$-current, then $T$ is locally exact and so is locally of the form $dd^c\varphi$ for some psh function $\varphi$; we call such $\varphi$ \emph{local potentials} of $T$. We say $T$ has \emph{bounded local potentials} if the $\varphi$ can be chosen to be bounded, and similarly with `continuous' or `smooth' instead of `bounded'. 
\end{remark}
When $X$ is compact, any global psh function on $X$ is constant. To
have a rich global theory we need to allow some
flexibility. 
\begin{definition}
  \begin{enumerate}
  \item An upper semicontinuous function
    $\varphi\colon X \to \bb R\cup \{-\infty\}$ is called
    \emph{quasi-plurisubharmonic (quasi-psh)} if $\varphi$ is locally of the
    form $u + f$, where $u$ is psh and $f$ is smooth.
  \item Let $\theta $ be a smooth closed $(1,1)$-form
  on $X$. 
  A measurable function $\varphi\colon X \to \bb R\cup \{-\infty\}$ is called \emph{$\theta$-psh} if $\varphi$ is 
  locally integrable, strongly upper semicontinuous and 
  $dd^c \varphi + \theta$ is a positive current.
  \end{enumerate}
\end{definition}
When $T$ is a current we write $T \ge 0$ to express that $T$ is positive.
\begin{lemma}\label{lemm:8}
  Let $\varphi\colon X\to \bb R\cup \{-\infty\}$ be a function.
  \begin{enumerate}
  \item If $\varphi$ is $\theta $-psh for some smooth closed $(1,1)$-form
    $\theta $ then $\varphi$ is quasi-psh. 
  \item If $X$ is compact K\"ahler with K\"ahler form $\omega $ and $\varphi$
    is quasi-psh then there is a real number $a>0$ such that $\varphi$
    is $(a\omega)$-psh. 
  \end{enumerate}
\end{lemma}
\begin{proof}
  For (1), assume that $\varphi$ is $\theta $-psh for some smooth closed
  $(1,1)$-form $\theta $. Then there is a covering of $X$ by
  open subsets $U_{i}$ 
  and on each $U_{i}$ there is a smooth function $f_{i}$ with $dd^{c}
  f_{i}=\theta $. Then, on $U_{i}$ we write $\varphi = (\varphi +
  f_{i}) - f_{i}$. The function $-f_{i}$ is smooth and the function
  $\varphi + f_{i}$ is locally integrable, strongly upper
  semi-continuous and satisfies
  \begin{displaymath}
    dd^{c}(\varphi + f_{i})= dd^{c}\varphi + \theta \ge 0 \, . 
  \end{displaymath}
  Therefore $\varphi + f_{i}$ is psh and $\varphi$ is quasi-psh.
  
  For (2), assume that $\varphi$ is quasi-psh and $X$ is
  compact K\"ahler with K\"ahler form $\omega $. There is a  finite open cover
  $\{U_{i}\}_{i}$ of $X$, and for each $i$ a decomposition
  $\varphi=f_{i}+\gamma_{i}$ with $f_{i}$ psh and $\gamma _{i}$
  smooth. After shrinking the $U_{i}$ if necessary we can assume that,
  for each $i$, 
  the form $dd^{c}\gamma _{i}$ can be extended to a smooth form on an open 
  neighbourhood of the compact set $\overline {U_{i}}$. Then there is
  a real number $a>0$
  such that $dd^{c}\gamma _{i}+a\omega|_{U_i} \ge 0$ for all $i$. It follows that
  \begin{displaymath}
    \left(dd^{c}\varphi + a\omega \right) |_{U_{i}}=dd^{c}f_{i}+ dd^{c}\gamma
    _{i}+a\omega|_{U_i} \ge 0
  \end{displaymath}
  for all $i$. Therefore, $\varphi$ is $(a\omega) $-psh.
 \end{proof}
In algebraic geometry it is often convenient to work with the related
concept of \emph{psh metrics} on a line bundle.  
Let $\mathcal{L}$ be a line bundle on $X$ and fix a trivialization
$\{(U_{i},s_{i})\}$ of $\mathcal{L}$ with
  transition functions $\{g_{ij}\}$.  
A \emph{hermitian metric} on $\mathcal{L}$ is a collection
of 
measurable functions 
\begin{displaymath}
h = \left\{\varphi_i \colon U_i \to \rr \cup
  \left\{\pm \infty\right\}\right\} \, ,
\end{displaymath}
 such that
  \begin{equation}\label{eq:10}
e^{-\varphi_i} = |g_{ij}|e^{-\varphi_j}
\end{equation}
on $U_i \cap U_j$.  The function $\varphi_{i}$ determines
the \emph{norm} $h(s_i)$ of the trivializing sections $s_{i}$ by the formula
\begin{displaymath}
  \varphi_{i}(z)=-\log h(s_{i}(z)),\quad z\in U_{i} \, .
\end{displaymath}
The condition \ref{eq:10} is equivalent to
\begin{displaymath}
  \log h(s_{i})- \log h(s_{j})=\log|s_{i}/s_{j}| \, .
\end{displaymath}
Then, the norm $h(s)$ of any section $s$ at a point $z$ is given by
\begin{displaymath}
  \log h(s(z))=\log|s(z)/s_{i}(z)| + \log h(s_{i}(z))=
  \log|s(z)/s_{i}(z)| - \varphi_{i}(z)  \, , 
\end{displaymath}
if $z\in  U_{i}$, and this is easily verified to be independent of the choice of $i$.
The functions $\varphi_i = -\log h(s_i)$ are called \emph{local
  potentials} of the metric $h$.
\begin{definition}\label{def:singular_metric}
The metric $h$ is called \emph{singular}
(resp.~\emph{psh}, \emph{continuous}, \emph{smooth}) if the local potentials
$\varphi_i$ are locally integrable (resp.~psh, continuous,
smooth). 
\end{definition}
The notion of singular (resp.~psh, continuous,
smooth) metric readily generalizes to the context of $\bb Q$-line bundles on $X$  (see \cite[Definition 2.10]{Biesel2014Neron-models-an} for a discussion of this terminology).
\begin{remark}\label{rem:bij-metrics}
  The global relation between psh metrics on line bundles on $X$ and $\theta $-psh functions on $X$
  is given as follows. 
  Choose a smooth reference  metric $h_{0}$ on $\mathcal{L}$ and write
  \begin{displaymath}
    \theta = c_{1}(\mathcal{L},h_{0}). 
  \end{displaymath}
  Then $\theta$ is a smooth closed $(1,1)$-form on $X$. 
  Note that for every trivializing open subset $U_i$ and trivializing section $s_i$
  of $\ca L$ over $U_i$ one has
  \begin{displaymath}
    \theta |_{U_i}=-dd^{c}\log(h_{0}(s_i)).
  \end{displaymath}
  Now choose a non-zero rational section $s$ of $\mathcal{L}$.
  Then the maps
  \[
  h \longmapsto \varphi = -\log (h(s)/h_{0}(s)) \; \text{ and }\; \varphi \longmapsto  \theta + dd^c\varphi
\]
do not depend on the choice of section. The first map is a 
  bijection between the set of psh metrics on $\mathcal{L}$ and the
  set of $\theta$-psh functions on $X$. If $X$ is compact connected, the
  second map induces a bijection between  the set of $\theta$-psh
  functions on $X$ up to constants and the set of positive $(1,1)$-currents in the
  cohomology class $c_1(\L)$ of $\L$.
  If $h$ is a psh metric on $\mathcal{L}$ and $\varphi = -\log (h(\cdot)/h_{0}(\cdot))$ the
  corresponding $\theta $-psh function,  
  we will use the notation
  $c_{1}(\mathcal{L},h)=\theta + dd^c\varphi$ for the first Chern
  current associated to $h$. 
\end{remark}

\subsection{Lelong numbers, equivalence of singularities, 
  multiplier ideals}\label{sec:lelong}
A first measure of the singularities of a psh function is given by its Lelong
numbers. 
Let $T$ be a closed positive $(1,1)$-current on the complex manifold~$X$.  
The \emph{Lelong number} $\nu(T,x)$ of $T$ at a point $x \in X$ is given by 
\[
\nu(T,x) = \lim_{r \to 0^+}\nu(T,x,r), 
\]
where $\nu(T,x,r)$ is computed in an open coordinate neighborhood of
$x$  as the integral
\[
\nu(T,x,r) = \frac{1}{(2\pi r^2)^{n-1}} \int_{B(x,r)}T(z)\wedge \left(i \partial \bar{\partial}|z|^2\right)^{n-1}.
\]
Here, $B(x,r)$ denotes the ball with center $x$ and radius $r$. Several important properties of Lelong numbers are stated in \cite[Section
2.B]{dem-anal}. In particular, they are non-negative real numbers
which are invariant under holomorphic changes of local coordinates. Further, the Lelong number $\nu(T,x)$ is
additive in $T$. If
$\varphi$ is a psh local potential of $T$ as in \ref{rem:local_potentials} then we write
\begin{displaymath}
  \nu (\varphi,x)=\nu (T,x).
\end{displaymath}
The Lelong numbers of psh functions have the following characterization.
\begin{proposition}
  Let $U\subset X$ be an open coordinate subset and let $\varphi$ be a psh
  function on $U$. Let $x \in U$.
  Then the equality
  \[
    \nu(\varphi,x) = \sup\left\{\gamma \geq 0
      \; \big{|} \; \varphi(z) \leq \gamma
      \log |z-x| + O(1) \; \text{near} \; x \right\}
  \]
  holds.
In particular, if $\varphi = \log|f|$ with $f \in 
\mathcal{O}_X(U)$ holomorphic, then $\nu(\varphi,x) = \on{ ord}_x(f)$,
where $\on{ ord}_x(f)$ is the largest power of the maximal ideal of $x$ which contains~$f$. 
\end{proposition}

\begin{definition}\label{def:modification}
A morphism $\pi\colon X' \to X$ of complex manifolds is called a \emph{modification} if it is proper and there exists a nowhere-dense analytic subset $Z \sub X$ such that the map $ \pi^{-1}(X \setminus Z) \to X \setminus Z$ given by restricting $\pi$ is an isomorphism, and such that $X' \setminus \pi^{-1} Z$ is nowhere dense in $X'$.

\end{definition}

Given a  modification $\pi \colon X' \to X $ and a psh function $\varphi$ on an analytic open subset
$U$ of $X$, the composition $\varphi \circ \pi$ is a psh function on
$\pi ^{-1}(U)$. Hence for any positive $(1,1)$-current $T$ on $X$ the Lelong number $\nu(T,x)$ at
 a point $x \in X'$ is well defined. 
Furthermore, one can define the Lelong number $\nu(T, P)$ at any prime
divisor $P$ on $X'$ by
\[
\nu(T,P) \coloneqq \nu(T,\eta),
\]
where $\eta$ is a very general point of $P$.
If $T$ is of the form $c_1(\L,h)$ for a psh metric $h$ on a line bundle $\L$, we write 
\[
\nu(h,P) \coloneqq \nu(T,P).
\]

\begin{remark}\label{rem:siu-dec}
Let $T$ be a closed positive $(1,1)$-current on $X$. We briefly recall the Siu decomposition of $T$ on $X$ (and refer to \cite[Section~2.2.1]{bouck} for details). This decomposes $T$ uniquely as a sum 
\[
T = R + \sum_k\nu\left(T, Y_k \right) \delta_{Y_k},
\]
where the sum is over an (at most countably infinite) family of $1$-codimensional subvarieties $Y_k$ of $X$. Here, $\delta_{Y_k}$ denotes the integration current determined by~$Y_k$, and $R$ is a closed positive $(1,1)$-current whose Lelong number on any prime divisor on $X$ is zero. 
In \ref{sec:b-div-psh} we will relate the Siu decomposition of $T$ to the so-called b-divisor associated to a psh metric (see \ref{rem:b-div-metric}). 
\end{remark}
The Lelong numbers allow us to classify singularities of psh functions in the sense that more singular
psh functions have bigger Lelong numbers. 
But in some situations, the classification of singularities by
Lelong numbers is too crude, and a more refined classification is needed. Such a more refined classification is given by
the concept of \emph{type of singularity} of a psh 
function.
\begin{definition}
  Let $U \subset X$ be an open subset and let $\varphi$, $\psi $
  be two psh functions on $U$. We say that $\varphi$ is \emph{more singular}
  than $\psi$ at a point $u \in U$ if there exists an open neighbourhood $u \in V \sub U$ and a constant $C \in \bb R$ such that 
 $\varphi\le \psi+C$ on $V$.  We say that $\varphi$ is more singular
  than $\psi$, denoted $\varphi\prec \psi $, if it is so at every $u \in U$. 

We write
  $\varphi\sim \psi $ if $\varphi\prec \psi $ and $\psi \prec
  \varphi$. In this case we say that $\varphi$ and $\psi $ have
  \emph{equivalent singularities}. This notion can be easily extended to the
  set of quasi-psh functions on $X$ and defines an
  equivalence relation on this set. Given a quasi-psh function
  $\varphi$ on $X$, we denote by $[\varphi]$ its equivalence
  class, and call $[\varphi]$ the \emph{type of singularity} of $\varphi$.
\end{definition}
The above classification of singularities of quasi-psh functions carries over to
closed positive $(1,1)$-currents
and to psh metrics on $\ca L$ by using local potentials.
If $T$ and $T'$ are two closed positive $(1,1)$-currents, we write
$T\prec T'$ if there is an open covering $\{U_{i}\}$ of $X$ and local
psh potentials 
$\varphi_{i}$ and $\varphi_{i}'$ of $T$ and $T'$ on $U_{i}$ such that
$\varphi_{i}\prec 
\varphi_{i}'$. Similarly, 
given two psh metrics $h$ and $h'$ on the line bundle
$\mathcal{L}$, we write
$h\prec h'$ if $c_{1}(\mathcal{L},h)\prec c_{1}(\mathcal{L},h')$.
\begin{remark}\label{rem:4}

Let $T\prec T'$ be two closed positive $(1,1)$-currents on $X$. Then their Lelong
numbers satisfy
\begin{displaymath}
  \nu (T,x) \ge \nu(T',x)
\end{displaymath}
for every point $x\in X'$ in any modification $\pi \colon X' \to X$ of $X$. In
particular if $T$ and $T'$ have equivalent singularities, then they have the same
Lelong numbers at every prime divisor on every modification of $X$. 

The converse does not hold as the following example shows. The example is local for
  ease of writing but can easily be made into a global one.
\end{remark}
\begin{example}\label{exa:sing-type}
Consider the function $f(z)=-\log(-\log(z\bar z))$ on the disk
  $z\bar z<1$. Then $\nu (f,q)=0$ for all $q$ in the disk, but $f$ is not bounded
  below. So $f $ has the same Lelong numbers as a constant
  function in any point but $f\not \sim 1$. This is a
  counterexample to  the converse of \ref{rem:4} because any
  modification of the disk is the disk itself. 
\end{example}

Two psh functions which have the same Lelong numbers on $X$ need not have the same Lelong numbers on a modification. In fact,
more is true: it can happen that $\nu(f,x) \le \nu(g,x)$ for all $x
\in X$ and even $\nu(f,x_{0}) < \nu(g,x_{0})$ for a point  $x_{0} \in X$, but $\nu(f,y) >
\nu(g,y)$ for some point $y$  above $x_{0}$ in some modification $X'$ of $X$.

\begin{example}
Consider the functions $f$ and $g$ in the polydisk 
  \begin{displaymath}
    E=\{(x,y)\in \mathbb C^{2}\mid x\bar x <1,\, y\bar y <1\}
  \end{displaymath}
  given by
  \begin{displaymath}
    f(x,y)=\frac{.9}{2}\log(x\bar x+(y\bar y)^2) \, ,\qquad
    g(x,y)=\frac{1}{2}\log((x\bar x)^2+y\bar y) \, ,
  \end{displaymath}
  Then $\nu (f,0)=.9 < 1.0 = \nu(g,0)$. Nevertheless
  \begin{displaymath}
    g(0,y)-f(0,y)=-.4\log(y\bar y)
  \end{displaymath}
  is not bounded above. So $\nu (f,q) \le \nu(g,q)$ for any point $q$ in
  the polydisk but $g\not \prec f$.

 However, if we consider the chart of the blow-up of the polydisk at
  $(0,0)$ with coordinates $(s,t)$ with  $x=st$, $y=s$ and if we let $p$ be the
  point $s=t=0$ on the blow-up, then $\nu (f,p)=1.8$ while $\nu (g,p)=1$.
 
  \end{example}

A convenient way to encode Lelong numbers is by means of 
multiplier ideal sheaves.  
\begin{definition} \label{def:multiplier_sheaf}
  Let $\varphi$ be a quasi-psh function on $X$. Then  
the \emph{multiplier ideal sheaf} $\mathcal{J}(\varphi)$ of $\varphi$ is the
coherent ideal sheaf of $\ca O_X$-modules given locally by those holomorphic functions $f$
such that $|f|^2e^{-2\varphi}$ is locally integrable. If $\mathcal{L}$
is a line bundle with a psh metric $h$ then the multiplier ideal sheaf $\mathcal{J}(h)$ of
$h$ is defined to be the multiplier ideal sheaf of the quasi-psh function $\varphi=-\log(h(s)/h_{0}(s))$ for any
smooth reference metric $h_{0}$ and any non-zero rational section $s$ of $\ca L$. 
\end{definition}
It turns out that the multiplier ideal sheaves of all multiples of a
  given quasi-psh function on $X$ give the same information as the Lelong
  numbers on all points of all modifications of $X$. See \ref{prop:alg-sing-type} below for a precise statement. Hence it follows from \ref{exa:sing-type} that for a general quasi-psh function also its multiplier ideal is not enough to recover the singularity type.

\subsection{Algebraic singularities}
\label{sec:algebr-sing-dema}
From now on we will assume that $X$ is a projective complex manifold. 
An important class of quasi-psh functions on $X$ consists of those having
algebraic singularities.
\begin{definition}\label{def:3}
A quasi-psh function
$\varphi$ on $X$ is said to have \emph{algebraic singularities} if
there is a constant $c\in \bb Q_{\ge 0}$ and $\varphi$
can be written locally as  
  \begin{equation}
\varphi = \frac{c}{2} \log \left(|f_1|^2 + \dotsm + |f_N|^2\right) +
\lambda,\label{eq:9}
\end{equation}
where $\lambda$ is a bounded function
and the $f_j$ are non-zero algebraic functions.
\end{definition}
\begin{remark}
  There is also the related notion of \emph{analytic}
  singularities. In the projective context the only difference with the notion of algebraic singularities is to
  allow the constant $c$ to be a real number.  
\end{remark}
Following \cite[1.10]{dem-anal}, if the quasi-psh function $\varphi$ on $X$ has
algebraic singularities with constant $c$, we can associate to it a
coherent sheaf of 
ideals $\mathcal{I}(\varphi/c)$, in the following way. Since $X$ is compact, we can
assume that there 
is a finite covering of $X$ and $\varphi$ has the form \ref{eq:9} on
each open of the covering.
Then $\mathcal{I}(\varphi/c)$ is defined to be the ideal sheaf
of holomorphic functions $h$ satisfying
\begin{displaymath}
  |h|\le C(|f_{1}|+\dots+|f_{N}|)
\end{displaymath}
for some constant $C$.
This is a globally defined ideal sheaf of $\ca O_X$ locally equal to the integral
closure of the ideal generated by $(f_{1},\dots,f_{N})$. Since $X$ is
assumed to be projective, the coherent ideal sheaf $\mathcal{I}(\varphi/c)$ is the analytification of an algebraic 
coherent ideal sheaf on $X$.

In contrast with what happens for arbitrary quasi-psh functions,  
for quasi-psh functions with algebraic singularities, we can recover
the singularity type from the multiplier ideal sheaf. Following
\cite[Remark 5.9]{dem-anal}, the multiplier ideal sheaf $\ca J(\varphi)$ of a
quasi-psh function $\varphi$ with algebraic singularities is easy to
describe. Assume first that there is a constant $c \in \bb Q_{\ge 0}$ and an effective Cartier divisor with simple
normal crossings $D=\sum_{i} \alpha _{i}D_{i}$ on $X$ where the $D_i$ are irreducible such that locally $\varphi$ can be written as
\begin{equation}\label{eq:8}
  \frac{c}{2}\log |g|^2 + \lambda ,
\end{equation}
where $\lambda $ is bounded and $g$ is a local equation for $D$. Then 
\begin{equation}
  \label{eq:17}
  \mathcal{J}(\varphi) = \mathcal{O}_X\Big(-\sum_{i}\lfloor c\alpha
  _{i}\rfloor D_{i}\Big).
\end{equation}
In particular $\mathcal{J}((k/c) \varphi)=\mathcal{I}(k\varphi/c)$
for every integer $k\ge 0$.

Assume now that $\varphi$ has algebraic singularities and, as before, let $c \in \bb Q_{\ge 0}$ be
the constant appearing in \ref{def:3}. Then there exists a
modification $\pi \colon X_{\pi }\to X$  such that $\pi ^{-1}\mathcal{I}(\varphi/c)\cdot \ca O_{X_\pi}=\mathcal{O}_{X_\pi}(-D)$
for a simple normal crossings divisor $D=\sum_{i}\alpha  _{i}D_{i}$ on
$X_{\pi }$. Let 
$R_{\pi }=\sum_{i}\rho _{i}D_{i}$ be the zero divisor of the Jacobian
function of $\pi $. Then combining \ref{eq:17} with the direct image
formula 
\cite[Proposition 5.8]{dem-anal} we obtain
\begin{equation}
  \label{eq:18}
  \mathcal{J}(\varphi) = \pi _{\ast}\mathcal{O}_{X_\pi}(R_\pi-\sum_{i}\lfloor c
  \alpha _{i}\rfloor D_{i})=\pi _{\ast}\mathcal{O}_{X_\pi}(\sum_{i}(\rho_i -\lfloor c
  \alpha _{i}\rfloor )D_{i}).
\end{equation}
Therefore, for every integer $k> 0$
\begin{equation} \label{multiples_multiplier}
  \mathcal{J}((k/c) \varphi)=\pi _{\ast}\big(\mathcal{O}_{X_\pi}(R_\pi
  )\otimes (\pi ^{-1}\mathcal{I}(k\varphi/c)\cdot \ca O_{X_\pi})\big).
\end{equation}

This description shows that, if $\varphi$ has algebraic singularities
with constant $c$, then the asymptotic properties of the family of
ideals $\mathcal{J}((k/c)\varphi)$, $k> 0$ and that of the family
$\mathcal{I}((k\varphi)/c)$, $k>0$ are similar. The following is an
example of this property.

\begin{lemma} \label{lemm:10}
  Let
  $\mathcal{L}$ be a line bundle on the projective complex manifold~$X$
  provided with a smooth reference metric $h_{0}$ and a psh metric
  $h$. Let $\theta =c_{1}(\mathcal{L},h_{0})$ and let
  $\varphi=-\log(h(s)/h_{0}(s))$ be the resulting $\theta $-psh function as in
  \ref{rem:bij-metrics}, where $s$ is any non-zero rational section of $\ca L$. If $\varphi$ has algebraic singularities with
  constant $c$, then 
  \begin{displaymath}
    \lim _{k\to \infty}\frac{h^{0}(X,\mathcal{L}^{\otimes k}\otimes
      \mathcal{J}(k\varphi))}{k^{n}/n!} =
    \lim _{\substack{k\to \infty\\ kc\in \bb Z}}\frac{h^{0}(X,\mathcal{L}^{\otimes k}\otimes
      \mathcal{I}((kc\varphi)/c))}{k^{n}/n!}.
  \end{displaymath}
\end{lemma}
\begin{proof}
  Let $\pi \colon X_{\pi }\to X$ be a modification such that $\pi ^{-1}\mathcal{I}(\varphi/c)\cdot \ca O_{X_\pi}=\mathcal{O}(-D)$, with $D=\sum_{i}\alpha
  _{i}D_{i}$ a simple normal crossings divisor on $X_\pi$. As before let $R_{\pi }=\sum_{i}\rho _{i}D_{i}$ be
  the zero divisor of the Jacobian function of $\pi $. Then by the
  descriptions \ref{eq:18} and \ref{multiples_multiplier} of the multiplier ideal in the case of
  algebraic singularities we deduce
  \begin{align*}
    \lim _{k\to \infty}\frac{h^{0}(X,\mathcal{L}^{\otimes k}\otimes
    \mathcal{J}(k\varphi))}{k^{n}/n!}
    & =
    \lim _{k\to \infty}\frac{h^{0}\big(X_{\pi },\pi ^{\ast}\mathcal{L}^{\otimes k}\otimes
      \mathcal{O}(\sum_{i} (\rho_{i} -\lfloor kc
      \alpha _{i}\rfloor) D_{i})\big)}{k^{n}/n!}\\
    & = \lim _{\substack{k\to \infty\\ kc\in \bb Z}}\frac{h^{0}\big(X_{\pi },\pi ^{\ast}\mathcal{L}^{\otimes k}\otimes
      \mathcal{O}(\sum_{i} (\rho_{i} -kc
      \alpha _{i}) D_{i})\big)}{k^{n}/n!}\\
       & = \lim _{\substack{k\to \infty\\ kc\in \bb Z}} \left( \frac{h^{0}\big(X_{\pi },\pi ^{\ast}\mathcal{L}^{\otimes k}\otimes
      \mathcal{O}(-\sum_{i}kc
    \alpha _{i}D_{i})\big) }{k^{n}/n!} + \frac{O\left(k^{n-1}\right)}{k^{n}/n!} \right) \\
    & = \lim _{\substack{k\to \infty\\ kc\in \bb Z}}\frac{h^{0}\big(X_{\pi },\pi ^{\ast}\mathcal{L}^{\otimes k}\otimes
      \mathcal{O}(-\sum_{i}kc
    \alpha _{i}D_{i})\big)}{k^{n}/n!}\\
    & = \lim _{\substack{k\to \infty\\ kc\in \bb Z}}\frac{h^{0}(X,\mathcal{L}^{\otimes k}\otimes
      \mathcal{I}((kc\varphi)/c))}{k^{n}/n!}.\qedhere
  \end{align*}
\end{proof}

The next lemma shows that, after a modification, a quasi-psh function with
algebraic singularities is always of the form \ref{eq:8}.

\begin{lemma}\label{lemm:1} Let $\varphi$ be a quasi-psh function on $X$ with
  algebraic singularities and let $c$ be the constant  in equation
  \ref{eq:9}. Let $\pi \colon X_{\pi }\to X$ be a
  modification such that $\pi ^{-1}\mathcal{I}(\varphi/c)\cdot \ca O_{X_\pi}$ is locally
  principal. Let $U\subset X_{\pi }$ be an open subset such that $\pi ^{-1}\mathcal{I}(\varphi/c)\cdot \ca O_{X_\pi}$ is generated by a holomorphic function
  $g$ on $U$. Then
  \begin{displaymath}
    \pi ^{\ast}\varphi-\frac{c}{2}\log|g|^{2}
  \end{displaymath}
  is locally bounded on $U$. 
\end{lemma}
\begin{proof}
  After shrinking $U$ if necessary, we can assume that there is an
  open set $V\subset X$ where $\varphi$ has the shape \ref{eq:9} and
  $U\subset \pi ^{-1}(V)$. To simplify the notation we will not
  distinguish between functions on $X$ and on $X_{\pi }$ as these
  spaces agree on a dense open subset. 

  Since $g$ is a generator of $\mathcal{I}(\varphi/c)$ and the
  functions $f_{i}$ belong to this ideal, there are holomorphic
  functions $b_{i}$ such that $f_{i}=b_{i}g$. Then
  \begin{displaymath}
    \varphi = \frac{c}{2}\log\left(|f_{1}|^{2}+\dots
      +|f_{N}|^{2}\right) + \lambda =
    \frac{c}{2}\log |g|^{2}+ \frac{c}{2}\log\left(|b_{1}|^{2}+\dots
      +|b_{N}|^{2}\right) + \lambda.
  \end{displaymath}
  Since $\lambda $ is locally bounded, it suffices to prove that
  $\log\left(|b_{1}|^{2}+\dots +|b_{N}|^{2}\right) $ is locally
  bounded. Assume that this is  not the case. Since $|b_{1}|^{2}+\dots
  +|b_{N}|^{2}$ is continuous, the only possibility  for the logarithm
  not to be locally bounded is that the functions $b_{i}$ have a
  common zero. Assume that $x$ is such that $b_{i}(x)=0$ for all
  $i$. Since $\mathcal{I}(\varphi/c)_{x}$ is the integral closure of
  the ideal  $I=(f_{1},\dots,f_{N})$, there exist an integer $r\ge 1$ and elements $\alpha
  _{j}\in I^{j}$ for $j=1,\ldots,r$ such  that
  \begin{equation}\label{eq:22}
    g^{r}+\alpha _{1}g^{r-1}+\dots+\alpha _{r}=0.
  \end{equation}
  Assume that $\on{ord}_{x}(g)=k$. Since $b_{i}(x)=0$ for all $i$, we have that
  $\on{ord}_{x}f_{i}> k$ for all $i$. Hence any element $\alpha
  _{j}\in I^{j}$ has $\on{ord}_{x}\alpha _{j} > jk$. Then condition
  \ref{eq:22} implies that $\on{ord}_{x}g^{r} > kr$, which contradicts
  the fact that $\on{ord}_{x}g^{r} = kr$. We conclude that the functions $b_{i}$
  do not have a common zero and $\log\left(|b_{1}|^{2}+\dots
    +|b_{N}|^{2}\right) $ is locally  bounded. 
\end{proof}

\subsection{Demailly's regularization theorem}

Next we discuss Demailly's regularization theorem for closed
positive $(1,1)$-currents. Roughly speaking it states that, for $X$ a projective complex manifold, any $\theta $-psh
function on $X$ can be approximated by functions with algebraic
singularities. There are several versions in the literature, depending on the properties we
want for the approximating functions. 
The version we use here can be obtained combining the local version in
\cite[Theorem 4.2]{dk} with the global version in \cite[Theorem
3.2]{DP} (see also \cite[Theorem 13.2]{dem-anal} and 
\cite[Theorem 1.1]{dem-reg} for the original statement and the heart
of the proof). We will use Demailly's regularization theorem in \ref{sec:criteria} to give criteria that ensure a psh metric has almost asymptotically algebraic singularities.
\begin{theorem}\label{thm:2} (Demailly's regularization theorem) Let $X$ be a projective complex manifold
  of dimension~$n$ with K\"ahler form $\omega $ and let $\theta $ be a
  smooth $(1,1)$-form on $X$. 
Let $T$ be a closed positive $(1,1)$-current in the same cohomology
class as $\theta$. Write $T = \theta + dd^c\varphi$ with $\varphi$ a
$\theta$-psh function.  Then there exists a sequence
$(\varphi_m)_{m \geq 1}$ of quasi-psh functions on $X$
satisfying the following properties.
\begin{enumerate}
\item\label{item:1} Each $\varphi_m$ has algebraic
  singularities. Moreover, for each $m$ there is a modification
  $\pi_{m} \colon X_{\pi _{m}}\to X$,  a simple normal crossings
  divisor $D_{m}$ on $X_{\pi _{m}}$, and a rational number $c
  _{m}>0$  such that, locally on $X_{\pi _{m}}$,
  \begin{displaymath}
    \varphi_{m}\circ \pi _{m}=c _{m}\log|g| + f
  \end{displaymath}
  where $g$ is a local equation for $D_{m}$ and $f$ is smooth (and not
  just locally bounded as in \ref{lemm:1}).
\item The sequence  $(\varphi_m)_{m \geq 0}$ is non-increasing and there is a sequence of positive real numbers
  $(a_{m})_{m\ge 0}$ converging monotonically to zero, such that
  $\varphi_{m}$ is $(\theta +a_{m}\omega) $-psh.
\item\label{item:2} For every coordinate open set $U$ and relatively
  compact open subset $V\subset\subset U$ there
  are constants $C_{1},C_{2}>0$ such that for all $m \in \bb Z_{+}$,
  $z\in V$ and $r \in \bb R_{+}$ with $r<d(z,\partial V)$, 
   \begin{equation}\label{eq:6}
    \varphi(z) -\frac{C_{1}}{m}\le \varphi_{m}(z)\le
    \sup_{\|x-z\|<r}\varphi(x)+\frac{1}{m}\log\frac{C_{2}}{r^{n}}. 
  \end{equation}
\item\label{item:3} For all $x\in X$, the Lelong numbers of $\varphi $
  and $\varphi _{m}$ satisfy the
  condition 
  \begin{equation}\label{eq:7}
    \nu (\varphi,x)-\frac{n}{m}\le \nu (\varphi _{m},x)\le \nu (\varphi ,x).
  \end{equation}
  In particular the Lelong numbers of the functions $\varphi _{m}$
  on the points of $X$ converge monotonically and uniformly to the
  Lelong numbers of the function $\varphi $.  
\end{enumerate}
\end{theorem}
We refer to any sequence of approximations $(\varphi_m)_{m \ge 1}$ with properties (1)--(4) from the theorem as a \emph{Demailly
  approximating sequence} of the function $\varphi$.
  Note that in particular $\varphi \prec \varphi_m$ for a Demailly approximating sequence.

\subsection{Non-pluripolar products}
Generalizing a construction due to Bedford and Taylor \cite{BT}, it has
been shown in \cite{begz} that given closed
positive $(1,1)$-currents $T_1, \dotsc, T_p$ on the projective complex manifold $X$ one has a
non-pluripolar product 
\[
\left\langle T_1 \wedge \dotsm \wedge T_p\right\rangle
\]
of these currents  
with good properties.
The result is a \emph{closed} positive $(p,p)$-current which does not
charge pluripolar sets. As before let $n = \dim X$.
\begin{definition}\label{def:7} Let $T_1, \dotsc, T_p$ be closed
  positive $(1,1)$-currents on $X$. The \emph{non-pluripolar product}
  \[
\left\langle T_1 \wedge \dots \wedge T_p\right\rangle
\]
is the $(p,p)$-current on $X$ determined as follows. 
For $i=1,\dots, p$ let $\theta _{i}$ be a smooth $(1,1)$-form in
the same cohomology class as $T_{i}$ and let $\varphi_{i}$ be a
$\theta _{i}$-psh function with $\theta
_{i}+dd^{c}\varphi_{i}=T_{i}$. For every $k\ge 0$ write $U_{k}$ for
the set
\begin{displaymath}
  U_{k}  = \{x\in X \, | \, \varphi_{i}(x)\ge -k \, ,\,  i=1,\dots,p\} \, .
\end{displaymath}
Then for every smooth $(n-p,n-p)$-form $\eta$ one sets
\begin{displaymath}
  \left\langle T_1 \wedge \dots \wedge T_p\right\rangle(\eta) =
  \lim_{k\to \infty} \int_{U_{k}}T_1 \wedge \dots \wedge T_p\wedge \eta.
\end{displaymath}
The current $\left\langle T_1 \wedge \dots \wedge T_p\right\rangle $ 
is independent of the choices of the $\theta_i$ and $\varphi_i$.
\end{definition}

\begin{remark}
  For the existence of the non-pluripolar product a hypothesis like
  $X$ being K\"ahler is needed (see \cite[Proposition 1.6]{begz}, and the
  examples before and the remark after it). 
\end{remark}

\begin{remark} \label{finite_improper_integral} It is clear from the definition that when $T_1,\ldots,T_n$ restrict to smooth differential forms
on the dense open $U \subset X$, we have
\[ \int_X \langle T_1 \wedge \dots \wedge T_n \rangle = \int_U T_1 \wedge \dots \wedge T_n \, . \]
In particular we see that in this case the improper integral $ \int_U T_1 \wedge \dots \wedge T_n$ is well-defined and finite.
\end{remark}

In the case of currents with locally bounded potentials, the
non-pluripolar product agrees with the  cohomology product, as is shown by
the next lemma.

\begin{lemma}\label{lemm:2}

  Let $T_{1}, \dots,T_{k}$ be closed positive $(1,1)$-currents on $X$
  with locally bounded potentials and let
  $\theta _{k+1},\dots ,\theta _{n}$ be smooth closed $(1,1)$-forms on
  $X$. Choose smooth closed $(1,1)$-forms $\theta _{i}$ in the
  cohomology class of $T_{i}$ for $i=1,\dots,k$. Let $\varphi_{i}$ for $i=1,\ldots,k$ be
  locally bounded quasi-psh functions satisfying $\theta
  _{i}+dd^{c}\varphi_{i}=T_{i}$. Then
  \begin{displaymath}
    \int_{X}\theta _{1}\wedge \dots \wedge \theta _{n}=
    \int_{X}T _{1}\wedge \dots \wedge T_{k}\wedge \theta _{k+1}\wedge
    \dots \wedge \theta _{n} =
        \int_{X}\langle T _{1}\wedge \dots \wedge T_{k}\wedge \theta
        _{k+1}\wedge \dots \wedge \theta _{n}\rangle, 
      \end{displaymath}
      where the middle product is defined by Bedford-Taylor theory \cite{BT82}.
      In particular, the cohomology classes $\class{\theta_i}=\class{T_i}$ are nef for $i=1,\ldots,k$. 
\end{lemma}
\begin{proof}
  Since the currents $T_{i}$ are positive, the functions $\varphi_{i}$
  are $\theta _{i}$-psh. Since they are locally bounded and $X$ is
  compact, they are bounded. Then, \ref{def:7} of the non-pluripolar
  product immediately implies the equality of the second and third
  integrals.

  We prove the first equality by induction on $k$. If $k=0$ there is
  nothing to prove. Assume that $k>0$ and that the result is true for
  $k-1$. Bedford-Taylor theory \cite{BT82} provides us with positive currents
  \begin{displaymath}
    T_{1}\wedge\dots\wedge T_{k-1}, \qquad T_{1}\wedge\dots\wedge T_{k},
  \end{displaymath}
  and a current $\varphi_{k}T_{1}\wedge\dots\wedge T_{k-1}$ satisfying
  \begin{displaymath}
    dd^{c}(\varphi_{k}T_{1}\wedge\dots\wedge T_{k-1})=
    T_{1}\wedge\dots\wedge T_{k}-
    T_{1}\wedge\dots\wedge T_{k-1}\wedge \theta _{k}.
  \end{displaymath}
Therefore, using the induction hypothesis and the fact that the
integral over $X$ of an exact current is zero we obtain
\begin{displaymath}
  \int_{X}T _{1}\wedge \dots \wedge T_{k}\wedge \theta _{k+1}\wedge
  \dots \wedge \theta _{n}
=  \int_{X}T _{1}\wedge \dots \wedge T_{k-1}\wedge \theta _{k}\wedge
\dots \wedge \theta _{n}
=  \int_{X}\theta _{1}\wedge \dots \wedge \theta _{n}.
\end{displaymath}
This proves the first equality.
Finally, let $C$ be a closed curve on $X$. Then for each $i=1,\ldots,k $ we have
\[ \class{ \theta_i } \cdot C =  \int_C \theta_i = \int_C \langle T_i \rangle \ge 0 \, . \]
This shows that the cohomology class $\class{\theta_i}=\class{T_i}$ is nef.
\end{proof}

The non-pluripolar product is clearly symmetric and it is multilinear
in the following sense.
\begin{proposition}[{\cite[Proposition 1.4]{begz}}]\label{prop:6}
  Let $T'_{1}$, $T_{1}$, $T_{2},\dots,T_{p}$ be closed positive $(1,1)$-currents. 
  Then for every pair of positive real numbers $\alpha$,
  $\beta $ the relation
  \begin{displaymath}
    \langle (\alpha T_{1}+\beta T'_{1})\wedge T_{2}\wedge \dots \wedge
    T_{p} \rangle =
    \alpha  \langle T_{1}\wedge T_{2}\wedge \dots \wedge
    T_{p} \rangle +
    \beta \langle T'_{1}\wedge T_{2}\wedge \dots \wedge
    T_{p} \rangle
  \end{displaymath}
  is satisfied.
\end{proposition}
The following is a monotonicity property of non-pluripolar
products with respect to singularity type. Following  \cite[Definition~1.2]{begz} we say that a quasi-psh function
$\varphi$ on $X$ has \emph{small unbounded locus} if there exists a (locally) complete
pluripolar closed subset $A$ of $X$ such that $\varphi$ is locally bounded on $X \setminus A$.

\begin{theorem}[{\cite[Theorem 1.16]{begz}}] \label{thm:1} Let $\theta $ be a smooth closed $(1,1)$-form.
 For $i=1,\dots,n$, let $\{\varphi_i\}$ and $\{\psi_i\}$ be two
 collections of $\theta$-psh functions with small unbounded locus such that $\varphi_{i}\prec
 \psi _{i}$ for all $i$.  Then 
  the non-pluripolar products satisfy  
  \begin{displaymath}
    \int_{X}\langle (\theta+dd^{c}\varphi_{1})\wedge\dots\wedge
    (\theta +dd^{c}\varphi_{n})\rangle \le
        \int_{X}\langle (\theta+dd^{c}\psi _{1})\wedge\dots\wedge
    (\theta+dd^{c}\psi _{n})\rangle.
  \end{displaymath}
\end{theorem}
The main technical difficulty we face now is that in Demailly's
regularization theorem (\ref{thm:2}) we have very good control on the Lelong
numbers of a quasi-psh function by means of approximating sequences, but not on
the type of singularity. By contrast, in order to apply \ref{thm:1} on the
monotonicity of 
non-pluripolar products  we need control on the  type of singularity
of our quasi-psh functions. In the long run this will force us to restrict
the space of quasi-psh functions we can consider.

\subsection{Algebraic singularity type}
\label{sec:algebr-type}
Multiplier ideals allow us to define the algebraic type of singularity of a quasi-psh function. The notion of algebraic type has been introduced in \cite{KS19}. We will follow \cite{DX} though as the main source for our discussion.
We continue to assume that $X$ is a projective complex manifold.
\begin{definition}\label{def:alg-sing-type}
  Let $\varphi$ and $\psi $ be two quasi-psh functions on $X$. Then $\varphi$
  is said to be \emph{algebraically more singular} than $\psi $ (denoted
  $\varphi\prec_{\mathcal{J}}\psi $) if for all real numbers $a>0$ the
  inclusion $\mathcal{J}(a\varphi)\subset \mathcal{J}(a\psi)$ holds. We say
  that $\varphi$ and $\psi $ have the same 
  \emph{algebraic singularity type}, denoted
  $\varphi\simeq_{\mathcal{J}}\psi $, if
  $\varphi\prec_{\mathcal{J}}\psi $ and $\psi\prec_{\mathcal{J}}\varphi $.   
\end{definition}
The algebraic singularity type is governed by the Lelong
numbers not just on $X$ but on all modifications of $X$, as the following result shows.
\begin{proposition}[{\cite[Corollary 2.16]{DX}}]\label{prop:alg-sing-type} Let 
  $\varphi$ and $\psi $ be two quasi-psh functions on~$X$. The
  following assertions are equivalent:
  \begin{enumerate}
  \item $\varphi\prec_{\mathcal{J}}\psi $;
  \item for every modification $Y\to X$  and for every 
    $y\in Y$ the inequality $\nu (\varphi,y)\ge \nu (\psi,y)$
    holds.
  \end{enumerate}
\end{proposition}
The type of singularity and the algebraic type of singularity allow us
to attach two envelopes to a quasi-psh function. The first one was
introduced in \cite{RWN14} and the second in \cite{DX}.
\begin{definition}\label{def:env}
  Let $\theta $ be a smooth closed $(1,1)$-form and $\varphi$ a
  $\theta $-psh function on~$X$. Then the \emph{envelope of the singularity
    type} of $\varphi$ is the function
  \begin{displaymath}
    P[\varphi]=\sup\{\psi \text{ $\theta$-psh } \mid \psi \prec \varphi, \psi
    \le 0\}^{\ast}
  \end{displaymath}
  on $X$, where $f^{\ast}$ denotes the upper semicontinuous regularization of $f$.
  The \emph{envelope of the algebraic singularity type} of $\varphi$
  is 
  \begin{displaymath}
    P[\varphi]_{\mathcal{J}}=\sup\{\psi \text{ $\theta $-psh} \mid \psi
    \prec_{\mathcal{J}} \varphi, \psi \le 0\}^{\ast}.
  \end{displaymath}
\end{definition}
The following are basic properties of the envelopes.
\begin{proposition}[{\cite[Proposition 2.19]{DX}}]\label{prop:eq-sing-type} Let $\varphi$ be a
  $\theta $-psh function on~$X$. Then
  \begin{enumerate}
  \item $P[\varphi]$ and $P[\varphi]_{\mathcal{J}}$ are $\theta $-psh functions on $X$.
  \item $P[\varphi]_{\mathcal{J}}\simeq_{\mathcal{J}}\varphi$.
  \item $\varphi \prec P[\varphi]\prec P[\varphi]_{\mathcal{J}}$. Moreover
    $P\left[P[\varphi]_{\mathcal{J}}\right]=P[\varphi]_{\mathcal{J}}$. 
  \end{enumerate}
\end{proposition}
The following result due to Darvas and Xia indicates that the difference between the envelope of the singularity type and the
envelope of the algebraic singularity type governs when the non-pluripolar product is well behaved.

Let $\mathcal{L}$ be a line bundle on~$X$
  provided with a smooth reference metric $h_{0}$ and a psh metric
  $h$. Let $\theta =c_{1}(\mathcal{L},h_{0})$ be the first Chern form and
  $\varphi=-\log(h(s)/h_{0}(s))$ the resulting $\theta $-psh function.

\begin{theorem}[{\cite[Theorem~5.5]{DX}}]\label{thm:3}  Assume that $\ca L$ is ample and $\theta$ is a K\"ahler form on $X$. Let $n = \dim X$. 
Then the limit
  \begin{equation} \label{eqn:arithm_volume} \vol_{\ca J}(\ca L,h) = \lim _{k\to \infty}\frac{h^{0}(X,\mathcal{L}^{\otimes k}\otimes
      \mathcal{J}(k\varphi))}{k^{n}/n!} 
   \end{equation}
  exists, and we have
  \begin{displaymath}
    \vol_{\ca J}(\ca L,h) =
    \int_{X}\langle (\theta +dd^{c}P[\varphi]_{\mathcal{J}})^{\wedge n}\rangle
    \ge \int_{X}\langle (\theta +dd^{c}\varphi)^{\wedge n}\rangle \, . 
  \end{displaymath}
If $\int_{X}\langle (\theta +dd^{c}\varphi)^{\wedge
    n}\rangle >0$ then the
  following assertions are equivalent:
  \begin{enumerate}
  \item$ \vol_{\ca J}(\ca L,h) =
    \int_{X}\langle (\theta +dd^{c}\varphi)^{\wedge n}\rangle $;
  \item $P[\varphi]=P[\varphi]_{\mathcal{J}}$.
  \end{enumerate}
\end{theorem}

The limit in \ref{eqn:arithm_volume} is called the \emph{multiplier ideal volume} of the pair $(\ca L,h)$, whereas the integral $    \int_{X}\langle (\theta +dd^{c}\varphi)^{\wedge n}\rangle$ is called the \emph{non-pluripolar volume} of $(\ca L,h)$.
In \cite[Theorem 5.5]{DX} one may find other statements equivalent to conditions 1.\ and 2.\ above. The precise statements will not be needed here but roughly speaking the equality $P[\varphi]=P[\varphi]_{\mathcal{J}}$ holds if and only if
$\varphi$ can be ``very well approximated'' (in what is called the
\emph{$d_{\mathcal{S}}$-distance}) by algebraic singularities. 

We take the above result as an indication that for a reasonable Chern--Weil and Hilbert-Samuel theory to hold, one should deal with quasi-psh functions that are well approximated by algebraic singularities in some sense. 
\begin{remark}
As was mentioned in the introduction, we note that the above result was very recently extended and generalized to the case where $\ca L$ is only assumed to be pseudoeffective, see \cite[Theorem 1.1]{DX21}.
\end{remark}

\subsection{Good psh metrics}

Based on this idea, in \ref{sec:psh-functions-that} we will introduce a large class of singularities where the equality
of non-pluripolar volume and multiplier ideal volume can be seen to be satisfied. For this class (the class of ``almost
asymptotically algebraic singularities'') we will be able to prove a reasonable Chern--Weil type statement. As a warm-up,
we examine here already the case of ``good'' metrics in the sense of Mumford \cite{hi} and in the next section
the case of algebraic singularities.
\begin{example}\label{exm:7}
  Consider  $X=\bb P^{1}$ with homogeneous coordinates $(x:y)$
  and absolute coordinate $t=x/y$. Let $\mathcal{L}=\mathcal{O}_X(1)$.
  The space of global sections of $\mathcal{L}$ can be identified
  with the space of linear forms in the variables $x,y$. We consider
  the psh metric $h$ on $\ca L$ given by
  \begin{displaymath}
    -\log h(y) =
    \begin{cases}
      1+\log |t|,& \text{ if }|t|\ge 1/e,\\
      -\log(-\log(|t|))& \text{ if }|t|\le 1/e.
    \end{cases}
  \end{displaymath}
  This metric is singular at the point $t=0$. The interest of this
  singularity is that (up to multiplying by a normalization factor) it
  is equivalent to the singularity of the invariant (i.e., Hodge) metric on the line
  bundle of modular forms on a modular curve at a cusp (see \ref{sec:Siegel-Jacobi-forms} for further discussion). 
  
  We choose now a smooth metric $h_{0}$ on $\bb P^{1}$. A canonical choice is the
  Fubini-Study metric, given by
  \begin{displaymath}
    -\log h_0(y) = \frac{1}{2}\log(1+|t|^{2}).
  \end{displaymath}
  Let $\omega $ be the first Chern form of $(\mathcal{L},h_{0})$. Then
  the function
  \begin{displaymath}
    \varphi(t)=
    \begin{cases}
      1+\log |t|-\log(1+|t|^{2})/2 & \text{ if } |t|\ge 1/e,\\
      -\log(-\log(|t|))-\log(1+|t|^{2})/2& \text{ if } |t|\le 1/e,
    \end{cases}
  \end{displaymath}
  is $\omega $-psh. All the Lelong numbers of
  the function $\varphi$ are zero. For $t\not =0$ this is clear because the
  function $\varphi$ is locally bounded in $\bb P^{1}-\{(0:1)\}$. And at
  $t=0$ this follows from the fact that $\varphi$ grows at most as the
  logarithm of the logarithm. It follows that $P[\varphi]_{\mathcal{J}}=0$. Also
  the growth of the function $\varphi$ at $t=0$ shows that the current 
  $dd^{c}\varphi$ does not charge any pluripolar set. Therefore
  \begin{displaymath}
    \int_{\bb P^{1}}\langle \omega +dd^{c}\varphi \rangle =
    \int_{\bb P^{1}}\omega +dd^{c}\varphi =
    \int_{\bb P^{1}}\omega =
    \int_{\bb P^{1}}\langle \omega +dd^{c} P[\varphi]_{\mathcal{J}}\rangle=1.
  \end{displaymath}
  From \ref{thm:3} we may deduce that $P[\varphi]=0$.
\end{example}

\begin{example}\label{exm:6}
  The previous example can be generalized to the setting of good metrics in
  the sense of Mumford  \cite{hi}. Let $X$ be a smooth complex variety with $\dim X=n$ 
  and $D\subset X$ a normal crossings divisor. Let $\L$ be a
  line bundle on $X$ and $h$ a singular metric on $\L$. The metric $h$
  is said to be \emph{good} if $h$ is smooth on $X\setminus D$, and for every
  holomorphic chart $U$ of $X$ with
  coordinates $z_{1},\dots, z_{n}$ in which $D$ has the equation
  $z_{1}\cdots z_{k}=0$, each generating section $s$ of $\L$ on $U$ and
  each vector field $v$ on $U$ there is a neighborhood $V$ of
  $(0,\dots,0)$ in which the estimates    
  \begin{enumerate}
  \item  $h(s), h^{-1}(s)\le C \left(\sum_{i=1}^{k}-\log
      |z_{i}|\right)^{2m}$ for some $C\in \bb{R}_{>0}$ and $m\in \bb{N}$,
  \item $\|h(s)^{-1}\partial h(s) (v)\|^{2}\le C
    \sum_{i=1}^{k}\frac{1}{|z_{i}|^{2}(\log|z_{i}|)^{2}}$  for some $C\in
    \bb{R}_{>0}$
  \end{enumerate}
    hold.
  Good metrics appear naturally when considering toroidal
  compactifications of locally symmetric spaces, see \cite{hi}.  
  
  Assume that $X$ is projective and that $h$ is, at the
  same time, good 
  and  psh. Choose a smooth reference metric $h_{0}$ on $\L$ and write $\omega 
  =c_{1}(\L,h_{0})$  and $\varphi=-\log\left(h(s)/h_0(s)\right)$, so that  $\varphi$ is $\omega $-psh. 
    
  By \cite[Proposition~1.2]{hi} we have that the Lelong numbers of
  $\varphi$ are zero on all points of all modifications of $X$
  and that 
 \begin{displaymath}
   \int_{X}\langle (\omega +dd^{c}\varphi)^{n}\rangle =
   \int_{X\setminus D} (\omega +dd^{c}\varphi)^{n} =  \int_X \omega^n = \deg(\L) \, .
 \end{displaymath}
 We see in particular that $P[\varphi]_{\mathcal{J}}=0$.
 
 Now assume that $\L$ is moreover ample. Then $\deg(\L)>0$ and hence we have 
 $\int_{X}\langle (\omega +dd^{c}\varphi)^{n}\rangle>0$.
 Also, since $\ca J(k \varphi)$ is an ideal sheaf
 \begin{displaymath}
   \lim_{k\to \infty }\frac{h^{0}(X,\mathcal{L}^{\otimes k}\otimes
     \mathcal{J}(k\varphi))}{k^{n}/n!}\le
   \lim_{k\to \infty }\frac{h^{0}(X,\mathcal{L}^{\otimes
       k})}{k^{n}/n!} = \deg(\L) \, . 
 \end{displaymath}
 We may deduce from \ref{thm:3} that equality holds in the latter and that
 $P[\varphi]=P[\varphi]_{\ca{J}}=0$.   
\end{example}

\subsection{Multiplier ideal volume equals non-pluripolar volume in the case of algebraic singularities} \label{sec:volumes_alg_sing}

The purpose of this section is to show that multiplier ideal volume equals non-pluripolar volume in the case of algebraic singularities. 

Let $\theta $ be a closed smooth $(1,1)$-form on the pure-dimensional projective complex manifold $X$ and let
  $\varphi$ be a $\theta $-psh function on $X$. Let $n = \dim X$.
  \begin{lemma}\label{prop:alg-sing}
Assume that $\varphi$ has algebraic
  singularities as in \ref{def:3}. 
 Let $c \in \bb Q_{\ge 0}$ be the constant associated to $\varphi$. 
 Take a modification $\pi \colon
  X_{\pi }\to X$ such that $\pi^{-1}\mathcal{I}(\varphi/c)\cdot \ca O_{X_\pi}=\mathcal{O}_{X_\pi}(-D)$ for an effective simple normal crossings divisor $D$. Let $\class{\pi^*\theta }$ and
$[D]$ denote the cohomology classes of the closed smooth $(1,1)$-form
$\pi^*\theta$ and the divisor $D$ on $X_\pi$. Then
$\class{\pi^*\theta} - c[D]$ is a nef class on $X_\pi$ and
   the equality
\begin{align}\label{eq:alg}
   \left(\class{\pi^*\theta} -c[D]\right)^{n} =  \int_{X}\left\langle \left(\theta +dd^{c}\varphi\right)^{\wedge n}\right\rangle   
  \end{align}
holds in $\qq_{\ge 0}$.
 \end{lemma} 
\begin{proof} The Siu decomposition (see \ref{rem:siu-dec}) of $\pi^*\left(\theta + dd^{c}\varphi \right)$ on $X_{\pi}$ is
\begin{align}\label{eq:siu}
\pi^*(\theta + dd^{c}\varphi) = T' + c \, \delta_D \, , 
\end{align}
with  $T'$  a closed positive $(1,1)$-current with locally bounded potentials representing the class $\class{\pi^*\theta} - c[D]$.   
It follows from \ref{lemm:2} that $\class{\pi^*\theta} - c[D]$ is a nef class. Next we note that
\[ \int_X \langle \theta + dd^{c}\varphi \rangle^n = \int_{X_\pi} \langle T' \rangle^n  \]
as $T'$ agrees with $\theta + dd^{c}\varphi$ on $X_{\pi}\setminus D$. 
As $T'$ represents the cohomology class $\class{\pi^*\theta} - c[D]$ on $X_\pi$ and moreover is positive with locally bounded potentials we have by \ref{lemm:2} that
\[  \int_{X_\pi} \langle T' \rangle^n =  \left(\class{\pi^*\theta} -c[D]\right)^{n} \, . \]
This proves the required equality.    It is clear that the degree $\left(\class{\pi^*\theta} -c[D]\right)^{n} $ is a rational number.
\end{proof}

The next result is proved in \cite[Theorem 2.26]{DX}, under the hypothesis that $\ca L$ is ample with a K\"ahler metric. We give an alternative proof, and remove the hypothesis that $\ca L$ is ample.

\begin{theorem}\label{th:asym_alg}
  Let
  $\mathcal{L}$ be a line bundle on~$X$
  provided with a smooth reference metric $h_{0}$ and a psh metric
  $h$. Let $\theta =c_{1}(\mathcal{L},h_{0})$ be the first Chern form and
  $\varphi=-\log(h(s)/h_{0}(s))$ the resulting $\theta $-psh function as in
  \ref{rem:bij-metrics}. If $\varphi$ has algebraic singularities then the equality
  \begin{displaymath}
    \lim _{k\to \infty}\frac{h^{0}(X,\mathcal{L}^{\otimes k}\otimes
      \mathcal{J}(k\varphi))}{k^{n}/n!} =
    \int_{X}\langle (\theta +dd^{c}\varphi)^{\wedge n}\rangle 
  \end{displaymath}
  holds in $\qq_{\ge 0}$.
\end{theorem}
\begin{proof} Let $c \in \qq_{\ge 0}$ be the constant associated to $\varphi$ as in \ref{def:3}. Let $\pi \colon X_{\pi }\to X$ be
  a modification such that $\pi
  ^{-1}\mathcal{I}(\varphi/c)\cdot\ca O_{X_\pi}=\mathcal{O}(-D)$ for $D=\sum_{i}\alpha
  _{i}D_{i}$ an effective simple normal crossings divisor on $X_{\pi}$. By \ref{prop:alg-sing} we have that $\pi^*\L \otimes \mathcal{O}(-cD)$ is nef, and the 
equality
\begin{align}\label{eq-al}
\deg(\pi ^{\ast}\mathcal{L} \otimes \mathcal{O}(-c
      D)) = \int_{X}\langle (\theta +dd^{c}\varphi)^{\wedge n}\rangle 
\end{align}
holds in $\qq_{\ge 0}$.
 Let $\ell$ be a positive integer such that $\ell c\in \bb Z$. Denote by $\vol(\ca M)$ the volume of a line bundle $\ca M$. By
  \ref{lemm:10}
  \begin{align*}
        \lim _{k\to \infty}\frac{h^{0}(X,\mathcal{L}^{\otimes k}\otimes
    \mathcal{J}(k\varphi))}{k^{n}/n!}
    & =
      \lim _{k\to \infty}\frac{h^{0}(X,\mathcal{L}^{\otimes k\ell}\otimes
      \mathcal{I}((k\ell c\varphi)/c))}{(k\ell)^{n}/n!}\\
    &=\lim _{k\to \infty}
      \frac{h^{0}\big(X_{\pi },\pi ^{\ast}\mathcal{L}^{\otimes k\ell}\otimes
            \mathcal{O}(-\ell c D)^{\otimes k}\big)}{(k\ell)^{n}/n!}\\
    &=\frac{\vol(\pi ^{\ast}\mathcal{L}^{\otimes \ell}\otimes \mathcal{O}(-\ell c
      D))}{\ell ^{n}}.
  \end{align*}
  
  As $\pi ^{\ast}\mathcal{L}^{\otimes \ell}\otimes \mathcal{O}(-\ell c D)$ is nef, we have the Hilbert--Samuel formula
      \begin{displaymath}
        \vol(\pi ^{\ast}\mathcal{L}^{\otimes \ell}\otimes \mathcal{O}(-\ell c
      D))= \deg(\pi ^{\ast}\mathcal{L}^{\otimes \ell}\otimes \mathcal{O}(-\ell c
      D)) \, , 
      \end{displaymath}
   and the result follows by applying \ref{eq-al}.

\end{proof}

\section{Almost asymptotically algebraic singularities}\label{sec:psh-functions-that}

We continue to assume that $X$ is a pure-dimensional projective complex manifold.

\subsection{Definition and first examples}

Let $\theta $ be a closed smooth
  $(1,1)$-form and $\omega $ a K\"ahler form on~$X$. The following terminology has been introduced by Rashkovskii
\cite{Rash:13} in the local case of isolated singularities, but can
be adapted 
to the global case of $\theta $-psh functions. 
\begin{definition}\label{def:4} 
  A $\theta$-psh function $\varphi$ on $X$ is said to have \emph{asymptotically
  algebraic singularities with respect to $\theta$} if there is a sequence of quasi-psh
  functions $(\psi_{m})_{m\ge 1}$ with algebraic singularities and a sequence of positive real numbers
  $(a_{m})_{m\ge 1}$
  converging monotonically to zero such that for each $m \ge 1$ the function
  $\psi_{m}$ is $(\theta +a_{m}\omega )$-psh
  and the inequalities
  \begin{displaymath}
    (1+\frac{1}{m})\psi _{m} \prec \varphi \prec
    (1-\frac{1}{m})\psi _{m}
  \end{displaymath}
  hold.
\end{definition}
We will omit the addition ``with respect to $\theta$'' from the terminology if the form $\theta$ is clear from the context. 

For instance, if $\varphi$ has isolated singularities, and is tame in
the sense of \cite{bfj-val} or is exponentially H\"older, then it has
asymptotically algebraic singularities, see \cite[Examples 3.6 and
3.7]{Rash:13}.

We introduce the following slightly weaker notion. 
\begin{definition}\label{def:6}
  A $\theta $-psh function $\varphi$ on $X$ is said to have \emph{almost asymptotically algebraic singularities with respect to $\theta$} if there exists a quasi-psh function $f$ with algebraic singularities,
 a sequence of quasi-psh
  functions $(\psi_{m})_{m\ge 1}$ with algebraic singularities and a sequence of positive real numbers
  $(a_{m})_{m\ge 1}$
  converging monotonically to zero such that for all $m \ge 1$ the function
  $\psi_{m}$ is $(\theta +a_{m}\omega )$-psh and the inequalities
  \begin{displaymath}
    \psi _{m} +\frac{1}{m}  f\prec  \varphi \prec
    \psi _{m }
  \end{displaymath}
  hold.
  A  psh metric $h$ on a line bundle $\mathcal{L}$ on $X$ such that $\theta$ represents the cohomology class $c_1(\ca L)$  is
  said to have almost asymptotically algebraic singularities if the
  corresponding $\theta $-psh function has almost asymptotically
  algebraic singularities with respect to $\theta$.

\end{definition}
Again, we will omit the addition ``with respect to $\theta$'' from the terminology if the form $\theta$ is clear from the context.

\begin{lemma} \label{ex:good_implies_aaas} If $h$ is a psh good metric (in the sense of \ref{exm:6}) on a line bundle $\ca L$, and $h_0$ is a smooth metric on $\ca L$, then $h$ has almost asymptotically algebraic singularities with respect to $\theta \coloneqq c_{1}(\mathcal{L},h_{0})$
  \end{lemma}

  \begin{proof}
Let $D$ be a divisor as in
  \ref{exm:6}, $s$ a non-zero rational section of $\mathcal{L}$, and set
  \begin{displaymath}
    \varphi=-\log(h(s)/h_{0}(s)).
  \end{displaymath}
  This is a $\theta $-psh function, with singularities contained in $D$.
  Choose an effective normal crossings divisor $A$ such that $|D|\subset |A|$ and such that $\ca O(A)$ admits a smooth psh metric (for example, we could take $A$ ample); choose one such metric. Let $1$ be the canonical section of $\mathcal{O}(A)$ (so
  $\divisor(1)=A$), and write $f=\log\|1\|$. We claim that for every $m>0$ the inequalities
  \begin{displaymath}
    \frac{1}{m}f\prec \varphi \prec 0
  \end{displaymath}
  hold. That $\varphi \prec 0$ follows from the assumption that $h$ be psh; we will use goodness to establish the other inequality. We work locally around a point $x \in X$, where we can assume that our rational section $s$ is in fact a generating section of $\ca L$, so that $\varphi \sim - \log h(s)$. We write $z_1, \dots, z_a$ for local defining equations for the branches of $A$ through $x$, and we order them so that $D$ is cut out by $\prod_{i=1}^b z_i$ for some $b \le a$. Then 
  \begin{equation}
  f \sim \sum_{i=1}^a a_i \log\abs{z_i}
  \end{equation}
  where the $a_i \in \bb Z_{>0}$ are the multiplicities of the branches cut out by $z_i$ in the divisor $A$. Now from the definition of a good metric we see that 
  \begin{equation}
  h(s) \le C \left(\sum_{i=1}^b - \log\abs{z_i}\right)^{2M}
  \end{equation}
  for some positive integer $M$ and positive real number $C$. Hence 
  \begin{equation}
  -2M \log\left(\sum_{i=1}^b - \log\abs{z_i}\right) \prec \varphi, 
  \end{equation}
  from which we see that 
  \begin{equation}
  \frac{1}{m} f \prec \varphi
  \end{equation}
  for all $m \ge 1$. 
  \end{proof}

\begin{remark} \label{tensor_toroidal}  \label{has_unbounded}
\begin{enumerate}
\item The quasi-psh functions with almost asymptotically algebraic singularities form a convex cone. More precisely, if $\varphi_1$ is a $\theta_1$-psh function with almost asymptotically algebraic singularities with respect to $\theta_1$, and $\varphi_2$ is a $\theta_2$-psh function with almost asymptotically algebraic singularities with respect to $\theta_2$, then $\varphi_1+\varphi_2$ is a $(\theta_1+\theta_2)$-psh function with almost asymptotically algebraic singularities with respect to $\theta_1+\theta_2$.

\item The notion of (almost) asymptotically algebraic singularities is birationally invariant. Hence all of the results concerning almost asymptotically algebraic singularities continue to hold if we pass to a modification of $X$.
\item  Any psh metric with almost asymptotically algebraic singularities  has
  small unbounded locus. Indeed, if $\psi + f \prec \varphi \prec \psi$ with $\psi, f$ quasi-psh with algebraic singularities, then $\varphi$ is locally bounded away from the singular loci of $\psi$ and $f$, which are proper Zariski closed, and in particular pluripolar, sets.

\end{enumerate}
\end{remark}

\begin{lemma} \label{lem:independent_Kahler} The notion of almost asymptotically algebraic
  singularities does not depend on the choice of K\"ahler metric
  $\omega $. Moreover we can choose the function $f$ in \ref{def:6} to be
  $\omega$-psh and the sequence  $(a_{m})$ to be  $a_m =
  \frac{1}{m}$. 
\end{lemma}
\begin{proof}
  Assume that $\varphi$ has almost asymptotically algebraic
  singularities with respect to $\theta$. Let $\omega '$ be another K\"ahler metric on $X$. Then by the compactness of
  $X$ there is a real number $a>0$ such that $\omega \le a\omega
  '$. Therefore, if $\psi _{m}$ is $(\theta +a_{m}\omega )$-psh, then
  it is also $(\theta +a_{m}a\omega ')$-psh. So after setting
  $(a_{m} a)$ instead of $(a_{m})$ we see that $\varphi$ satisfies
  \ref{def:6} for $\omega '$.

  Since $f$ is quasi-psh, by \ref{lemm:8} there is a real number $b>0$
  such that $f$ is $(b\omega)$-psh and hence  $f/b$ is $\omega $-psh.
  Choose an increasing sequence $m_{k}$ of integers satisfying
  \begin{equation}\label{eq:5}
    a_{m_{k}}\le \frac{1}{k}, \qquad \text{ and } \qquad m_{k}>kb.
  \end{equation}
  This can easily be achieved as $(a_{m})$ converges to zero.
  Writing $\psi '_{k}=\psi _{m_{k}}$ and $f'=f/b$ we have that $\psi
  '_{k}$ is $(\theta +a_{m_{k}}\omega) $-psh. By the first
  condition in \ref{eq:5}, the function $\psi'_k$ is $(\theta +(1/k)\omega) $-psh. Using
  the second condition in \ref{eq:5} we find
  \begin{displaymath}
    \psi '_{k}+\frac{1}{k}f' =
    \psi _{m_{k}}+\frac{1}{bk}f \prec
    \psi _{m_{k}}+\frac{1}{m_{k}}f\prec \varphi.
  \end{displaymath}
 It follows that the functions $\psi '_{k}$ and $f'$ satisfy the conditions of
  \ref{def:6} where $f'$ is  $\omega $-psh  and  where $a_{k}=1/k$. 
\end{proof}

\begin{lemma}\label{ex:aa}
If $\varphi$ has asymptotic algebraic singularities with respect to $\theta $, then it has almost asymptotically algebraic singularities with respect to $\theta$. 
    \end{lemma}
    \begin{proof}
  Let $\psi _{m}$ be a sequence of functions satisfying
  \ref{def:4}. For $m\ge 2$ we have the chain of inequalities
  \begin{displaymath}
    \frac{3}{2}\psi _{2}\prec \varphi \prec \left (
      1-\frac{1}{m}\right)\psi _{m}\prec \frac{1}{2}\psi _{m},
  \end{displaymath}
  which implies $\psi _{m}\succ 3\psi _{2}$. Therefore we have the
  chain of inequalities
  \begin{displaymath}
    \left(1-\frac{1}{m}\right)\psi _{m}+\frac{1}{m}6\psi _{2}
    \prec
    \left(1-\frac{1}{m}\right)\psi _{m}+\frac{2}{m}\psi _{m}
    \prec \varphi \prec \left(1-\frac{1}{m}\right)\psi _{m} \, , 
  \end{displaymath}
showing that $\varphi$ has almost asymptotically algebraic singularities. 
\end{proof}

    The converse is not true. We will illustrate this with the
    function $\varphi$ of \ref{exm:7}. We use the notations in
    that example. Note that  the function
    \begin{displaymath}
      f(t)=\log|t|-\log(1+|t|^{2})/2
    \end{displaymath}
 is $\omega $-psh. For all $m \in \zz_{>0}$ we have
    \begin{displaymath}
      \frac{1}{m} f \prec \varphi \prec 0.
    \end{displaymath}
    Indeed the only point where $f$ and $\varphi$ are not locally
    bounded is the point $t=0$. Close to this point $\varphi$ has a
    singularity of the shape $-\log(-\log|t|)$, while $f/m$ has a
    singularity of the shape $\log(|t|)/m$ which is more singular for
    all values of $m$. 
    So, taking $\psi _{m}=0$ in \ref{def:6}, we see that
    $\varphi$ has almost asymptotically algebraic 
    singularities. Nevertheless $\varphi$ does not have asymptotically
    algebraic singularities. Assume that it satisfies \ref{def:4} for a family of functions $\psi _{m}$ with algebraic
    singularities. Then for $m=2$ we have
    \begin{displaymath}
      \frac{3}{2}\psi _{2}\prec \varphi \prec \frac{1}{2} \psi _{2}
    \end{displaymath}
    where $\psi_2$ has
    algebraic singularities.
    Since the Lelong numbers of $\varphi$ are zero, the right
    inequality implies that $0\prec
    \psi _{2}$. But this contradicts the left inequality as $0\not
    \prec \varphi$.

 \subsection{Some criteria for almost asymptotically algebraic singularities} \label{sec:criteria}
 
 The purpose of this section is to exhibit some useful criteria that allow to verify that a quasi-psh function
 has almost asymptotically algebraic singularities. Our results are based on Demailly's regularization theorem (\ref{thm:2}). We continue
 with our assumption that $X$ is a projective complex manifold. We fix a background K\"ahler form $\omega$ on $X$.
 
 \begin{definition} \label{def:1} Let $U$ be a Euclidean open subset of $X$ with norm $\|\cdot\|$.
  An upper semi-continuous function $\varphi \colon U \to \bb R\cup
  \{-\infty\}$ is said to be
  \emph{meromorphically 
  Lipschitz} if there exists a finite open coordinate covering $\{U_{i}\}$
  of $U$ and for each open $U_{i}$ there is a regular algebraic function $f_{i}$ on $U_i$ 
  such that
  \begin{displaymath}
    \varphi(x)-\varphi(y)\le\frac{\| x-y\|}{|f_{i}(y)|} \, , \quad x, y \in U_i \, . 
  \end{displaymath}
\end{definition}

\begin{lemma}\label{lemm:9}
  Assume we are given a finite open covering
  $\{U_{i}\}$ of $X$ and for each $i$ a regular algebraic function $f_{i}$ on
  $U_{i}$. For each $i$ let $V_{i}\subset \subset U_{i}$ be a relatively
  compact open subset such that the collection $\{V_{i}\}$ is still an open
  cover of $X$.  Then there exists a
  quasi-psh function $\varphi$  on $X$ with algebraic singularities such that,
  for all $i$, the inequality
  \begin{equation}\label{eq:14}
    \varphi|_{V_{i}}\le \log |f_{i}|
  \end{equation}
  holds.
\end{lemma}
\begin{proof}
  For each $i$ let $D_{i}$ be the divisor of $f_{i}$ on $U_{i}$. Let
  $E$ be an effective divisor on $X$ with
  \begin{equation}
    \label{eq:12}
    E|_{U_{i}}\ge D_{i} \quad \text{for all} \, i \, . 
  \end{equation}
  Choose a smooth hermitian metric $h_0$ on $\mathcal{O}(E)$ and let $s$ be
  a global section of $\mathcal{O}(E)$ with $\divisor(s)=E$. By
  condition \ref{eq:12}, the function $s/f_{i}$ is regular on $U_{i}$. Therefore
  $h_0(s)/\abs{f_i}$ is continuous on $U_{i}$, hence
  $\log(h_0(s))-\log|f_{i}|$ is bounded above on
  the compact subset $\overline {V_{i}}$. Let $M$ be a real number such that for all
  $i$ we have the bound
  $\log(h_0(s))-\log|f_{i}|\le M$ on  $\overline {V_{i}}$. Let $\theta
  =c_{1}(\mathcal{O}(E),h_0)$ be the first Chern form of
  $\mathcal{O}(E)$ with smooth metric $h_0$. Then the function $\varphi=\log(h_0(s))-M$ is
  $\theta $-psh, has algebraic singularities  and satisfies the
  inequalities \ref{eq:14}.
\end{proof}

Let $\theta$ be a smooth $(1,1)$-form on $X$.
\begin{proposition} \label{prop:4}
  Let $\varphi$ be a $\theta $-psh function on $X$ that can be written locally as a
  sum $\varphi=\phi +\gamma $ with $\phi $ meromorphically
  Lipschitz and $\gamma $ bounded. Then $\varphi$ has almost asymptotically algebraic singularities with respect to $\theta$. 
\end{proposition}
\begin{proof} 
Let $n=\dim X$, and choose finite
  open coordinate coverings  $\{U_{i}\}$ and $\{V_{i}\}$ of $X$ with
  $V_{i}\subset \subset U_{i}$. 
  By Demailly's regularization theorem (\ref{thm:2}) there exist constants
  $C_{1}$ and $C_{2}$, a sequence of functions $\varphi_{m}$ with algebraic singularities on $X$ satisfying, in each $V_{i}$,
  \begin{equation}
    \label{eq:11}
    \varphi(z) -\frac{C_{1}}{m}\le \varphi_{m}(z)\le
    \sup_{\|x-z\|<r}\varphi (x)+\frac{1}{m}\log\frac{C_{2}}{r^{n}}
  \end{equation}
  and a sequence of positive real numbers $(a_{m})_{m> 0}$ converging
  monotonically to zero such that each function $\varphi_{m}$ is $(\theta +a_{m}\omega) $-psh.
  Note that $\varphi \prec \varphi_m$.

  By our assumption on $\varphi$, after taking a finite refinement of the open cover we can
  assume that there exist functions $f_{i}$ that are regular on
  $U_{i}$ and such that on each $U_{i}$ the
  estimate
  \begin{displaymath}
    \varphi(x)-\varphi(z)\le \frac{\|x-z\|}{|f_{i}(z)|}+ C_{3}
  \end{displaymath}
  holds for a constant $C_{3}$. By \ref{lemm:9} there is a
  quasi-psh function $\psi $ on $X$ with algebraic singularities and
  such that $\psi |_{V_{i}}\le \log|f_{i}|$ for all $i$.
  
  Taking now $r=|f_{i}(z)|$, we deduce from~\ref{eq:11}
  that 
  \begin{align*}
    \varphi_{m}(z)&\le
    \varphi (z)+
    \frac{r}{|f_{i}(z)|}+C_{3}+\frac{1}{m}\log\frac{C_{2}}{r^{n}}\\
    &\le \varphi (z)+ 1+C_{3}+\frac{1}{m}\log C_{2}-\frac{n}{m}\log
      |f_{i}(z)|\\
    &\le \varphi(z) + C_{
    4}- \frac{n}{m}\psi(z) .
  \end{align*}
  The function
  \begin{math}
    f(z) = n\psi (z)
  \end{math}
  is quasi-psh with algebraic singularities and we have found that for
  every $m$ the estimate
  \begin{displaymath}
    \varphi_{m}+ \frac{1}{m}f \prec \varphi \prec \varphi_{m} 
  \end{displaymath}
  holds. This concludes the proof of the proposition.
\end{proof}
We can use \ref{prop:4} to see that toroidal singularities are almost asymptotically algebraic.
The following definition is inspired by the fact that if $g \colon \rr_{>0}^k \to \rr$ is a bounded-above convex function, then the function
\[ g(-\log|z_1|, \ldots, -\log|z_k|) \]
is a psh function on $D(1)^k$. Here $D(1)$ is the open unit disk $\{ z \in \cc \, | \, |z| <1 \}$.

\begin{definition}  \label{def:toroidal_sing}
 Let $U\subset X$ be Zariski open with $D=X\setminus U$ a
 normal crossings divisor. A quasi-psh function $\varphi$ on $X$ is
 said to have \emph{toroidal singularities} (with respect to~$D$) if $\varphi$ is locally bounded on $U$
 and there exists an open coordinate covering $\{V_{i}\}$ of $X$ such that on
 each $V_{i}$ the divisor $D$ has equation $z_{1}\cdots
 z_{k_{i}}=0$ and the restriction $\varphi|_{V_{i}\cap U}$ can be written as
 \begin{equation} \label{def:toroidal}
   \varphi|_{V_{i}\cap U}= \gamma +g(-\log|z_{1}|,\dots ,-\log|z_{k_{i}}|) \, ,
 \end{equation}
 where $\gamma $ is locally bounded on $V_{i}$ and $g$ is a
 bounded above convex Lipschitz continuous function defined on a quadrant
 $\{u_{j}\ge M_{j}| j=1,\dots,k_{i}\}\sub \bb R^{k_i}$.
 We say that a  psh metric on a line bundle $\mathcal{L}$ on $X$  has toroidal singularities if any
  corresponding $\theta $-psh function has toroidal singularities. Here  $\theta$ is any smooth closed differential form representing the cohomology class $c_1(\ca L)$.  Note that this is equivalent to any local psh potentials being of the form \ref{def:toroidal}.
\end{definition}

\begin{proposition} \label{prop:toroidal_implies_aaas}
  A toroidal $\theta$-psh function on $X$ has almost asymptotically algebraic
  singularities with respect to $\theta$.  
\end{proposition}
\begin{proof} By \ref{prop:4}
  it suffices to check that, if $g$ is a bounded above convex Lipschitz
  continuous function defined on a quadrant 
 $K=\{u_{j}\ge M_{j}| j=1,\dots,k\}$, then 
  $g(-\log|z_{1}|,\dots ,-\log|z_{k}|)$ is meromorphically
  Lipschitz in the polydisk $|z_{j}|< e^{-M_{j}}$.
  We first observe that, if $s,t>0$ are real numbers, then
  \begin{equation}
    \label{eq:15}
    \log(s)-\log(t)\le \frac{|s-t|}{t}.
  \end{equation}
  Indeed, if $s\le t$ then the right hand side is greater or equal to
  zero while the left hand side is smaller or equal than zero. And
  if $t< s$ then
  \begin{displaymath}
    \log(s)-\log(t)=\int _{t}^{s}\frac{d\xi}{\xi}\le \frac{s-t}{t}.
  \end{displaymath}
  Next we see that since $g$ is bounded above convex Lipschitz continuous on the quadrant $K$ the function
  $g$ is non-increasing on each semi-line of the form $u+\lambda
  w$, $\lambda \ge 0$, 
  for $u\in K$ and $w$ a  vector with non-negative entries.
  Therefore, if  
  $u=(u_{1},\dots,u_{k})$ and $v=(v_{1},\dots,v_{k})$ are points
  in $K$ then
  \begin{equation}\label{eq:16}
    g(u) - g(v) \le  \sum_{j\colon v_{j}>u_{j}} C(v_{j}-u_{j}),  
  \end{equation}
  where $C$ is the Lipschitz constant of $g$.
  Now using equations \ref{eq:15} and \ref{eq:16} we obtain for $z, x \in \cc^k$ with $|z_j|, |x_j| < e^{-M_j}$ for $j=1,\ldots,k$ that
  \begin{align*}
   & g(-\log|z_{1}|,\dots ,-\log|z_{k}|)
                          -g(-\log|x_{1}|,\dots
    ,-\log|x_{k}|)\\
    &\le \sum _{j\colon |z_{j}|>|x_{j}|} C(\log|z_{j}|-\log|x_{j}|)\\
    &\le \sum _{j\colon |z_{j}|>|x_{j}|}
      C\frac{|x_{j}-z_{j}|}{|x_{j}|}\\
    &\le C'\frac{\|x-z\|}{\prod |x_{j}|}
  \end{align*}
  for some constant $C'$, proving the claim.
\end{proof}

\subsection{Multiplier ideal volume equals non-pluripolar volume in the case of almost asymptotically algebraic singularities}

We continue to assume that $X$ is a projective pure-dimensional complex manifold. Let $n=\dim X$.

The main result of this section is that for a quasi-psh function with almost asymptotically algebraic singularities the 
multiplier ideal volume and the non-pluripolar volume are equal. See \ref{th:asym}. 

\begin{lemma} \label{lemm:11}
  Let $D\subset X$ be a smooth
  divisor, $A$ an ample divisor and $\ell \in \zz_{>0}$. We denote by
  $\mathcal{O}_{\ell D}$ the coherent sheaf on $X$ defined by the exact
  sequence
  \begin{displaymath}
    0\longrightarrow \mathcal{O}_{X}(-\ell D) \longrightarrow \mathcal{O}_{X}\longrightarrow
    \mathcal{O}_{\ell D}\to 0 \, .
  \end{displaymath}
  Then for all $k \in \zz_{>0}$ we have
  \begin{displaymath}
    h^{0}(X,\mathcal{O}_{\ell D}(kA))\le \ell
    k^{n-1}D\cdot A^{n-1}. 
  \end{displaymath}
\end{lemma}
\begin{proof}
  For every  $0\le j\le \ell-1$ there is an exact sequence
  \begin{displaymath}
    0\longrightarrow \frac{\mathcal{O}(-(j+1) D)}{\mathcal{O}(-\ell D)}
    \longrightarrow \frac{\mathcal{O}(-j D)}{\mathcal{O}(-\ell D)}
    \longrightarrow \mathcal{O}_{D}(-jD)\longrightarrow 0.
  \end{displaymath}
  Adding up, these exact sequences imply that
  \begin{align*}
     & h^{0}(X,\mathcal{O}_{\ell D}(kA)) \le
    \sum _{j=0}^{\ell-1}h^{0}(X,\mathcal{O}_{D}(kA-jD))\\
    & \le \sum _{j=0}^{\ell-1}h^{0}(X,\mathcal{O}_{D}(kA)))
    \le \ell k^{n-1}D\cdot A^{n-1}. \qedhere
  \end{align*}

\end{proof}
  Let
  $\mathcal{L}$ be a line bundle on~$X$
  provided with a smooth reference metric $h_{0}$ and a psh metric
  $h$. Let $\theta =c_{1}(\mathcal{L},h_{0})$ be the first Chern form and
  $\varphi=-\log(h(s)/h_{0}(s))$ the associated $\theta $-psh function. As before we let
  \[ \vol_{\ca J}(\ca L,h) =  \lim _{k\to \infty}\frac{h^{0}(X,\mathcal{L}^{\otimes k}\otimes
      \mathcal{J}(k\varphi))}{k^{n}/n!}  \]
be the multiplier ideal volume of $(\ca L,h)$. The following can be seen as a Hilbert--Samuel type formula. 
\begin{theorem}\label{th:asym}
 If $\varphi$ has almost
  asymptotically algebraic singularities with respect to $\theta$ then the equality
  \begin{equation}\label{eq:lin_eq_int}
     \vol_{\ca J}(\ca L,h)  =
    \int_{X}\langle (\theta +dd^{c}\varphi)^{\wedge n}\rangle 
  \end{equation}
  holds in $\rr_{\ge 0}$.
\end{theorem}
\begin{proof}

    By \ref{lem:independent_Kahler} we can find an ample line bundle
    $\mathcal{O}(D)$ on $X$ with a smooth metric such that its
    first Chern form is  a K\"ahler form $\omega $, an $\omega $-psh
    function $f$ with algebraic singularities, and a sequence of
    quasi-psh functions $\varphi_{m}$ with algebraic singularities
    such that $\varphi_{m}$ is $(\theta +(1/m)\omega )$-psh and   
    \begin{equation}\label{eq:almost-asy}
    \varphi_m + \frac{1}{m}f \prec \varphi \prec \varphi_m \, .
  \end{equation}
  Also, for each $m>0$ each of the three functions separated by the two inequalities in \ref{eq:almost-asy} is $(\theta +(2/m)\omega )$-psh by \ref{lemm:8}. We recall from  \ref{has_unbounded} that almost asymptotically  algebraic singularities have small unbounded locus. Thus we can apply the monotonicity property of the
    non-pluripolar product (\ref{thm:1}) and conclude that for each fixed $m>0$ that
    \begin{align}\label{eq:19}
      &\int _{X}\langle (\theta +\frac{2\omega
      }{m}+dd^{c}\varphi_{m}+\frac{1}{m}dd^{c}f)^{n}\rangle
      \le \int _{X}\langle (\theta +\frac{2\omega
      }{m}+dd^{c}\varphi)^{n}\rangle \nonumber \\
      &\le
      \int _{X}\langle (\theta +\frac{2\omega
      }{m}+dd^{c}\varphi_{m})^{n}\rangle. 
    \end{align}
    Now set
    \[ A_{m}\coloneqq \int _{X}\langle (\theta +\frac{2\omega
      }{m}+dd^{c}\varphi_{m}+\frac{1}{m}dd^{c}f)^{n}\rangle
      - \int _{X}\langle (\theta +\frac{\omega
      }{m}+dd^{c}\varphi_{m})^{n}\rangle \, . \]
    By the multi-additivity of the non-pluripolar product (\ref{prop:6}) we have 
    \begin{align*}
A_{m}&=\sum_{k=1}^{n}\binom{n}{k}\int_{X} \langle (\theta +\frac{\omega
      }{m}+dd^{c}\varphi_{m})^{n-k}(\frac{\omega
        }{m}+\frac{1}{m}dd^{c}f)^{k}\rangle\\
       &=\frac{1}{m}\sum_{k=1}^{n}\frac{1}{m^{k-1}}\binom{n}{k}\int_{X}
         \langle (\theta +\frac{\omega
      }{m}+dd^{c}\varphi_{m})^{n-k}(\omega
        +dd^{c}f)^{k}\rangle \, . 
    \end{align*}
    Note that $A_m \ge 0$ since all the summands on the right of the above expression are positive. 
    We see that there is a constant $C$ such that $0\le A_{m}\le C/m$ for all $m$
    so 
    \begin{displaymath}
      \lim_{m\to \infty} A_{m}=0.
    \end{displaymath}
    Similarly,
    \begin{align*}
      &\lim_{m\to \infty }
      \int _{X}\langle (\theta +\frac{2\omega
      }{m}+dd^{c}\varphi_{m})^{n}\rangle-\int _{X}\langle (\theta +\frac{\omega
        }{m}+dd^{c}\varphi_{m})^{n}\rangle=0,\\
      &\lim_{m\to \infty }
      \int _{X}\langle (\theta +\frac{2\omega
      }{m}+dd^{c}\varphi)^{n}\rangle-\int _{X}\langle (\theta +dd^{c}\varphi)^{n}\rangle=0.
    \end{align*}
    We conclude that when $m \to \infty$ each of the three terms
    in the inequality \ref{eq:19} converges to  
    \begin{displaymath}
    \int_X \langle (\theta + dd^c \varphi)^n\rangle .
  \end{displaymath}
      Next, for each $m\ge 0$ and each $k>0$ with $m\vert k$, by the monotonicity
      of multiplier ideals with the type of singularity,  we obtain the
      inequalities 
      \begin{equation}\label{eq:20}
      \begin{aligned}
        h^{0}\big(X,\mathcal{L}^{k}\otimes \mathcal{O}(2(k/m)D)\otimes
        & \mathcal{J}(k(\varphi_{m}+f/m)) \big) \\
        & \le
        h^{0}\left(X,\mathcal{L}^{k}\otimes \mathcal{O}(2(k/m)D)\otimes
        \mathcal{J}(k\varphi) \right)\\
        &\le 
        h^{0}\left(X,\mathcal{L}^{k}\otimes \mathcal{O}(2(k/m)D)\otimes
        \mathcal{J}(k(\varphi_{m}) \right).
    \end{aligned}
      \end{equation}
    Since $\varphi_{m}$ and $f$ have algebraic singularities, by \ref{th:asym_alg} we have the statement of the theorem, that is, the equality in  \ref{eq:lin_eq_int},
    for the terms at the left and at the right of the chain of inequalities \ref{eq:20}. Using the
    convergence to $\int_X \langle (\theta + dd^c \varphi)^n\rangle$ of all three terms in the chain \ref{eq:19} we deduce that
    \begin{displaymath}
      \lim_{m\to \infty}\lim _{\substack{k\to \infty\\m\vert k}}
      \frac{h^{0}\left(X,\mathcal{L}^{k}\otimes \mathcal{O}(2(k/m)D)\otimes
          \mathcal{J}(k\varphi) \right)}{k^{n}/n!}
      =\int_X \langle (\theta + dd^c \varphi)^n\rangle.
    \end{displaymath}
    It remains to show that
    \begin{equation}\label{eq:21}
      \lim_{m\to \infty}\lim _{\substack{k\to \infty\\m\vert k}}
      \frac{h^{0}\left(X,\mathcal{L}^{k}\otimes \mathcal{O}(2(k/m)D)\otimes
          \mathcal{J}(k\varphi) \right)}{k^{n}/n!}=
      \lim _{k\to \infty}
      \frac{h^{0}\left(X,\mathcal{L}^{k}\otimes
          \mathcal{J}(k\varphi) \right)}{k^{n}/n!}.
    \end{equation}
    
        Fix $k=\ell m$. 
    Replacing $D$ by a positive multiple we may assume that $D$ is effective and that $\ca L(D)$ is ample. 
    Then the exact sequence 
    \begin{displaymath}
     0\to \mathcal{O}(-2\ell D)\to \mathcal{O}\to \mathcal{O}_{2\ell D}\to 0 
   \end{displaymath}
 implies that
   \begin{align*}
     h^{0}\big(X,\mathcal{L}^{\ell m}\otimes \mathcal{O}(2\ell D)\otimes
      & \mathcal{J}(\ell m\varphi) \big) -
     h^{0}\big(X,\mathcal{L}^{\ell m}\otimes 
       \mathcal{J}(\ell m\varphi) \big) \\
    & \le h^{0}\left(X,\mathcal{L}^{\ell m}\otimes \mathcal{O}_{2\ell D}(2\ell D)\otimes
       \mathcal{J}(\ell m\varphi) \right) \\
   &  \le h^{0}\left(X,\mathcal{L}^{\ell m}\otimes \mathcal{O}_{2\ell D}(2\ell D)\right),
   \end{align*}

where the last inequality follows from the fact that $\mathcal{J}(\ell
m\varphi)$ is an ideal sheaf.

By Bertini we can further assume that $D$ is smooth. Then writing $\ca L = \ca O_X(E)$ for a divisor $E$ we calculate
\begin{equation}
\begin{split}
h^0\left(X,\ca O_{2 \ell D}(\ell m E + 2 \ell D)\right) & \le h^0\left(X,\ca O_{2 \ell D}(\ell m (E + D))\right) \\
& \le 2\ell(\ell m)^{n-1} D \cdot (E + D)^{n-1}, 
\end{split}
\end{equation}
where the last inequality is an application of \ref{lemm:11} with $A = E + D$.  In particular, 
\begin{equation}
\frac{h^0\left(X,\ca O_{2 \ell D}(\ell m E + 2 \ell D)\right) }{(\ell m)^n/n!} \to 0
\end{equation}
as $k = \ell m \to \infty$. This implies equation \ref{eq:21}.

\end{proof}
Combining \ref{thm:3} and \ref{th:asym} we obtain the following corollary.
\begin{corollary}\label{cor:env}
Let the situation be as in  \ref{th:asym}. Suppose that moreover we have $\int_X\langle ( \theta + dd^c\varphi)^n\rangle >0$, and that $\ca L$ is ample.  Then $P[\varphi] = P[\varphi]_{\mathcal{J}}$, where $P[\varphi]$ and $P[\varphi]_{\mathcal{J}}$ are the envelopes of singularity type from \ref{def:env}.
\end{corollary}
We do not expect that conversely, the equality $P[\varphi] = P[\varphi]_{\mathcal{J}}$ for a quasi-psh function $\varphi$ implies that $\varphi$ has almost asymptotically algebraic singularities.

\begin{remark} \label{upper_bound_by_multiplier_volume}
Let $h$ be a psh metric on $\ca L$ with small unbounded locus. We do not assume that $h$ has almost asymptotically algebraic singularities. Then we still have the inequality
  \[
  \vol_{\ca J}(\ca L,h) \geq \int_X \langle c_1(\ca L, h)^n \rangle.
      \]
 Indeed, consider a Demailly approximation sequence $\{\varphi_m\}_{m \in \bb N}$ for $\varphi$. Then $\varphi \prec \varphi_m$ and in the same way as in the proof of \ref{th:asym}, the inequality follows. This recovers the inequality in the result of Darvas and Xia in \ref{thm:3}, without the assumption that $\ca L$ is ample.
\end{remark}
\begin{remark}\label{rem:mult-vol-cont}
It follows from the monotonicity properties of the non-pluripolar product and \ref{th:asym} that the multiplier ideal volume has the following continuity property on the space of psh metrics with almost asymptotically algebraic singularities: let $\varphi$ be a quasi-psh function with almost asymptotically algebraic singularities with respect to $\theta$. Let $\varphi_m$ an approximating sequence of quasi-psh functions satisfying the conditions of \ref{def:6}. Then 
\[
\lim_{m \to \infty} \vol_{\ca J}(\L, \varphi_m) = \vol_{\ca J}(\L, \varphi).
\]

\end{remark} 

\section{b-divisors}\label{sec:b-div}
In this section we discuss Weil and Cartier $\rr$-b-divisors on compact algebraic complex manifolds. This is essentially Shokurov's notion of birational divisiors, or b-divisors, see  \cite{sh-prel}. For more background concerning b-divisors we refer to \cite{bff} and \cite{bfj-val}; see also \cite{BoteroBurgos} and
\cite{bo} for a discussion of the toroidal and the toric cases, respectively. 
\subsection{Basic definitions}
Throughout this section $X$ is a compact algebraic complex manifold (this section is purely algebraic, so if preferred the reader can work with finite-type algebraic varieties over any field of characteristic zero). We write $\Div\left(X\right)$ for the set of Weil divisors on $X$ with real coefficients, viewed as a real vector space (generally of infinite dimension). We endow it with
the direct limit topology with respect to its finite dimensional
subspaces. Explicitly, a sequence of divisors $(D_{i})_{i\ge
  0}$ converges to a divisor $D$ in $\Div\left(X\right)$ if there is a divisor $A$,
such that $\on{supp} (D_{i})\subset \on{supp}(A)$ for all $i\ge 0$ and
$(D_{i})_{i\ge 0}$ converges to $D$ in the finite dimensional vector
space of real divisors with support contained in $\on{supp}(A)$.

In \ref{def:modification} we defined a modification of complex manifolds. We note that if $\pi\colon X' \to X$ is a modification then $X'$ is also a compact algebraic complex manifold. 
\begin{definition}\label{def:modi}
The set of \emph{models of $X$} is
\[
R(X) \coloneqq \left\{\pi \colon X_{\pi}\to X \; \big{|} \; \pi \text{ is a modification}\right\}.
\]
We view $R(X)$ as a full subcategory of the category of complex manifolds over $X$, in particular morphisms are over $X$. Maps of models are unique if they exist, and are necessarily proper and bimeromorphic. 
\end{definition}

Hironaka's resolution of singularities implies that $R(X)$ is a directed set, where we set $\pi' \geq \pi$ if there exists a morphism $\mu \colon X_{\pi'} \to X_{\pi}$.

Consider a pair $\pi' \geq \pi$ in $R(X)$, and let $\mu \colon X_{\pi'}\to X_{\pi}$
be the corresponding modification. We have a pullback map
\[
\mu^* \colon \Div\left(X_{\pi}\right) \longrightarrow \Div\left(X_{\pi'}\right)
\]
and a pushforward map 
\[
\mu_* \colon \Div\left(X_{\pi'}\right) \longrightarrow \Div\left(X_{\pi}\right)
\]
of divisors. Both  maps are continuous.  

\begin{definition}\label{b-divisor}
The group of \emph{Cartier $\rr$-b-divisors on} $X$ is the direct limit 
\[
\CbDiv(X) \coloneqq \varinjlim_{\pi \in R(X)} \Div\left(X_{\pi}\right),
\]
taken in the category of topological vector spaces, with maps given by the pullback maps. The resulting topology
is called the \emph{strong} topology. The group of \emph{Weil $\rr$-b-divisors on} $X$ is the inverse limit 
\[
\WbDiv(X) \coloneqq \varprojlim_{\pi \in R(X)} \Div\left(X_{\pi}\right),
\]
taken in the category of topological vector spaces, with maps given by the pushforward maps. The resulting topology
is called the \emph{weak} topology.
\end{definition}

\begin{remark} \label{rem:2}
  As a set, $\CbDiv(X)$ can be seen as the disjoint union of the sets
  $\Div\left(X_{\pi}\right)$ modulo the equivalence  relation which sets two
  divisors equal if they coincide after pullback to a common modification.
  The set $\WbDiv(X)$ can be seen as the subset of $\prod_{\pi\in R(X)
  }\Div(X_{\pi})$ given by the elements $\bb D=\left(D_{\pi }\right)_{\pi \in
    R(X)}$ satisfying 
  the compatibility condition that for each $\pi '\ge \pi $ we have $\mu _{\ast}D_{\pi '}=D_{\pi }$, where $\mu
  $ is the corresponding modification.
\end{remark}
\begin{definition}
Let $\bb D$ be a Cartier $\rr$-b-divisor. A \emph{determination} of
$\bb D$ is a representative $D$ of the equivalence class given by $\bb D$ as described in
\ref{rem:2}. If $X_\pi$ is the modification where $D$ lives, we say
that $\bb D$ is \emph{determined} in $X_\pi$.  
If $\bb D=(D_{\pi })_{\pi \in R(X)}$ is a Weil $\rr$-b-divisor, then  
for $\pi \in R(X)$, the divisor $D_{\pi}$ is called the
\emph{incarnation} of $\bb D$ on $X_{\pi}$.
\end{definition}
\begin{remark}\label{rem:3}\label{b-prop}
\begin{enumerate}
\item
  For every modification $\mu \colon X_{\pi '}\to X_{\pi }$,
  the identity $\mu_*\mu^* = \operatorname{Id}$ holds. Therefore, 
  the natural map $\CbDiv(X) \to \WbDiv(X)$ is injective. This
  map can be described as follows. Let $D\in \Div(X_{\pi_{0} })$ be any
  determination of a Cartier b-divisor. Then for each $\pi \in R(X)$
  we choose any element $\pi '\in R(X)$ such that $\pi '\ge \pi _{0}$
  and $\pi '\ge \pi $. Let $\mu _{0}\colon X_{\pi '}\to X_{\pi_{0}} $ and
  $\mu \colon X_{\pi '}\to X_\pi $ be the corresponding modifications. Then
  $D_{\pi }\coloneqq \mu _{\ast}\mu_{0}^{\ast}D$ does not depend on
  the choice of $\pi '$ and 
  the image of the Cartier $\rr$-b-divisor given by $D$ is the Weil $\rr$-b-divisor  $(D_{\pi })_{\pi
  }$. From now on we will identify, as a set, $\CbDiv(X) $ with its
  image in $\WbDiv(X)$, and by \emph{$\bb R$-b-divisor} we will mean a Weil
  $\rr$-b-divisor. 
\item The injection $\CbDiv(X)\hra \WbDiv(X)$ as described above is continuous, but is \emph{not} a homeomorphism onto its image. In
  fact, a net of Cartier $\rr$-b-divisors $\left\{\bb D_i\right\}_{i
    \in I}$ converges in $\CbDiv(X)$ to a Cartier $\rr$-b-divisor $\bb D$ if and only if the following is satisfied. There
  exists a model $\pi $ such that  $\bb D$ and all the $\bb D_{i}$
  are determined in $\pi $, and if $D,D_{i}\in \Div\left(X_{\pi}\right)$ are
  determinations of $\bb D$ and $\bb D_{i}$, respectively, then   
  \[
    D = \lim_{i \in I}D_{i}
  \]
  in $\Div(X_{\pi})$. 
  On the other hand, a net of  
  $\rr$-b-divisors $\left\{\bb D_i\right\}_{i 
    \in I}$ converges in $\WbDiv(X)$ to an $\rr$-b-divisor $\bb D$
  if and only if for each model $\pi \in R(X)$, we have that 
  \[
    D_{\pi} = \lim_{i \in I}D_{i, \pi}
  \]
  in $\Div\left(X_{\pi}\right)$. 
\item Any $\rr$-b-divisor $\bb D = \left(D_{\pi}\right)_{\pi \in
    R(X)}$ is the limit of its incarnations
  $D_{\pi}$. It follows that $\CbDiv(X)$ is dense in $\WbDiv(X)$.  
\item It is natural in many situations to consider also integral or
  rational coefficients. The definitions are then easily adapted.
\end{enumerate}
\end{remark}

\subsection{Nef and approximable nef b-divisors}
\begin{definition}[{\cite{bff}}]
A Cartier $\rr$-b-divisor $\bb D \in \CbDiv(X)$ is \emph{nef} if $D_{\pi} \in \Div\left(X_{\pi}\right)$ is nef for one (and hence for every) determination $\pi \in R(X)$ of $\bb D$. 
A Weil $\rr$-b-divisor $\bb D \in \WbDiv(X)$ is \emph{nef} if it is a limit (in the weak topology) of a net of nef Cartier $\rr$-b-divisors. 
\end{definition}
\begin{remark} 
It is a priori not clear that if a nef Weil $\rr$-b-divisor is Cartier, then it is nef as a Cartier $\rr$-b-divisor. This is known to be true in the toroidal setting \cite[Lemma~4.24]{BoteroBurgos} and if we work with algebraic varieties over a countable field (instead of complex manifolds), see \cite[Corollary 4]{Da-Fa20}.

\end{remark}
\newcommand{\approach}{\on{Nef-b-Div}_{\rr}}
\begin{definition}[{\cite[Section 2]{Da-Fa20}}]\label{nef-app}
Let $\bb D$ be a nef $\rr$-b-divisor on $X$. Then $\bb D$ is called \emph{approximable nef} if $\bb D$ can be written as a limit
\[
\bb D = \lim_{i \in \bb N}\bb E_i
\]
of a sequence of nef Cartier $\rr$-b-divisors satisfying the monotonicity property $\bb E_i \geq \bb E_j$ whenever $i \leq j$. We call such a sequence an \emph{approximating sequence}. 

\end{definition}
\begin{remark}
In \cite[Theorem 5]{Da-Fa20} Dang and Favre show that in the case of algebraic varieties over a countable field, \emph{any} nef Weil $\bb R$-b-divisor is approximable nef. We will show in \ref{sec:psh-app} that any b-divisor that comes from a psh metric in a suitable sense is approximable nef. 
\end{remark}

\subsection{Intersection products of approximable nef b-divisors}\label{sec:top_intersection_products}

Let $n = \dim X$, let $\bb E_1, \dotsc, \bb E_{n-1}$ be Cartier $\rr$-b-divisors on $X$, and let $\bb D$ be a Weil $\rr$-b-divisor on $X$. 
If $\pi \in R(X)$ is such that all $\bb E_i$ are determined in $\pi$, then the real-valued intersection number 
\[
E_{1,\pi} \dotsm E_{n-1, \pi} \cdot D_{\pi}
\]
is independent of the choice of $\pi$, by the projection formula. This yields a well-defined intersection product
\begin{equation}\label{eqn:int-cl}
\CbDiv(X) \times \dotsc \times \CbDiv(X) \times \WbDiv(X) \longrightarrow \rr.
\end{equation}
Extending this to an intersection product on all Weil $\bb R$-b-divisors seems too much to ask for, in general. However, if $\bb D_1 , \dots, \bb D_n$ are \emph{approximable nef} $\rr$-b-divisors and $(D_{i, r})_r$ are approximating sequences for the $\bb D_i$, then one can show that the limit 
\begin{equation}
\lim_{r \to \infty} D_{1,r} \cdots D_{n, r}
\end{equation}
exists in $\rr_{\ge 0}$ and is independent of the choice of approximating sequences. This yields a top intersection product on the space of approximable nef b-divisors, which is continuous with respect to approximating sequences. The details can be found in \cite[Section~3]{Da-Fa20}. We mention that in \cite{Da-Fa20} the authors work over a countable ground field, but one can either check that this part of their argument does not use the countability assumption, or apply the following lemma: 
\begin{lemma}
Let $\bb D\in \WbDiv(X)$ be a limit of a sequence of Cartier divisors. Then there exists a countable subfield $L \sub \bb C$ over which both $X$ and $\bb D$ are defined. 
\end{lemma}

\section{The b-divisor associated to a psh metric}\label{sec:b-div-psh}

\subsection{The definition of the b-divisor}
\label{sec:definition-b-divisor}

Let $X$ be a projective complex manifold of dimension~$n$. Let $\ca L$ be a line bundle on $X$, and let $h$ be a psh metric on $\ca L$ (see \ref{def:singular_metric}). If $\pi\colon X_\pi \to X$ is a model in the category $R(X)$ (see \ref{def:modi}), we define the anti-effective $\bb R$-divisor
\begin{equation}
Z(\L,h)_{\pi} \coloneqq \sum_P - \nu(h,P)P 
\end{equation}
on $X_\pi$,
where the sum is over all prime divisors $P$ on $X_\pi$.  In general this sum need not be finite, so we make the following definition.
\begin{definition}\label{def:zar-unb}
 The psh metric $h$ on $\L$ has \emph{Zariski unbounded locus} if there exists  a non-empty Zariski open subset $U \subseteq X$ such that the local potentials of $h$ are locally bounded on $U$.
\end{definition}
 
\begin{example} By
\ref{has_unbounded}, any psh metric with almost asymptotically algebraic singularities in the sense of \ref{def:6} has
  Zariski unbounded locus. 
  \end{example}
\begin{lemma} \label{check_unbounded}
Assume that the psh metric $h$ on $\ca L$ has Zariski unbounded locus. Then, given any model $\pi \colon X_{\pi} \to X$ in $R(X)$, for only finitely many prime divisors $P$ on $X_\pi$ we have that $\nu(h,P)$ is non-zero. Moreover, if $\mu\colon X_{\pi'} \to X_{\pi}$ is a map of models then we have an equality of Weil $\bb R$-divisors
\begin{equation}
\mu_*Z(\L,h)_{\pi'} = Z(\L,h)_{\pi}. 
\end{equation}
on $X_{\pi'}$.
\end{lemma}
\begin{proof}
This is an easy consequence of the definitions.
\end{proof}
\begin{definition}\label{b-div-lelong}

Assume that the psh metric $h$ has Zariski unbounded locus. The Weil $\rr$-b-divisor $Z(\L,h) \in \WbDiv(X)$ is defined by 
\[
Z(\L,h) \coloneqq \left(Z(\L,h)_{\pi}\right)_{\pi \in R(X)}. 
\]
It follows from \ref{check_unbounded} that this is indeed a Weil $\rr$-b-divisor.

Let $s$ be a non-zero rational section of $\L$ and write $D = \on{div}(s)$, seen as a Cartier $\rr$-b-divisor on $X$. Then we define the Weil $\rr$-b-divisor 
\[
D(\L,h,s) \coloneqq D + Z(\L,h). 
\]
\end{definition}
We observe that the formation of $D(\L,h,s)$ is multiplicative in the sense that
$D(\L_1 \otimes \L_2, h_1 \otimes h_2, s_1 \otimes s_2) = D(\L_1,h_1,s_1) +  D(\L_2,h_2,s_2)$ whenever $(\L_1,h_1)$, $(\L_2,h_2)$ are line bundles with psh metrics with Zariski unbounded locus on $X$ and $s_1$, $s_2$ are non-zero rational sections of $\L_1$ resp.\ $\L_2$.

We write 
\[
D(\L,h,s) = \left(D(\L,h,s)_\pi\right)_{\pi \in R(X)},
\]
where 
\begin{equation*}
D(\L,h,s)_\pi = \pi^*D + Z(\L,h)_\pi = \on{div}_{X_\pi} (s) + Z(\L,h)_\pi \, . 
\end{equation*}

\begin{example}\label{exm:4}
  Let notations be as above and assume that $h$ has algebraic singularities. Then $
  D(\L,h,s)$ belongs to $\CbDiv(X)$. Indeed, it can be computed as
  follows. Fix a reference smooth metric $h_{0}$ on $\ca L$ and write $\varphi = -
  \log(h(s)/h_{0}(s))$ and $\theta = c_1(\ca L, h_0)$. As in \ref{rem:bij-metrics} we have that $\varphi$
  is a $\theta $-psh function on $X$ and  it has algebraic
  singularities by assumption.   
  Let $c$ be the rational constant from \ref{def:3} for $\varphi$ and write
  $\mathcal{I}=\mathcal{I}(\varphi/c)$. It follows \ref{lemm:1} from there is a model $\pi \colon X_{\pi}\to X$ in $R(X)$ such that $\pi^{-1}\mathcal{I}\cdot\ca O_{X_\pi} =\mathcal{O}(-D)$ for an effective simple
  normal crossings divisor $D$, and such that
  \begin{displaymath}
    D(\L,h,s) = \divisor( s) - cD
  \end{displaymath}
  as Cartier b-divisors. Note that $D(\L,h,s)$ is actually a $\qq$-b-divisor in this case.
\end{example}
Note that also a psh metric that is good in the sense of Mumford gives rise to a Cartier b-divisor (all Lelong numbers are zero on all modifications). This shows that the converse to the statement in the example does not hold as good metrics do not necessarily have algebraic singularities. 
\begin{remark}\label{rem:b-div-metric}
Let $T$ be the closed positive $(1,1)$-current on $X$ given as $T = c_1(\mathcal{L},h)$ and let $s$ be a non-zero rational section of $\mathcal{L}$. We can relate the b-divisor $D(\L,h,s)$ to the Siu decomposition of $T$ on $X$ (see \ref{rem:siu-dec}). Recall that this decomposes $T$ uniquely as a sum 
\[
T = R + \sum_k\nu\left(T, Y_k \right) \delta_{Y_k},
\]
where the sum is over an at most countable family of $1$-codimensional subvarieties $Y_k$ of $X$ and $R$ is a closed positive $(1,1)$-current whose Lelong number on any divisor is zero. The sum $\sum_k\nu\left(T, Y_k \right) \delta_{Y_k}$ is called the \emph{divisor part of $T$}. If $h$ has Zariski unbounded locus, then the family $\left\{Y_k\right\}$ is finite and the divisor part agrees with $-Z(\mathcal L,h)_{\on{id}}$. 

Also, for each $\pi \in R(X)$ we may consider the closed positive current $T_{\pi} = c_1\left(\pi^*\mathcal{L}, \pi^*h\right)$ on $X_{\pi}$, where $\pi^*h$ denotes the pullback metric whose local potentials are defined by pulling back the local psh potentials of~$h$. Then the divisor part of $T_{\pi}$ agrees with $-Z(\mathcal L,h)_{\pi}$. 

We see that the b-divisor $D(\L,h,s)$ encodes the divisor parts of the Siu decomposition of the pullbacks of $T$ on all models $\pi \in R(X)$, and hence it also encodes in some sense the parts of higher codimension of the singular locus of the metric. 
\end{remark}

\begin{example} Let $h$ and $h'$ be two psh metrics on $\L$ with Zariski unbounded locus. Let $s$ be a non-zero rational section of $\L$ and denote by $D(\L,h,s)$ and $D(\L,h',s)$ the associated $\bb R$-b-divisors on $R(X)$. 
Fix a smooth reference metric $h_0$ on $\L$ and write $\varphi = -\log\left(h(s)/h_0(s)\right)$ and $\varphi' = -\log\left(h'(s)/h_0(s)\right)$. Let $\theta = c_1(\L,h_0)$. Then $\varphi$ and $\varphi'$ are $\theta$-psh functions on $X$. 
Recall the notion of algebraic singularity type and the equivalence relation $\prec_{\mathcal J}$ from \ref{def:alg-sing-type}. We have
\begin{align}
D(\L,h,s) \geq D(\L,h',s) \quad \text{iff}\quad \varphi' \prec_{\mathcal J}\varphi.
\end{align}
This is just a rephrasing of \ref{prop:alg-sing-type}. 
\end{example}

If $h$ is a psh metric on $\ca L$ with Zariski unbounded locus and
$U\subset X$ is a dense Zariski open subset such that the local potentials
of $h$ are locally bounded on $U$, then, for any rational section $s$
of $\ca L$, the b-divisor $D(\ca L,s,h)$ only depends on the
restriction of $(\ca L,s,h)$ to $U$. More concretely,

\begin{proposition}\label{prop:3} Let $U\subset X$ be a dense Zariski
  open subset and $\ca L$ a line bundle on $U$ with a psh metric
  $h$. Let $X_{1}$ and $X_{2}$ be two compactifications of $U$ and
  $\ca L_{1}$ and $\ca L_{2}$ line bundles on $X_{1}$ and $X_{2}$ together with isomorphisms $\ca L_{1}|_{U} \isom \ca L^{\otimes e_{1}}$ and $\ca L_{2}|_{U}\isom \ca
  L^{\otimes e_{2}}$ for some integers $e_{1}, e_{2}>0$. Assume
  moreover that the
  metric $h$ extends to singular psh metrics $h_{1}$ and $h_{2}$ on
  $\ca L_{1}$ and $\ca L_{2}$, respectively. Let $s$ be a non-zero rational section of $\ca
  L$, so $s^{\otimes e_{i}}$ is a non-zero rational section of $\ca L_{i}$. Then
  \begin{displaymath}
    \frac{1}{e_{1}}D(\ca L_{1},s^{\otimes e_{1}},h_{1})=
    \frac{1}{e_{2}}D(\ca L_{2},s^{\otimes e_{2}},h_{2}).
  \end{displaymath}
\end{proposition}
\begin{proof}
By considering a high enough model dominating both $X_1$ and $X_2$, and high enough powers of the $\ca L_i$, we may reduce to the case $X = X_1=X_2$, $e_1 = e_2=1$. Then by symmetry we may further reduce to the case where $\ca L_1 = \ca L_2 \otimes \ca O(D)$ with $D$ an effective divisor such that $\on{supp}(D) \subset X \setminus U$. Let $P$ be a prime divisor on any modification $X'$ of $X$; we need to show that 
\begin{align}\label{eq:lelong}
\on{ord}_P(s, \ca L_1) - \nu(P,h_1, \ca L_1) = \on{ord}_P(s, \ca L_2) - \nu(P,h_2, \ca L_2). 
\end{align}
Let $r = \on{coeff}_P(D)$. Then 
\begin{align}\label{eq:lelong2}
\on{ord}_P(s, \ca L_1)- \on{ord}_P(s, \ca L_2) = r.
\end{align}
Let $g$ be a local equation of $P$ and for $i = 1,2$, let $s_i$ be local invertible sections of the $\ca L_i$'s at $P$. We can write $s_2 = s_1 \cdot g^r$. Hence 
\[
\log \|s_2 \| = \log \|s_1\| + r \log \| s_1 \|
\]
which implies that 
\begin{align}\label{eq:lelong3}
r = \nu(P,h_1, \ca L_1) - \nu(P,h_2, \ca L_2).
\end{align}
Putting \ref{eq:lelong2} and \ref{eq:lelong3} together yields \ref{eq:lelong}. 
\end{proof}

\subsection{b-divisors coming from psh metrics are approximable nef}\label{sec:psh-app}
In this section we show that the $\bb R$-b-divisors associated to psh metrics with Zariski unbounded locus are approximable nef (see \ref{nef-app}). 

As before let $X$ be a  projective complex manifold of dimension $n$, let $\ca L$ on $X$ be a line bundle, and let $h$ be a psh metric on $\ca L$. We choose a canonical divisor $K_X$ on $X$, a very ample line bundle $B$ on $X$, and a smooth hermitian metric $g$ on $B$. Let $\omega \coloneqq c_1(B,g)$ and assume it is a K\"ahler form on $X$. We essentially globalize some of the arguments found in \cite[Sections~5.1 and 5.2]{bfj-val}. 

We start with some preparatory results (the first of which only needs $X$ to be a compact K\"ahler manifold). 
\begin{theorem} [{Analytic Nadel Vanishing, cf. \cite[Theorem 9.4.21]{Laz}}] \label{nadel}
Assume that $c_1(\ca L,h) \ge \epsilon \omega$ for some $\epsilon > 0$. Then
\begin{equation*}
H^i(X, \ca O_X(K_X ) \otimes  \ca L \otimes \ca J(h)) =0 \;\; \text{for} \;\; i > 0. 
\end{equation*}
\end{theorem}
\begin{theorem} [{\cite[Lecture 14]{mum-lec}}] \label{th:glo}
Let $\ca F$ be a coherent sheaf on $X$ such that 
\begin{equation*}
H^i(X, \ca F \otimes B^{\otimes(k-i)}) = 0 \; \text{for all} \; i > 0 \; \text{and} \; k \ge 0. 
\end{equation*}
Then $\ca F$ is globally generated. 
\end{theorem}
\begin{corollary}\label{cor:global_gen}
Assume that $c_1(\ca L,h) \ge \epsilon  \omega$ for some $\epsilon > 0$. Then 
\begin{equation*}
\ca O(K_X) \otimes \ca L \otimes \ca J(h) \otimes B^{\otimes n}
\end{equation*}	
is globally generated. 
\end{corollary}
\begin{proof}
Let $\ca F = \ca O(K_X) \otimes \ca L \otimes \ca J(h) \otimes B^{\otimes n}$. Then $\ca F$ is coherent. By \ref{th:glo} it suffices to show that 
\begin{equation*}
H^i(X, \ca O(K_X) \otimes\ca L \otimes \ca J(h) \otimes B^{\otimes n} \otimes B^{\otimes(k-i)}) = 0 \; \textrm{for} \; i > 0 \; \textrm{and} \; k \ge 0. 
\end{equation*}
For $i > n$ this follows from the fact that $\dim X=n$,  and if $i \le n $ the vanishing follows from \ref{nadel} applied to $\ca L \otimes B^{\otimes (n+k-i)}$. 
\end{proof}

\begin{definition} \label{def:cartier_divisor_from_ideal}
Let $\ca J$ be a non-zero coherent sheaf of ideals on $X$, with its canonical rational section $1$. If $X_\pi$ is a model on which the inverse image ideal sheaf $\pi^{-1}\ca J\cdot\ca O_{X_\pi}$ of $\ca J$ is invertible, then the pullback of $1$ in $\pi^{-1}\ca J\cdot\ca O_{X_\pi}$ defines a divisor on $X_\pi$. This determines an anti-effective Cartier b-divisor on $X$, independent of the choice of $X_\pi$, which we denote by $Z(\ca J)$.

\end{definition}

\begin{lemma}\label{lem:g_gen_nef}
In the notation of \ref{def:cartier_divisor_from_ideal}, suppose that $\ca L \otimes \ca J$ is globally generated. Then the line bundle $\pi^*\ca L \otimes \left(\pi^{-1}\ca J\cdot\ca O_{X_\pi}\right)$ on $X_\pi$ is globally generated, in particular, it is nef. 

\end{lemma}
\begin{proof}
The sheaf $\pi^*(\ca L\otimes \ca J)$ is globally generated, and the canonical epimorphism $\pi^*(\ca L\otimes \ca J ) \to \pi^*\ca L \otimes \left(\pi^{-1}\ca J\cdot\ca O_{X_\pi}\right)$ shows that the latter is globally generated. It is therefore the pullback of $\ca O(1)$ under a morphism to projective space, and so it is nef. 

\end{proof}
\begin{definition}
If $\pi\colon X_\pi \to X$ is a model in $R(X)$, then the relative dualising sheaf is \emph{canonically} trivial over the locus where the map $\pi$ is an isomorphism, and so comes with a canonical rational section, which is in fact a regular section, and defines an effective relative canonical \emph{divisor} $K_{\pi}$. 
If $\mu\colon X_\pi \to X_{\pi'}$ is a map of models then $K_{\pi} =
\mu^*K_{\pi'} + K_{\mu}$ (see also \cite[Section 3.1]{bfj-val}), so
the $K_\pi$ assemble to an effective $\rr$-b-divisor on $X$.

\end{definition}
We mention that $K_\pi$ can equivalently be computed by taking the jacobian ideal sheaf $\on{Jac}(\pi)$. 
\begin{lemma} \label{lem:sandwich}

Assume that the psh metric $h$ has Zariski unbounded locus, and let
$D$ be a reduced effective divisor on $X$ outside of which $(\ca L,h)$
has bounded local potentials. For $\pi \in R(X)$, we let $\pi^{\ast}
D_{\text {\rm red}}$ denote the reduced divisor with same support as $\pi^{\ast}
D$. Let $\bb A = (A_\pi)_\pi $ be the $\bb R$-b-divisor given on
each model $\pi \in R(X)$ by 
\begin{equation*}
A_\pi = \pi^{\ast} D_{\text {\rm red}} + K_\pi. 
\end{equation*}
Then we have the following inequalities of $\rr$-b-divisors 
\begin{equation}\label{eq:sandwich}
Z(\ca J(h)) - \bb A \le Z(\ca L, h) \le Z(\ca J(h)).
\end{equation}
\end{lemma}

We refer to \ref{def:multiplier_sheaf} for the definition of the
multiplier ideal sheaf $\ca J(h)$ associated to the metric~$h$.

\begin{proof}

We follow the idea of the proof of \cite[Theorem~5.1]{bfj-val}.
  Let $\pi \in R(X)$ and let $U$ be a small ball in $X_\pi$, which we identify with a small ball centered at the origin in $\bb C^n$. We assume that $\pi^*\ca L$ admits a generating section $\xi$ over $U$, and we let $\varphi = -\log h(\xi)$, which is then a psh function on $U$. We need to check the inequalities of divisors
\[ Z( \ca J(\varphi)) - A_\pi|_U \le Z(\varphi)  \le Z(\ca J(\varphi))  \]
on $U \subset X_\pi$. 

The inequality $Z(\varphi) \le Z(\ca J(\varphi))$ is standard and follows from the Ohsawa-Takegoshi extension theorem (e.g.~see the proof of the first inequality in part a) of \cite[Theorem 13.2]{dem-anal}). 

It thus remains to show the inequality $Z(\ca J(\varphi)) -  A_\pi|_U \le Z(\varphi) $. 
We first check this result `before blowing up'; in other words, we let $P$ be a prime divisor of $X$. If $P$ is not contained in the support of $D$ then both sides of the inequality vanish at $P$ and we are done. So assume that $P$ is contained in the support of $D$, so $\on{ord}_P \bb A = 1$. Then $\ca J(\varphi)$ is principal at the generic point of $P$, say generated by $f$. Then 
\begin{equation*}
\on{ord}_PZ(\ca J(\varphi)) = - \on{ord}_P f = -\nu(\log\abs{f}, p) \, , 
\end{equation*}
where $p$ is a generic point of $P$. 
To prove the claim it suffices then to show that 
\[
\nu(\log\abs{f}, p) \ge \nu(\varphi, p) -1.
\] 
This comes down to showing that $\varphi(z) \le c \log \abs{z} + O(1)$ implies that we have $c \le 1 + \on{ord}_p(f)$. 
Hence assume that $\varphi(z) \le c \log \abs{z} + O(1)$.
Then $e^{2\varphi} \ll \abs{z}^{2c}$
 and thus $\abs{f}^2e^{-2\varphi} \gg \abs{z}^{2\on{ord}_pf} \abs{z}^{-2c}.$
  Since $\abs{f}^2e^{-2\varphi}$ is assumed to be integrable it follows that $\abs{z}^{2\on{ord}_p(f) -2c}$ is integrable. Hence $2\on{ord}_p(f) -2c \ge -2$ as required. 
It remains to treat the case where $P$ is an exceptional prime divisor
on some $X_\pi$. We denote by $\on{Jac}(\pi)$ the jacobian ideal sheaf of $\pi$, and recall
that the relative canonical divisor $K_\pi$ is the divisor cut out by
$\on{Jac}(\pi)$. On a small disc in $X_\pi$ centered at a generic point $p$
of $P$, write $j$ for a generator of
$\on{Jac}(\pi)$. Then if $f$ on $X$ is a function such that
$\abs{f}^2e^{-2\phi}$ is locally integrable on $X$ we have that
$\abs{f}^2e^{-2\phi}\abs{j}^2$ is locally integrable on $X_\pi$. Hence
we see that  
\begin{equation*}
\on{ord}_PZ(\ca J(\phi \circ \pi)) \le \on{ord}_PZ(\ca J(\phi)) + \on{ord}_P K_\pi. 
\end{equation*}
Using this inequality, the
claim follows as in the case where $P$ was a prime divisor on~$X$.
\end{proof}
Note that the $\rr$-b-divisor $\bb A$ only depends on the divisor $D$ on which $h$ has singularities; in particular, if we replace $(\ca L,h)$ by some positive power, then $\bb A$ need not change. Hence, from the previous lemma we obtain the following corollary. 
\begin{corollary}\label{cor:conv-lelong}
\[
Z(\L,h) = \lim_{k \to \infty}\frac{1}{k}Z\left(\mathcal{J}(h^{\otimes k})\right).
\]
\end{corollary}
\begin{remark}
Stated in a different language \ref{cor:conv-lelong} can also be found in \cite[Proposition 2.14]{DX}, with essentially the same proof. 
\end{remark}
We are now ready to state and prove the main result of this section.
\begin{theorem}\label{thm:psh_implies_approximable}
Let $X$ be a projective complex manifold of dimension~$n$, let $\ca L$ on $X$ be a line bundle, and let $h$ be a psh metric on $\ca L$ with Zariski unbounded locus. Let $s$ be a non-zero rational section of $\ca L$. Then the associated $\rr$-b-divisor $D(\ca L,h,s) = \on{b-div}(s) + Z(\ca L,h)$ is approximable nef. 
\end{theorem}
In particular, by what was said in \ref{sec:top_intersection_products} the $\rr$-b-divisor $D(\ca L,h,s)$ has a well-defined degree $D(\ca L,h,s)^n \in \rr_{\ge 0}$.

\begin{proof}

We choose a canonical divisor $K_X$ on $X$, a very ample line bundle $B$ on $X$, and a smooth hermitian metric $g$ on $B$. Let $\omega \coloneqq c_1(B,g)$ which we assume to be a K\"ahler form on $X$.
Let $k \in \zz_{> 0}$.
By \ref{cor:global_gen} applied to $\ca L^{\otimes k} \otimes B$, using that  $c_1(\ca L^{\otimes k} \otimes B) \ge \omega$ and that $\ca J(h^k) = \ca J(h^k g)$, we see that the sheaf on $X$ given by
\begin{equation}\label{eq:globally_gen}
\ca O(K_X) \otimes \ca L^{\otimes k} \otimes \ca J(h^k) \otimes B^{\otimes n+1}
\end{equation}	
is globally generated. Choose a non-zero rational section $b$ of $B$, determining a Cartier $\rr$-b-divisor $\on{div}(b)$ on $X$. We also write $K_X$ for the Cartier $\rr$-b-divisor determined by $K_X$. Recall that every $Z(\ca J(h^{\otimes k}))$ is Cartier. 
By \ref{lem:g_gen_nef} and the fact that \ref{eq:globally_gen} is globally generated, the Cartier $\rr$-b-divisor given by
\begin{equation*}
K_X + k\on{div}(s) + Z(\ca J(h^k)) + (n+1)\on{div}(b)
\end{equation*}
 is nef. Dividing by $k$ we find
\begin{equation}
\frac{1}{k}K_X + \on{div}(s)  + \frac{1}{k}Z(\ca J(h^k)) + \frac{n+1}{k}\on{div}(b)
\end{equation}	
is again a nef Cartier $\rr$-b-divisor. Let $C := K_X + (n+1)\on{div}(b)$, a divisor on $X$. Write $C = E - P$ where $E$ is effective and $P$ is nef. Then for each $k$ we have 
\begin{equation*}
\frac{1}{k}K_X + \on{div}(s) + \frac{1}{k}Z(\ca J(h^k)) + \frac{n+1}{k}\on{div}(b)  + \frac{P}{k} =  \on{div}(s) + \frac{1}{k}Z(\ca J(h^k)) + \frac{E}{k}. 
\end{equation*}
Hence we see that 
 \[ \on{div}(s) + \frac{1}{k}Z(\ca J(h^k)) + \frac{E}{k} \]
is a nef $\rr$-b-Cartier divisor, as it is a sum of such. Since $E$ is effective, the sequence $E/k$ is decreasing as $k \to \infty$. Next, we have that a subsequence of $Z(\ca J(h^k))/k$ is decreasing. Indeed, this follows from the fact that for all $k$, we have
  \[ \frac{1}{k}Z(\ca J(h^k))\ge \frac{1}{2k}Z(\ca J(h^{2k})), \]
which follows from the additivity property stated in \ref{lem:subadditivity} below. The proof is then concluded by \ref{cor:conv-lelong}.
\end{proof}

\begin{lemma}[{\cite[Theorem 14.2]{dem-anal}}] \label{lem:subadditivity}
Let $\varphi_1$ and $\varphi_2$ be psh functions. Then $\ca J(\varphi_1 + \varphi_2) \subseteq \ca J(\varphi_1)\ca J(\varphi_2)$. 
\end{lemma}

\subsection{A Chern--Weil type result for psh metrics with almost asymptotically algebraic singularities}
\label{sec:proof-refth-weil}

Let $X$ be a projective complex manifold of dimension $n$ and let $\L$ be a line bundle on $X$. Let $h$ be a psh metric on $\L$, let $s$ be a non-zero rational section of $\ca L$  and let $D(\L,h,s)$ be the associated $\bb R$-b-divisor on $X$ as in \ref{b-div-lelong}.  We fix a smooth reference metric $h_0$ on $\L$ and write $\theta = c_1(\ca L,h_0)$.

\begin{theorem}\label{thm:b-div-chern-weil}
Assume that the psh metric
$h$ has almost asymptotically algebraic singularities with respect to $\theta$.  Then the equality
\[
D(\ca L,h,s)^n = \int_X \langle c_1(\mathcal{L},h)^n \rangle 
\]
holds in $\rr_{\ge 0}$. Here $\langle \cdot \rangle$ denotes the non-pluripolar product of closed positive $(1,1)$-currents, and $D(\ca L,h,s)^n$  is the degree of the approximable nef $\bb R$-b-divisor $D(\ca L,h,s)$.
\end{theorem}
\begin{proof}
Write $\bb D= D(\ca L,h,s)$ and $T=c_1(\ca L,h)$. If $h$ has algebraic singularities then the result follows from  combining \ref{prop:alg-sing} and \ref{exm:4}.

Assume now that $h$ has almost asymptotically algebraic singularities. Let $B$ be a very ample line bundle on $X$, $g$ a K\"ahler metric on $B$, and set $\omega = c_1(B, g)$ and $\varphi = -\log\left(h(s)/h_0(s)\right)$.  By \ref{lem:independent_Kahler} we may choose an $\omega$-psh function $f$ with algebraic singularities on $X$, and for every integer $m>0$ an $(\frac{1}{m}\omega + \theta)$-psh function $\varphi_{m}$ with algebraic singularities on $X$ such that
\begin{displaymath}
\varphi_{m}+ \frac{1}{m} f \prec \varphi \prec\varphi_{m} \, .
\end{displaymath}
Each $\varphi_{m}$ defines a psh metric $h_m$ with algebraic singularities on $\mathcal{L}\otimes B^{\otimes \frac{1}{m}}$, and each $\varphi_m + \frac{1}{m}f$ defines a psh metric $h_m'$ with algebraic singularities on $\mathcal{L}\otimes B^{\otimes \frac{2}{m}}$. 
Choose a non-zero rational section $t$ of $B$. 
Since the Lelong numbers of $\varphi_{m}$ converge to the Lelong
numbers of $\varphi$, the $ D(\mathcal{L}\otimes B^{\otimes \frac{1}{m}},h_{m},s t^{\frac{1}{m}})$
converge to $\bb D$, hence the $ D(\mathcal{L}\otimes B^{\otimes \frac{1}{m}},h_{m},s t^{\frac{1}{m}})^n$ converge
to $\bb D^{n}$. Similarly the $ D(\mathcal{L}\otimes B^{\otimes \frac{2}{m}}, h'_{m}, s t^{\frac{2}{m}})^n$ converge
to $\bb D^{n}$. 
The case of algebraic singularities yields that 
\begin{equation*}
\int_{X}\langle (\theta + \frac{1}{m} \omega + dd^c \varphi_m )^{n} \rangle =  D(\mathcal{L}\otimes B^{\otimes \frac{1}{m}},h_{m},s t^{\frac{1}{m}})^n
\end{equation*}
and
\begin{equation*}
  \int_{X}\langle (\theta + \frac{2}{m} \omega + dd^c (\varphi_m  + \frac{1}{m} f))^{n} \rangle = D(\mathcal{L}\otimes B^{\otimes \frac{2}{m}},h'_{m}, s t^{\frac{2}{m}})^n
\end{equation*}
for each $m$. 
Since the $\theta$-psh function $\varphi$ is also $\frac{1}{m}\omega + \theta $-psh, we can apply monotonicity of non-pluripolar products (\ref{thm:1}) to deduce 
\begin{equation*}
 \int_{X}\langle (T + \frac{1}{m} \omega)^{n} \rangle =  \int_{X}\langle (\theta + \frac{1}{m} \omega + dd^c \varphi)^{n} \rangle \le \int_{X}\langle (\theta + \frac{1}{m} \omega + dd^c \varphi_m )^{n}\rangle 
\end{equation*}
for each $m$, and similarly 
\begin{equation*}
  \int_{X}\langle (\theta + \frac{2}{m} \omega + dd^c( \varphi_m  + \frac{1}{m} f))^{n} \rangle \le   \int_{X}\langle (\theta + \frac{2}{m} \omega + dd^c \varphi)^{n} \rangle =  \int_{X}\langle (T + \frac{2}{m} \omega) ^{n}\rangle. 
\end{equation*}
Now by the multi-linearity of the non-pluripolar product as expressed in \ref{prop:6} we see
\begin{equation}
 \int_{X}\langle (T + \frac{1}{m} \omega) \rangle^{n} = \sum_{i=0}^n \binom{n}{i}
 \frac{1}{m^i} \int_X \langle T^{n-i} \omega^i \rangle, 
 \end{equation}
 so that 
\begin{align}\label{eqn:ineq}
 \int_X \langle T^n \rangle   &= \lim_{m \to \infty}  \int_{X} \langle (T + \frac{1}{m} \omega) \rangle^{n} \le \lim_{m \to \infty} \int_{X}\langle (\theta + \frac{1}{m} \omega + dd^c \varphi_m )^{n}\rangle \nonumber \\
 &= \lim_{m \to \infty} D(\mathcal{L}\otimes B^{\otimes \frac{1}{m}},h_{m},s t^{\frac{1}{m}})^n = \bb D^n \, . 
\end{align}
Similarly
 \begin{equation*}
\begin{split}
 \bb D^n  = \lim_{m \to \infty} & \bb D(\mathcal{L}\otimes B^{\otimes \frac{2}{m}},h'_{m},s t^{\frac{2}{m}})^n = \lim_{m \to \infty}   \int_{X}\langle (\theta + \frac{2}{m} \omega + dd^c (\varphi_m  + \frac{1}{m} f))^{n} \rangle \\
 &\le \lim_{m \to \infty} \int_{X}\langle (T + \frac{2}{m} \omega) ^{n}\rangle  =  \int_X \langle T^n \rangle . \qedhere
\end{split}
\end{equation*}

  \end{proof}

\begin{example} (Chern--Weil for good and psh metrics) Suppose that
  $h$ is, at the same time, good and psh as in \ref{exm:6}. Then $h$
  has almost asymptotically algebraic singularities by
  \ref{ex:good_implies_aaas}. Since all the Lelong numbers at all
  points on all modifications of $X$ are zero, we have that $D(\ca L,
  h, s) = \operatorname{div}(s)$, seen as a Cartier
  b-divisor. Following the computations in \ref{exm:6}, we have that
\[
D(\ca L, h,s)^n = \on{div}_X(s)^n = \deg(\ca L) = \int_X\langle c_1(\ca L, h)^n \rangle \, ,
\]
and this verifies the Theorem in the case of good and psh metric.
\end{example}

\begin{remark} \label{rem:ineq}
Let $h$ be a psh metric on $\ca L$ with Zariski unbounded locus. We do not assume that $h$ has almost asymptotically algebraic singularities. Then we still have the inequality
  \[
  D(\L,h,s)^n \geq \int_X \langle c_1(\ca L, h)^n \rangle.
      \]
Indeed, consider a Demailly approximation sequence $\{\varphi_m\}_{m \in \bb N}$ for $\varphi$. Then $\varphi \prec \varphi_m$ and in the same way as in the proof of \ref{thm:b-div-chern-weil}, see in particular equation \ref{eqn:ineq}, the inequality follows. 
 \end{remark}
 Combining \ref{th:asym} and \ref{thm:b-div-chern-weil} we obtain the following b-divisorial analogue of the classical Hilbert--Samuel type statement for nef line bundles.
 \begin{corollary}\label{cor:hs}
 Let assumptions be as in \ref{thm:b-div-chern-weil}. Then the equality
\[ D(\ca L,h,s)^n = \vol_{\ca J}(\ca L,h) \]
holds in $\bb R_{\ge 0}$.
\end{corollary}

\section{The line bundle of Siegel--Jacobi forms} \label{sec:Siegel-Jacobi}

The purpose of this section is to exhibit an application of our results  in the context of the line  bundle of Siegel--Jacobi forms on the universal abelian variety over the fine moduli space $A_{g,N}$ of principally polarized complex abelian varieties of dimension~$g$ with level~$N$ structure. The results in this section form a generalization of the main results of \cite{bkk}.

\subsection{The biextension metric on the Poincar\'e bundle} \label{sec:PBTor}

We start by showing that the Poincar\'e bundle on an abelian scheme has a natural psh extension over any smooth toroidal compactification of the abelian scheme.

Let $S$ be a smooth complex algebraic variety, and let $\pi \colon U
\to S$ be an abelian scheme with zero section $e \colon S \to U$. That
is, the morphism $\pi$ is proper and smooth, and the fibers of $\pi$
are abelian varieties with origin determined by the section~$e$. We
have a tautological Poincar\'e line bundle $\ca P$ on the fiber
product $U  \times_S U^\lor$, where $\pi^\lor \colon U^\lor \to S$
denotes the dual abelian scheme. The Poincar\'e bundle $\ca P$ comes
equipped with a rigidification along the zero section and with a
canonical smooth hermitian metric $h_{0}$ as explained in
\cite[Section~3]{mb}. When $\lambda \colon U  \to U^\lor$ is a
polarization of abelian schemes over $S$, we define $\ca B_\lambda$ to
be the line bundle $(\id,\lambda)^*\ca P$ on $U$, equipped with the
pullback metric, which we denote by~$h_\lambda$. 

The line bundle $\ca B_\lambda$ is an example of a \emph{biextension line bundle} associated to a polarized variation of pure Hodge structures of weight~$-1$, as discussed in \cite[Sections 6 and 7]{hrar}. In our case, the underlying variation of pure Hodge structures of weight $-1$ is given by the local system $\ca H = R^1 \pi_*\zz_U(1)$, and the polarization of the variation $\ca H$ is induced by the polarization  $\lambda$. 

As a special case of general results such as \cite[Theorem~13.1]{hain_normal} or \cite[Theorem~8.2]{pp}, or alternatively by a computation in local coordinates using the explicit formulas for $h_{0}$ in  \cite[Section~2]{bghdj_sing}, we have that the metric $h_\lambda$ on the biextension line bundle $\ca B_\lambda$ is  \emph{semipositive}.

Let $X \supset U$ be a smooth projective compactification of $U$, and assume that $D= X \setminus U$ is a normal crossings divisor on $X$. We will assume throughout that the pullback of the variation of pure Hodge structure $\ca H$ from $S$ to $U$ has unipotent monodromy around each local branch of $D$. 

Let $D^{\mathrm{sing}}$ denote the singular locus of $D$. 
As a special case of \cite[Theorems~24 and 27]{bp} we have the
following result.

\begin{theorem} \label{bp_general}There exists a positive integer $e$ such that:
\begin{itemize}
\item[(i)]
The semipositive line bundle $\ca B_\lambda^{\otimes e}$ on $U$ has  a unique extension $\widetilde{\ca B^{\otimes e}_\lambda}$ as a continuously metrized line bundle over the locus $X \setminus D^{\mathrm{sing}}$. 
\item[(ii)] The continuously metrized line bundle $\widetilde{\ca B^{\otimes e}_\lambda}$ on $X \setminus D^{\mathrm{sing}}$ has a unique extension as a line bundle with a psh metric $\overline{\ca B^{\otimes e}_\lambda}$ over $X$.
\end{itemize}
\end{theorem}

\begin{remark}
To avoid carrying the exponent $e$ around in what follows, we can think of $\left(\overline{\ca B^{\otimes e}_\lambda}\right)^{\otimes \frac{1}{e}}$ as an extension of $\ca B_\lambda$ as a metrised $\bb Q$-line bundle over $X$, and we simply denote it $\overline{\ca B}_\lambda$ (which is independent of the choice of $e$). We say a metric on a $\bb Q$-line bundle $\ca L$ on a projective complex manifold is psh (or toroidal, or has almost asymptotically algebraic singularities, ...) if some positive tensor power of $\ca L$ which is a line bundle has this property. 
\end{remark}

 Following terminology introduced in \cite{hrar} we call the $\bb Q$-line bundle $\overline{\ca B}_\lambda$ the \emph{Lear extension} of $\ca B_\lambda$ over $X$. We will denote by $\overline{h}_\lambda$ the natural psh metric on $\overline{\ca B}_\lambda$ which is given by \ref{bp_general}. Note that the singularities of the psh metric $\overline{h}_\lambda$ are contained in the codimension two locus $D^{\mathrm{sing}}$ of $X$.
 
 The main result of this section is as follows.

\begin{theorem} \label{thm:biext_toroidal} 
The $\bb Q$-line bundle $\overline{\ca B}_\lambda$ with psh  metric $\overline{h}_\lambda$ on $X$ has toroidal singularities in the sense of \ref{def:toroidal_sing}. In particular, the psh metric $\overline{h}_\lambda$ has almost asymptotically algebraic singularities.
\end{theorem}
Note that the second part of the theorem follows from the first by \ref{prop:toroidal_implies_aaas}.

We start with two lemmas.

\begin{lemma}\label{lem:lipschitz_on_open}
Let $V \sub \bb R^k$ be any subset, and let $g \colon V \to \bb R$ be a Lipschitz continuous function with Lipschitz constant $c$. Then $g$ can be extended to a Lipschitz continuous function $\bb R^k \to \bb R$ with constant $c$. 
\end{lemma}
\begin{proof}
One extension is given by  $y \mapsto \inf_{x \in V}( g(x) + c\abs{y - x})$. 
\end{proof}

\begin{lemma} \label{lem:is_Lipschitz_cts} Let $r \in \zz_{> 0}$.
Let $A_1,\ldots,A_k$ be positive semi-definite $r \times r$-matrices
such that for all $x_1,\ldots,x_k \in \rr_{>0}$ we have $\sum_{i=1}^k
x_i A_i $ of full rank. Let $c_1,\ldots, c_k \in \rr^r$ be column
vectors. Then the smooth function
\begin{equation}
g = \left(\sum_{i=1}^k x_{i}A_i c_i \right)^t \cdot \left(\sum_{i=1}^k x_{i}A_i  \right)^{-1} \cdot \left(\sum_{i=1}^k x_{i}A_i c_i  \right)
\end{equation}
on $\rr_{>0}^k$ has a unique continuous extension to the whole of $\rr_{\ge0}^k$, and the extension is Lipschitz continuous. 
\end{lemma}
\begin{proof}
Define $P = \{(x_1, \dots, x_k) \in \rr_{>0}^k | x_1 < x_2 < \cdots < x_k\}$. Then the union of the translates of $P$ by the action of the symmetric group $\mathfrak{S}_k$ is dense in $\rr_{>0}^k$; by symmetry and \ref{lem:lipschitz_on_open} it suffices to show that $g$ is Lipschitz continuous on $P$. Writing $y_1=x_1$, $y_i=x_i - x_{i-1}$ for $i=2,\ldots,k$ we find that $x_i = \sum_{j=1}^i y_j$ and that $P$ is parametrized by $y_1>0, \dots, y_k >0$. Note that
\begin{displaymath}
  \sum_{i=1}^k x_i \tmatrix_i = \sum_{i=1}^k y_i \sum_{j=i}^k \tmatrix_j
\quad\text{and}\quad
\sum_{i=1}^k x_i \tmatrix_i\tvector_i = \sum_{i=1}^k y_i \sum_{j=i}^k
\tmatrix_j\tvector_j.
\end{displaymath}
By \cite[Lemma 3.5]{bghdj_sing} we have, writing $\tilde{\tmatrix}_i =
\sum_{j=i}^k \tmatrix_j$, that the \emph{flag condition}
$\mathrm{Ker}\, \tilde{\tmatrix}_i \subseteq \mathrm{Ker}\,
\tilde{\tmatrix}_{i+1}$ holds for $i=1,\ldots,k-1$. Moreover we have $\Im(\tilde \tmatrix_{i})=\sum_{j=i}^{k} \Im( \tmatrix_j )$. 
It follows that
there exist vectors $\tilde{\tvector}_i \in \bb{R}^r$ such that
\[ \sum_{j=i}^k \tmatrix_j\tvector_j = \tilde{\tmatrix}_i \tilde{\tvector}_i \, . \]

Replacing $\tmatrix_i$ by $\tilde{\tmatrix}_i$, $x_i$ by $y_i$ and $\tvector_i$
by $\tilde{\tvector}_i$ we are reduced to proving the Lipschitz continuity of $g$ on $\rr_{> 0}^k$ under the extra hypothesis that the matrices $\tmatrix
_{1},\dots,\tmatrix_{k}$ satisfy the above introduced flag condition. We do this by showing that the partial derivatives of the smooth
function $g$ are bounded on $\rr_{> 0}^k$. To shorten the notation we
set $z = \sum_{i=1}^k x_{i}A_i c_i $ and $\Omega = \sum_{i=1}^k
x_{i}A_i $. Then $g = z^t\Omega^{-1} z$ and a small computation shows
\begin{equation}
\frac{\partial g}{\partial x_i} = 2 (A_i c_i)^t \Omega^{-1} z - z^t
\Omega^{-1} A_i \Omega^{-1} z \, .  
\end{equation}
We see that it suffices to show that $\Omega^{-1}z$ is bounded on $\rr_{> 0}^k$. Write $r_j = \mathrm{rk}(A_j)$ for $j=1,\ldots,k$. We have
$ r =r_1 \ge \cdots \ge r_k \ge 1$ by the flag condition. Now \cite[Lemma 3.6]{bghdj_sing} states the existence of a constant $c_1$ such that for all integers $1 \le \alpha, \beta \le r$ we have 
\begin{equation}
\left(\Omega^{-1}\right)_{\alpha, \beta} \le \frac{c_1}{\sum_{j : r_j \ge \min(\alpha, \beta)} x_j} \le \frac{c_1}{x_1} \, .
\end{equation}
Hence there exists a constant $c_2$ such that 
\begin{equation}
\left(\Omega^{-1}\right)_{\alpha, \beta} z_\beta \le \frac{c_2 \sum_{j : r_j \ge \beta} x_j}{\sum_{j : r_j \ge \min(\alpha, \beta)} x_j} \le c_2 \, , 
\end{equation}
showing that $\Omega^{-1}z$ is bounded on $\rr_{> 0}^k$. 
\end{proof}

\begin{proof}[Proof of \ref{thm:biext_toroidal}]
Let $D(\epsilon)$ denote the open disk in $\cc$ with radius $\epsilon$.
Write $n=\dim X$. Choose a point  $p \in D$. By \cite[Theorem~81]{bp}
there exist a small positive number $\epsilon$, an open neighborhood
$V$ of $p$ in $X$, a coordinate chart  $V \isom D(\epsilon)^n$
 with $p
\mapsto 0$, and a generating section $s$ of $\ca B_\lambda$ over $U
\cap V$. Assume that $D \cap V$ is given by the equation $z_1 \cdots
z_k = 0$. Our task is to show that the psh function $-\log h_\lambda(s)$ can be written
on $U \cap V$ as a sum  
\begin{equation} \label{eqn:expansion_metric}
- \log h_\lambda(s) = \gamma + g(-\log|z_1|,\ldots,-\log|z_k|) \, , 
\end{equation}
 where $\gamma $ is locally bounded on $V$ and $g$ is a
  convex Lipschitz continuous function defined on some quadrant
 $\{x_{j}\ge M_{j} \, | \, j=1,\dots,k\}\sub \bb R^{k}$.

For this we invoke  \cite[Theorem~1.1]{bghdj_sing} and its proof in \cite[Section~4]{bghdj_sing}; these give us the expansion \ref{eqn:expansion_metric} for $-\log h_\lambda(s)$ with $\gamma $ locally bounded on $V$ and with $g$ convex on $\rr_{>0}^k$ and given as
\begin{equation}
g = \left(\sum_{i=1}^k x_{i}A_i c_i \right)^t \cdot \left(\sum_{i=1}^k x_{i}A_i 
\right)^{-1} \cdot \left(\sum_{i=1}^k x_{i}A_i c_i\right)
\end{equation}
up to a linear form in the $x_i$.
Here $A_1,\ldots,A_k$ are positive semi-definite $r \times r$-matrices
for some $r \in \zz_{>0}$ such that for all $x_1,\ldots,x_k \in
\rr_{>0}$ we have $\sum_{i=1}^k A_i x_i$ of full rank and where
$c_1,\ldots, c_k \in \rr^r$. Then by \ref{lem:is_Lipschitz_cts} the
function $g$ is Lipschitz continuous on any closed quadrant contained
in $\rr_{>0}^k$.

\end{proof}

Note that  the psh metric $\overline{h}_{\lambda}$ has its unbounded locus supported on the boundary divisor $D=X \setminus U$. Let $s$ be a non-zero rational section of the line bundle $\ca B_\lambda  $. 
We have a natural associated $\rr$-b-divisor $D(  \overline{\ca B}_\lambda, \overline{h}_{\lambda},s )$ on the category $R(X)$. We note that $D(  \overline{\ca B}_\lambda, \overline{h}_{\lambda},s )$ is actually a $\qq$-b-divisor, since on each model in $R(X)$ its incarnation is given by the Lear extension of $\ca B_\lambda  $, which is a $\bb Q$-line bundle. By~\ref{prop:3}, the b-divisor $D(  \overline{\ca B}_\lambda, \overline{h}_{\lambda},s )$ is independent of the choice of the chosen compactification $X$ of $U$.

Let $\vol_{\ca J}( \overline{\ca B}_\lambda, \overline{h}_{\lambda})$ denote the multiplier ideal volume of the psh line bundle $( \overline{\ca B}_\lambda, \overline{h}_{\lambda})$ on $X$. By combining  \ref{thm:biext_toroidal}  with the Hilbert--Samuel formula in \ref{th:asym} and the Chern-Weil  formula in  \ref{thm:b-div-chern-weil} we 
find the following result.
\begin{theorem} \label{thm:CW_for_biext}
Let $n = \dim U$. The equalities
\[ \vol_{\ca J}( \overline{\ca B}_\lambda, \overline{h}_{\lambda})  =\left(  D(  \overline{\ca B}_\lambda, \overline{h}_{\lambda},s ) \right)^n   
 =  \int_{X } \langle c_1( \overline{\ca B}_{\lambda}, \overline{h}_{\lambda} )^n \rangle = \int_{U} c_1(\ca B_\lambda, h_\lambda)^n \]
hold in $\rr_{\ge 0}$.
\end{theorem}

\subsection{The line bundle of Siegel--Jacobi forms and its invariant metric} \label{sec:Siegel-Jacobi-forms}

The aim of this final section is to discuss a variant of \ref{thm:CW_for_biext} in the context of Siegel--Jacobi forms. Let $g \in \zz_{\ge 1}$ and $N \in \zz_{\ge 3}$ and let $A_{g,N}$ denote the fine moduli space of principally polarized complex abelian varieties of dimension~$g$ and level~$N$. This is a smooth quasi-projective complex algebraic variety of dimension $g(g+1)/2$. Let  $\pi \colon U_{g,N} \to A_{g,N}$ be the universal family of abelian varieties. Following the constructions in \ref{sec:PBTor}, the tautological polarization on $U_{g,N}$ gives rise to a canonical biextension line bundle $\ca B$ on $U_{g,N}$ equipped with a smooth  semipositive hermitian metric $h$.

We write $\ca H$ for the tautological polarized variation of pure Hodge structure $ R^1 \pi_*\zz_{U_{g,N}}(1)$ of weight~$-1$ on $A_{g,N}$. We denote by $\ca F = \ca F^0(\ca H_{\bb C} \otimes \ca O_S)$ the associated Hodge bundle on $A_{g,N}$, and write $\M = \bigwedge^g \ca F$ for its determinant. The tautological polarization of $\ca H$ induces a smooth hermitian metric $h^{\mathrm{inv}}$ on $\M$ called the \emph{invariant metric} or \emph{Hodge metric}. By \cite[Lemmas~7.18 and 7.19]{schmid}, the metric $h^{\mathrm{inv}}$ is semi-positive. 
  
\begin{definition} 
Let $k, m \in \zz_{\geq 0}$. The line bundle 
\[ \L_{k,m} = \pi^*\M^{\otimes k} \otimes \ca B^{\otimes m} \]
on $U_{g,N}$ is called the \emph{line bundle of Siegel--Jacobi forms of weight~$k$ and index~$m$}. 
\end{definition}
Note that each line bundle $\L_{k,m}$  has a natural smooth hermitian metric $h_{k,m}$ obtained by taking appropriate  tensor product combinations of the metrics $h^{\mathrm{inv}}$ on $\M$ and $h$ on  $\ca B$. 
As the metrics $h^{\mathrm{inv}}$ and $h$ are semipositive, it is clear that for each $k, m \in \zz_{\geq 0}$ the smooth hermitian metric $h_{k,m}$ is a semipositive metric on $\L_{k,m}$. 

The work done in \cite[Chapter~VI]{fc} and \cite[Section~3]{hi} allows to choose:
\begin{itemize}
\item a projective smooth toroidal compactification $\overline{A}_{g,N}$ of $A_{g,N}$, 
\item an extension $\overline{\M}$ of $\M$ over $\overline{A}_{g,N}$ as a line bundle, 
\item a projective smooth toroidal compactification $\overline{U}_{g,N}$ of $U_{g,N}$, and 
\item a map $\overline{\pi} \colon \overline{U}_{g,N} \to \overline{A}_{g,N}$ extending the projection map $\pi \colon U_{g,N} \to A_{g,N}$, 
\end{itemize}
such that 
\begin{itemize}
\item the metric $h^{\mathrm{inv}}$ extends as good psh metric $\overline{h}^{\mathrm{inv}}$ over $\overline{\M}$, and 
\item the pullback of the variation $\ca H$ from $A_{g,N}$ to $U_{g,N}$ has unipotent monodromy around each local branch of the normal crossings boundary divisor $D = \overline{U}_{g,N} \setminus U_{g,N}$. 
\end{itemize}

As the psh metric $\overline{h}^{\mathrm{inv}}$ is good, it has in particular almost asymptotically algebraic singularities, see  \ref{ex:good_implies_aaas}. 
Further it follows  from \ref{bp_general} that the biextension line bundle $\ca B$ on $U_{g,N}$ has a Lear extension $\overline{\ca B}$  over $\overline{U}_{g,N}$, and from \ref{thm:biext_toroidal} that the metric $h$ has a natural extension $\overline{h}$ over the $\qq$-line bundle $\overline{\ca B}$ as a psh metric with toroidal, and in particular  almost asymptotically algebraic, singularities. 

We define
\[ \overline{\L}_{k,m} = \overline{\pi}^* \overline{\M}^{\otimes k} \otimes \overline{\ca B}^{\otimes m} \, , \]
a $\qq$-line bundle on $\overline{U}_{g,N}$ extending the line bundle $\L_{k,m}$,
and denote by $\overline{h}_{k,m}$ the metric on $\overline{\L}_{k,m}$ induced from $\overline{h}^{\mathrm{inv}}$ on $\overline{\M}$ and $\overline{h}$ on $\overline{\ca B}$ by taking the appropriate tensor product combinations.

As both $\overline{h}^{\mathrm{inv}}$ and $\overline{h}$ have almost asymptotically algebraic singularities, we conclude by \ref{tensor_toroidal} that the natural metric $\overline{h}_{k,m}$ on the $\qq$-line bundle $\overline{\L}_{k,m}$  is a psh metric with almost asymptotically algebraic singularities.  
Further we note that  the unbounded locus of the psh metric $\overline{h}_{k,m}$ has support contained in the boundary divisor $D = \overline{U}_{g,N} \setminus U_{g,N}$. 

Let $g \in \zz_{\ge 1}$, $N \in \zz_{\ge 3}$, $k, m \in \zz_{\ge 0}$. Let $s$ be a non-zero rational section of the line bundle $\L_{k,m}$ of Siegel--Jacobi forms on $U_{g,N}$ of weight~$k$ and index~$m$. Let $D(  \overline{\L}_{k,m}, \overline{h}_{k,m},s )$ be the $\rr$-b-divisor associated to the Lear extension $ \overline{\L}_{k,m}$ over $\overline{U}_{g,N}$, its canonical psh metric $\overline{h}_{k,m}$ and the rational section $s$. As follows from \ref{prop:3}, the b-divisor $D(  \overline{\L}_{k,m}, \overline{h}_{k,m},s )$ is independent of the choice of compactification $\overline{U}_{g,N}$. 

Applying \ref{th:asym} and \ref{thm:b-div-chern-weil} we obtain the following result. 
 
\begin{theorem}\label{thm:siegel-jacobi-CW} Let $n = \dim U_{g,N} = g + g(g+1)/2$. Let $\vol_{\ca J}( \overline{\L}_{k,m}, \overline{h}_{k,m} )$ denote the multiplier ideal volume of the psh line bundle $( \overline{\L}_{k,m}, \overline{h}_{k,m} )$. 
 Then  the equalities
\[  \begin{split} \vol_{\ca J}( \overline{\L}_{k,m}, \overline{h}_{k,m} )  & = \left(  D(  \overline{\L}_{k,m}, \overline{h}_{k,m},s ) \right)^n  \\
& =  \int_{\overline{U}_{g,N} } \langle c_1( \overline{\L}_{k,m}, \overline{h}_{k,m} )^n \rangle \\
& = \int_{U_{g,N}} c_1(\ca L_{k,m},h_{k,m})^n \end{split} \]
hold in $\rr_{\ge 0}$.
\end{theorem}
In the case of the modular curve $Y(N)=A_{1,N}$ and for $ k=m=4$ this reproduces the main result of \cite{bkk}. We observe that in \cite{bkk}, volume, degree and integral were each calculated independently, and the outcomes were seen to be equal by inspection. Here we see a more intrinsic approach.

\begin{remark} The above Chern--Weil type result can be extended, with proofs carrying over mutatis mutandis, to the more general setting of a finite self-product of the universal abelian scheme over any PEL type Shimura variety. The necessary smooth projective toroidal compactifications are provided in this setting by the work of Lan \cite{lan}.
\end{remark}

\appendix

\section{On the non-continuity of the volume function}
\label{ex:non-cont}
\newcommand{\D}{\mathbb{D}}
\newcommand{\bdiv}{\on{b\text{-}div}}
\newcommand{\Q}{\mathbb{Q}}
\newcommand{\caD}{\mathcal{D}}
\newcommand{\caR}{\mathcal{R}}
\newcommand{\N}{\mathbb{N}}

In \ref{thm:3} the multiplier ideal volume of
a line bundle with a psh metric is defined. This volume function has nice
properties for metrics with almost asymptotically algebraic singularities. For instance, it
agrees with the non-pluripolar volume (\ref{th:asym}) and it is 
continuous with respect to  good algebraic approximations (see
\ref{rem:mult-vol-cont} for the precise statement).

Another approach to define the volume of a line bundle with a psh
  metric, assuming it has Zariski unbounded locus, is to measure the
  abundance of global sections for multiples of the associated
  b-divisor. In fact, if the b-divisor is Cartier, this is the
  classical definition of the volume of a divisor. In this appendix we
  give the definition of the volume of a b-divisor that is the
  analogue of the volume function of a classical divisor. We can then
  ask whether this volume function is continuous with respect to the
  weak topology. The results of \cite[\S 5]{BoteroBurgos} show that
  this question has an affirmative answer in the toroidal
  case. Nevertheless, we give an example showing that, 
  in general, this volume function is not continuous for approximable
  nef b-divisors, even in the big 
case. We don't know if the constructed example can be realized as the
b-divisor associated to a psh metric. 
 
 We start with some definitions.
 Let $X$ be a smooth projective complex variety with function field $F$. 
Let $\D = (D_\pi)_{\pi \in R(X)}$ be a Weil $\mathbb R$-b-divisor on
$X$ as in \ref{sec:b-div}. We define 
\begin{align*}
  \ca L(\D)&=\{0\not = f \in F\mid \forall \pi \in R(X) \, \colon \, D_{\pi }+\on{div}(f) \ge 0  
               \}\cup \{0\} \, . 
\end{align*}

The \emph{volume} $\vol(\D)$ of the b-divisor $\D$ is defined via the formula
\begin{equation}\label{eq:vol-b}
\vol(\D) \coloneqq \limsup_{\ell} \frac{\dim \ca L(\ell\D)}{\ell^{d}/d!}.
\end{equation} 
We say $\D$ is \emph{big} if  $\vol(\D)>0$. 

The example is constructed in two steps. The first one is a preparatory step.

\begin{figure}[H]
\begin{center}
\begin{tikzpicture}[scale = 1.5]
\draw (-5.6,1.4) node{$X_0$};
\draw (-5.2,0) to  [bend left = 25] (-4,0) node[below left]{$A$};
\draw (-5,-0.2) to  [bend right = 25] (-5,1) node[below left]{$B$};
\draw (-4.7,0.2) node[font = \scriptsize]{$p_0$};
\draw (-4.88,0.1) node[font=\scriptsize]{$\bullet$};
\draw (-3,1.4) node{$X_1$};
\draw (-2.7,0) to  [bend left = 25] (-1.5,0) node[below left]{$\widehat{A}$};
\draw (-2.5,-0.2)  to  [bend right = 25] (-2.5,1);
\draw (-2.45, 0.4) node[right, font = \small]{$E_1$};
\draw (-2.7,0.7) to  [bend right = 25] (-1.5,1.2) node[left]{$\widehat{B}$};
\draw (-2.25, 0.8) node[font = \scriptsize]{$p_1$};
\draw (-2.4,0.71) node[font=\scriptsize]{$\bullet$};
\draw (-0.4,1.4) node{$X_2$};
\draw (-0.1,-0.5) to  [bend left = 25] (1.1,-0.5) node[below left]{$\widehat{A}$};
\draw (0.1,-0.2)  to  [bend right = 25] (0.1,1);
\draw (0.1,-0.7)  to  [bend right = 30] (0.1,0.2);
\draw (0.15, 0.4) node[right, font = \small]{$E_2$};
\draw (-0.1,0.7) to  [bend right = 25] (1.1,1.2) node[left]{$\widehat{B}$};
\draw (0.15, -0.15) node[right, font = \small]{$\widehat{E_1}$};
\draw (0.35, 0.83) node[font = \scriptsize]{$p_2$};
\draw (0.2,0.71) node[font=\scriptsize]{$\bullet$};
\draw (2.2, 0) node{$\dotsc$};
\draw (3.1,1.4) node{$X_b$};
\draw (3.4,0.7) to  [bend right = 25] (4.5,1.2) node[left]{$\widehat{B}$};
\draw (3.85, 0.85) node[font = \scriptsize]{$p_b$};
\draw (3.75,0.72) node[font=\scriptsize]{$\bullet$};
\draw (3.6, 0.4) node[right, font = \scriptsize]{$E_b$};
\draw (3.6,0.2)  to  [bend right = 30] (3.6,1.2);
\draw (3.6,0.1) node{$\vdots$};
\draw (3.6,-0.2)  to  [bend left = 30] (3.6,-0.6);
\draw (3.6, -0.4) to [bend left = 30] (3.6,-0.8);
\draw (3.5,-0.75) to [bend left = 25] (4.5,-1) node[above]{$\widehat{A}$};
\draw (3.75,-0.6) node[font = \tiny]{$\widehat{E_1}$};
\draw (3.75,-0.3) node[font = \tiny]{$\widehat{E_2}$};
\end{tikzpicture}
\end{center}
\caption{The surfaces $X_b$}
\label{fig1}
\end{figure}
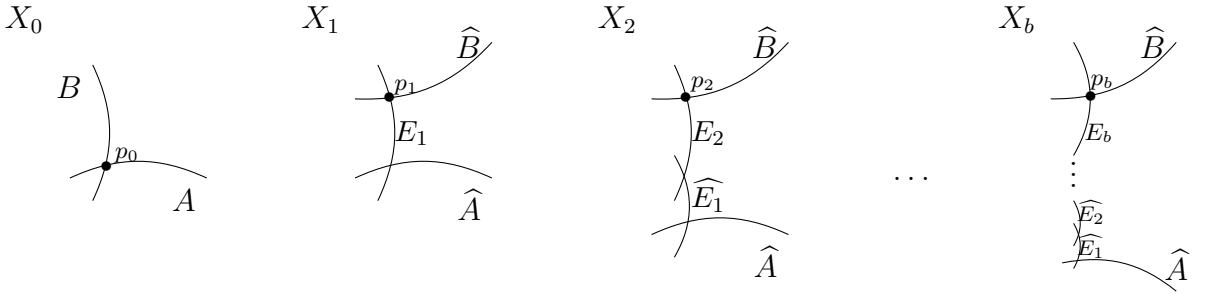

 \textbf{Step 1:}\\
 Let $X$ be a smooth projective surface, and $A$, $B$ divisors on $X$
 meeting transversely at a point $p$ and let $b \geq 1$ be an
 integer. Setting $X_0 = X$ and $p_{0} = p$ we make the following
 recursive definition, which is illustrated in  \ref{fig1}:  
 \begin{enumerate}
 \item
   For $0\le i < b$,
 let $X_{i+1}$ be the blow up of $X_i$ at $p_i$, and $E_{i+1}$ the
 exceptional locus of the blowup;
 \item For any $i \in \{1, \dotsc, b\}$ and any divisor $F$ on $X_j$
   for $0 \leq j \leq i$, we write $\widehat{F}$ for the strict
   transform of $F$ on $X_{i}$ (note that $F = \widehat{F}$ if $i
   =j$).
  \item Let $p_{i+1}$ be the unique point of intersection of
    $\widehat{B}$ with $E_{i+1}$ on $X_{i+1}$. 
 \end{enumerate}
 
 We easily compute some intersection numbers on the surface $X_b$ for any $b>0$:
 \begin{itemize}
 \item $\widehat{E_i} \cdot \widehat{E_j} = 1 $ if $\abs{i-j} = 1$ and 0 if $\abs{i-j} \ge 2$; 
 \item $\left(E_b\right)^2 = -1$;
 \item For all $i < b$ we have $\left(\widehat{E_i}\right)^2 = -2$;
 \item $\left(\widehat{A}\right)^2 = A^2 - 1$ and $\left(\widehat{B}\right)^2 = B^2 - b$. 
 \end{itemize}
 
 Now let $D$ be a nef divisor on $X$ whose support does not contain $p$.
 
 On $X_b$ we define the $\Q$-divisor
   \begin{displaymath}
     D_b = \pi^*D - \sum_{i =1}^{b} \frac{i}{b}\widehat{E_i},
   \end{displaymath}
 where $\pi \colon X_b \to X$ denotes the corresponding sequence of blow-ups.

 \begin{lemma}\label{lem:non-cont}
 For each integer $b \geq 1$, the following equalities hold true:
 \begin{enumerate}
 \item $\left(D_b\right)^2 = D^2 - \frac{1}{b}$,
 \item $D_b \cdot \widehat{A} = D \cdot A - \frac{1}{b}$,
 \item $D_b \cdot \widehat{E_i} = 0 \text{ for }i < b$, 
 \item $D_b \cdot E_b = \frac{1}{b}$, 
 \item  $D_b\cdot \widehat{B} = D \cdot B - 1$.
 \end{enumerate}

 \end{lemma}
 \begin{proof}
For $1.$ first observe that $\pi^*D \cdot \widehat{E_i} = 0$ for all $i$ since $D$ does not pass through $p$. We then calculate
\begin{equation*}
\begin{split}
\left(D_b\right)^2 & = D^2 - 2\sum_{i=1}^{b-1} \frac{i^2}{b^2} - 1 + 2\sum_{i=1}^{b-1} \frac{i(i+1)}{b^2}\\
& = D^2 - 1 + 2\sum_{i=1}^{b-1} \frac{i}{b^2} \\ & = D^2 - \frac{1}{b}. 
\end{split}
\end{equation*}
Similarly, $\widehat{A}\cdot \widehat{E_i} = 0$ for $1 < i \le b$ and $\widehat{B} \cdot \widehat{E_i} = 0$ for $1 \le i < b$, so we compute
\begin{equation*}
D_b \cdot \widehat{A} = D \cdot A - \frac{1}{b}\widehat{E_1} \cdot \widehat{A} = D \cdot A - \frac{1}{b}, 
\end{equation*}
\begin{equation*}
D_b \cdot \widehat{B} = D \cdot B - E_b \cdot B = D\cdot B -1, 
\end{equation*}
which gives $2.$ and $5.$ On the other hand, we have 

\begin{equation*}
D_b \cdot \widehat{E_1} = \frac{-1}{b} \left(\widehat{E_1}\right)^2 - \frac{2}{b} \widehat{E_2} \cdot \widehat{E_1} = \frac{2}{b} - \frac{2}{b} = 0, 
\end{equation*}
and for $2 \le i \le b-1$ we get
\begin{equation*}
D_b \cdot \widehat{E_i} = - \frac{i-1}{b} \widehat{E_{i-1}} \cdot \widehat{E_i} - \frac{i}{b} \left(\widehat{E_i}\right)^2 - \frac{i+1}{b} \widehat{E_{i+1}} \cdot \widehat{E_i} = -\frac{i-1}{b} + \frac{2i}{b} - \frac{i+1}{b} = 0,
\end{equation*}
which gives $3.$ Finally, we compute
\begin{equation*}
D_b \cdot E_b = -\frac{b-1}{b} \widehat{E_{b-1}} \cdot E_b - \left(E_b\right)^2 = -\frac{b-1}{b} + 1 = \frac{1}{b},
\end{equation*}
giving $4.$
 \end{proof}
 
 The second step is the main construction.

 \textbf{Step 2:}\\
 Now we let $X = \bb P^2$, and we choose a line $L$ on $X$ and
 $(p_k)_{k \ge 0}$ an ordered countable set of distinct points on
 $L$. Let $H$ be a line not passing through any $p_k$, and let $D=
 2H$, a divisor on $X$. For each $k$ we choose $B_k$ a line through
 $p_k$ distinct from $L$.
  
 Set $X'_{0} = X$ and $D'_0 = D$. For any $k$ we let $X'_{k+1}$ and $D'_{k+1}$ be the result of applying Construction 1 to $X'_{k}$, $A = \widehat{L}$, $B= B_k$, $p = p_k$ and $D'_k$ with $b = 2^k$. 
By \ref{lem:non-cont}, for any $k$, we find that
\[
\left(D'_{k}\right)^2 = D^2 - \frac{1}{2} - \frac{1}{4} -\cdots - \frac{1}{2^k} = 3 + \frac{1}{2^k}. 
\]
 
 \begin{lemma}
 $D'_k$ is nef. 
 \end{lemma}
 \begin{proof}
 First we compute (using \ref{lem:non-cont}) that 
 \[
 D'_k \cdot \widehat{L} = D \cdot A - \sum_{i=1}^k \frac{1}{2^i} = 1 + \frac{1}{2^k} \geq 0, 
 \]
 \[
 D'_k \cdot \widehat{B_i} = D \cdot B_i - 1 = 1 \geq 0, 
 \]
 and $D'_k \cdot E \ge 0$ for any irreducible exceptional divisor $E$. Now let $C \subset X'_k$ be any irreducible curve distinct from $\widehat{L}$, $\widehat{B_i}$, and exceptional curves. Denote by $E_{i,j}$ the $i$'th exceptional divisor on $X'_j$. Then 
 \begin{equation*}
 \sum_{j=1}^k \sum_{i=1}^{2^j} \widehat{E_{ij}} \cdot C \le L \cdot \pi_*C,
 \end{equation*}
 where $\pi \colon X_k' \to X$ denotes the corresponding sequence of blow-ups. Indeed, 
 \begin{equation*}
 \begin{split}
 L \cdot \pi_*C & = \pi^*L \cdot C\\
 & = (\widehat{L} + \sum_{i,j} \widehat{E_{i,j}}) \cdot C\\
 & = \widehat{L} \cdot C + \sum_{i,j} \widehat{E_{ij}} \cdot C, 
 \end{split}
 \end{equation*}
 each summand of which is non-negative. We then compute
 \begin{align*}
 D'_k \cdot C & = \pi^*D \cdot C - \sum_{i,j} \frac{i}{2^j} \widehat{E_{i,j}} \cdot C\\
 & \ge \pi^*D \cdot C - \sum_{i,j}  \widehat{E_{i,j}} \cdot C\\
 & \ge \pi^*D \cdot C - \pi^*L \cdot C\\
 &= L \cdot \pi_*C \geq  0, 
 \end{align*}
 where the first equality on the last line holds because $D \sim 2L$. 
 \end{proof}

 Now, the sequence $\{D_k'\}_{k \in \N}$ converges in the weak
 topology and defines an approximable nef b-divisor on $R(X)$ which we denote by $\D$. 
 
Given $k \ge 0$, let $f$ be a rational function on $X$ such that 
 \[
 k\D + \operatorname{div} f \ge 0.
 \]
Since $f$ must cancel infinitely many exceptional  divisors over $L$ with multiplicity $k$ we see that $\operatorname{ord}_L f \ge k$, and in fact this condition is equivalent to canceling all the negative multiples of exceptional divisors. Since the sum of these exceptional parts is by definition given by $2H - \D$ we see that 
\[
\vol(\D) = \vol(2H - L) = 1/2. 
\]
This implies in particular that $\D$ is big. On the other hand, since the $D'_k$'s are nef, we have 
  \[
  2\vol(D_k') = \left(D_k'\right)^2 = 3 + \frac{1}{2^k}. 
  \]
 Hence $\lim_k  \vol(D_k') = 3/2$.  
In conclusion, we see that the volume function is not continuous in
the space of approximable big and nef b-divisors. In consequence, for
an approximable nef b-divisor, the volume and the degree 
do not agree necessarily.
 
 \begin{remark}\label{rem:app}
 \begin{enumerate}
   \item In \cite[Theorem 5.13]{BoteroBurgos} it is shown
   that, for toroidal nef b-divisors, the degree and the volume
   agree. In that paper it was also announced that the volume was not
   necessarily continuous in the non-toroidal case. This is the
   promised example.
\item To a psh metric with almost asymptotically algebraic singularities, we can associate an
   approximable nef b-divisor (\ref{thm:psh_implies_approximable}) and
   the multiplier ideal volume of the metric agrees with the degree of
   this b-divisor (\ref{cor:hs}).   Then, \cite[Theorem
   5.13]{BoteroBurgos} and \ref{cor:hs} imply that, for toroidal psh metrics, the multiplier
   ideal volume of the psh metric agrees with the volume of the
   associated b-divisor.   
 \item  It would be interesting to know whether the equality between the
  multiplier ideal volume of a psh metric and the volume of the
  associated b-divisor continues to be true in the general case of almost asymptotically
  algebraic psh metrics. Note that if one can show that the above
  example can be realized as the b-divisor of an almost asymptotically
  algebraic psh metric, we obtain a negative answer to this
  question. 
  \end{enumerate}
  \end{remark}

\end{document}